\documentclass[11pt]{article}
\usepackage[T1]{fontenc}
\usepackage{lmodern,amsmath,amssymb,amsthm,mathtools,listings}
\usepackage[margin=1in]{geometry}
\usepackage[colorlinks=true,linkcolor=blue,citecolor=blue,urlcolor=blue]{hyperref}

\newtheorem{theorem}{Theorem}[section]
\newtheorem{lemma}[theorem]{Lemma}
\newtheorem{proposition}[theorem]{Proposition}
\newtheorem{corollary}[theorem]{Corollary}
\newtheorem{conjecture}[theorem]{Conjecture}
\DeclareMathOperator{\Li}{Li}
\DeclareMathOperator{\CT}{CT}
\DeclareMathOperator{\Rea}{Re}
\DeclareMathOperator{\Ima}{Im}
\title{Real zeros of Durfee polynomials: Proof of the Canfield--Corteel--Savage conjecture}
\author{Jing Huang$^1$ and Guanghui Zhang$^2$\footnote{E-mail addresses: {\it
jhuangmath@gzhu.edu.cn (J. Huang), zghui@squ.edu.cn (G. Zhang).
}}
}

\date{\small
$1.$ School of Mathematics and Information Science,
Guangzhou University, Guangzhou 510006, China\\
$2.$ School of Mathematics and Physics, Suqian University, Suqian 223800, China
}

\begin{document}
\maketitle

\begin{abstract}
We prove the conjecture of Canfield, Corteel, and Savage that for
every positive integer $n$, the Durfee polynomial
$D_n(x)=\sum_{d=1}^{\lfloor\sqrt n\rfloor}p(n,d)x^d$ counting partitions
of $n$ by the side length of their Durfee square has only real zeros.
More precisely, $D_n$ has a simple zero at the origin and exactly
$\lfloor\sqrt n\rfloor-1$ simple negative zeros. Consequently, for every
$n\ge1$, the sequence $p(n,1),\ldots,p(n,\lfloor\sqrt n\rfloor)$ is strictly log-concave and unimodal.
As a further consequence,  the Durfee square side length of a uniformly chosen
partition of $n$ has the same distribution as $1$ plus a sum of independent Bernoulli
random variables.


\end{abstract}

\section{Introduction}

A partition of a nonnegative integer $n$ is a finite weakly decreasing
sequence $\lambda=(\lambda_1,\ldots,\lambda_\ell)$ of positive integers
with sum $n$; the empty partition is allowed when $n=0$.
We write $\lambda\vdash n$.
In the Ferrers diagram of $\lambda$, the $i$th row has $\lambda_i$ cells.
Its Durfee square is the largest square
in the upper left corner, with side length
$d(\lambda)=\max\{1\le j\le\ell:\lambda_j\ge j\}$ for a nonempty
partition, with $d(\varnothing)=0$; see Andrews~\cite{Andrews}. The Durfee polynomial is
\[
 D_n(x)=\sum_{\lambda\vdash n}x^{d(\lambda)}
       =\sum_{d=1}^{\lfloor\sqrt n\rfloor}p(n,d)x^d,
\]
where
$p(n,d)$ is the number of partitions of $n$ whose Durfee square
has side length $d$;
we set $D_0(x)=1$. The coefficient $p(n,d)$ is positive for
$1\le d\le \lfloor\sqrt n\rfloor$, and $p(n,1)=n$. Hence $D_n$ has degree
$\lfloor\sqrt n\rfloor$ and a simple zero at the
origin.

In 1998, Canfield, Corteel, and Savage introduced Durfee polynomials
for several families of partitions. For ordinary partitions, their
computations found only real zeros for $1\le n\le1000$ and strict
log-concavity of the coefficients for $1\le n\le5000$~\cite{CCS}.
Real-rootedness would imply log-concavity and unimodality by Newton's
inequalities; see Stanley~\cite{Stanley}. These observations led them
to the following conjecture, which we state for ordinary partitions.

\begin{conjecture}[Canfield--Corteel--Savage {\cite{CCS}}]
\label{conj:CCS}
For every integer $n\ge1$, the Durfee polynomial $D_n$ has only real
zeros.
\end{conjecture}

Canfield, Corteel, and Savage obtained an
asymptotic formula for $p(n,d)$ and proved that the Durfee square side
length, after centering and scaling, has a limiting normal
distribution~\cite{CCS};
in the same paper,  they also proved that, for each fixed $0<\varepsilon<1/2$, the
coefficients are strictly log-concave throughout
$\varepsilon\sqrt n\le d\le(1-\varepsilon)\sqrt n$ when $n$ is
sufficiently large. Mutafchiev proved a local
limit theorem for the side length in 2002~\cite{Mutafchiev}.
Boyer and Goh restated the conjecture in their study of the
asymptotics and zeros of partition polynomials~\cite{BG}.

In 2012, Sagan and Savage~\cite{SS} showed that Foata's bijection,
followed by the lattice-path encoding of partitions, sends the major
index of a binary word to partition size and its descent number to the side
length of the Durfee square. Together with Simion's real-rootedness theorem for multiset Eulerian polynomials~\cite{Simion}, this gives real-rootedness of the
Durfee generating polynomial over all partitions contained in a fixed rectangle.
Haglund and Visontai~\cite{HV} established an
equivalent formulation of Conjecture~\ref{conj:CCS}: the real-rootedness
of these
multiset Eulerian
polynomials is preserved even when restricting to words with a prescribed major index.

At the 2016 AIM workshop on \emph{Polyhedral geometry and partition
theory}, participants investigated a coefficient recurrence and an
approach through Eulerian polynomials. The workshop
report~\cite{AIM2016} records no proof from either approach.
Savage again listed the conjecture among three open questions in her
interview with Mansour~\cite{Mansour}, published online in
2021 and included in the 2022 volume of \emph{Enumerative
Combinatorics and Applications}. In this paper,
we prove Conjecture~\ref{conj:CCS} and show, in addition, that all
zeros are simple.

\begin{theorem}\label{thm:main}
For every integer $n\ge1$, the Durfee polynomial $D_n$ has a simple
zero at the origin and exactly $\lfloor\sqrt n\rfloor-1$ simple
negative zeros.
\end{theorem}

By Newton's inequalities~\cite{Stanley}, Theorem~\ref{thm:main} implies
that the sequence $p(n,1),\ldots,p(n,\lfloor\sqrt n\rfloor)$ is strictly log-concave and
unimodal for every $n\ge1$. This extends the log-concavity result of
Canfield, Corteel, and Savage~\cite{CCS} to every $n$ and the full range
$1\le d\le \lfloor\sqrt n\rfloor$.
Theorem~\ref{thm:main} also provides a finite-$n$ probabilistic
structure underlying the central limit theorem of Canfield, Corteel, and
Savage~\cite{CCS} and the local limit theorem of
Mutafchiev~\cite{Mutafchiev}. Write $p(n)=D_n(1)$ for the number of
partitions of $n$, and for $n\ge1$ let $\Lambda_n$ be chosen uniformly
from these partitions, so that $\mathbb P(d(\Lambda_n)=d)=p(n,d)/p(n)$.

\begin{corollary}\label{cor:bernoulli}
For every integer $n\ge1$, let $-r_{n,1},\ldots,-r_{n,\lfloor\sqrt n\rfloor-1}$ with
$r_{n,j}>0$ be the nonzero zeros of $D_n$. Then $d(\Lambda_n)$ has the
same distribution as $1+\sum_{j=1}^{\lfloor\sqrt n\rfloor-1}\xi_{n,j}$, where
$\xi_{n,1},\ldots,\xi_{n,\lfloor\sqrt n\rfloor-1}$ are independent Bernoulli random
variables with $\mathbb P(\xi_{n,j}=1)=1/(1+r_{n,j})$. In particular,
\[
 \mathbb E[d(\Lambda_n)]=1+\sum_{j=1}^{\lfloor\sqrt n\rfloor-1}\frac1{1+r_{n,j}},
 \qquad
 \operatorname{Var}(d(\Lambda_n))
 =\sum_{j=1}^{\lfloor\sqrt n\rfloor-1}\frac{r_{n,j}}{(1+r_{n,j})^2}.
\]
\end{corollary}

 Corollary~\ref{cor:bernoulli} follows immediately from the
factorization $D_n(x)/p(n)=x\prod_{j=1}^{\lfloor\sqrt n\rfloor-1}(r_{n,j}+x)/(r_{n,j}+1)$.
By the general results for sums of independent Bernoulli random
variables reviewed in Pitman~\cite[Section~2]{Pitman}, this
representation yields central and local normal limit theorems
whenever $\operatorname{Var}(d(\Lambda_n))\to\infty$.
The same representation places every mode of $p(n,d)$ within
distance~$1$ of $\mathbb E[d(\Lambda_n)]$.

To prove Theorem~\ref{thm:main}, we estimate $D_n(-r^2)$ ($r>0$)
across four analytic regimes to produce alternating signs for all
$n\ge6401$, and then apply the intermediate value theorem.
Section~\ref{sec:preliminaries} collects the notation and identities
used throughout; in particular, Euler's identity writes
$F(x,q)=\sum_{n\ge0}D_n(x)q^n$ as the constant term
$\CT_z[(-z;q)_\infty(-xq/z;q)_\infty]$, giving a double integral of an
infinite product for $D_n(-r^2)$. Section~\ref{sec:central} treats the
near-origin and intermediate regimes. Near the origin
($0<r\le16\pi^2/n$), we approximate $D_n(-r^2)$ in terms of
$\sigma=\sqrt{nr}$ by Bessel products: at $2\sigma=k\pi$
($1\le k\le8$), the first correction
$-\frac{\sigma^4}{2n^2}J_2(2\sigma)I_2(2\sigma)$ has the same sign
$(-1)^k$ as the leading term
$\frac{\sigma}{n}(J_0(2\sigma)I_1(2\sigma)-J_1(2\sigma)I_0(2\sigma))$.
In the intermediate range $16\pi^2/n\le r\le e\sqrt n$, two complex
conjugate saddle points yield an effective oscillatory bound
$|D_n(-r^2)-M_n(r)\cos\Psi_n(r)|<cM_n(r)$ ($c<1$), with amplitude
$M_n(r)>0$ and phase $\Psi_n(r)$ governed by the dilogarithm via
$\sqrt{-2\Li_2(-ir)}$, so $\operatorname{sgn}D_n(-r^2)=(-1)^k$ whenever
$\Psi_n(r)=k\pi$. For $r\le1$, the logarithm of the product expands
directly. To cross the barrier $r=1$, we unfold the angular integral
onto $\mathbb R$ via Jacobi's theta and Dedekind's $\eta$
transformations, isolating two conjugate branches and exponentially
suppressing the rest; the minor arcs are bounded via Fourier estimates
on finitely many intervals.

Larger radii are treated in Section~\ref{sec:outer} through the
two-colour partition identity
\[
p(n,d)=[q^{n-d^2}](q;q)_d^{-2}
=[q^{n-d^2}](q;q)_\infty^{-2}\prod_{\nu>d}(1-q^\nu)^2.
\]
On $\log(n-1/12)\le2\log r\le4(n-1/12)^{1/6}$, we expand the part-size
restriction $\prod_{\nu>d}(1-q^\nu)^2$ and apply the parameter-$2$
Rademacher-type formula of Iskander, Jain, and
Talvola~\cite[Theorem~1.1]{IJT} to replace the unrestricted
two-colour partition numbers $p_2(t)$ by a compactly supported Bessel
kernel with an explicit bilateral Laplace transform. Poisson summation
then isolates the dominant frequencies $\pm\pi$, while the restriction
corrections are controlled by exponential weighting and Euler's
identity. Beyond $2\log r=4(n-1/12)^{1/6}$, uniform bounds on $p_2(t)$
yield monomial dominance
$p(n,d)r^{2d}>\sum_{\ell\ne d}p(n,\ell)r^{2\ell}$ at a suitable radius
$r$ for each remaining $d\le \lfloor\sqrt n\rfloor-1$. Matching the continuous phase
$\Psi_n(r)$ to the discrete sign counts at both transition endpoints
yields all $\lfloor\sqrt n\rfloor-1$ simple negative zeros for $n\ge6401$. For
$1\le n\le6400$, we certify the required sign changes in Python using
exact integer arithmetic at dyadic rational points; see the end of
Section~\ref{sec:outer}. 
In addition, Theorem~\ref{thm:main}  has been
formally verified in Lean~4 using only Lean's standard axioms
(\texttt{propext}, \texttt{Classical.choice}, and \texttt{Quot.sound});
the formalization and accompanying verification script are available at
\url{https://github.com/tfyf103/durfee-ccs-lean}.

\section{Preliminaries}\label{sec:preliminaries}
We fix the notation and record the generating functions and
analytic identities used in the proof of Theorem~\ref{thm:main}.
Given a formal power series $H$ in an indeterminate $q$, let $[q^j]H$
denote the coefficient of $q^j$ (with $[q^j]H=0$ when $j<0$), and
similarly for other indeterminates.
For a Laurent series in the variable $z$, let $\CT_z$ denote its constant term.
Unless otherwise specified, $\log$ and $\arg$ have their
principal values on $\mathbb C\setminus(-\infty,0]$, with
$-\pi<\arg z<\pi$, and
$\sqrt z=\exp(\tfrac12\log z)$.
For an integer $d\ge0$, $|q|<1$, and $z\in\mathbb C$, write
\[
 (q;q)_d=\prod_{j=1}^d(1-q^j),\qquad
 (z;q)_\infty=\prod_{j\ge0}(1-zq^j),
\]
where an empty product is $1$.

The Durfee-square decomposition~\cite[Section~2]{CCS}
gives a bijection between partitions of $n$ whose Durfee square has side
length $d$ and pairs of partitions of total size $n-d^2$, each with at
most $d$ parts: remove the square and transpose the partition below it.
Thus, for $|q|<1$,
\begin{equation}\label{eq:F}
 \begin{aligned}
 F(x,q)
 &=\sum_{n\ge0}D_n(x)q^n
   =\sum_{d\ge0}\frac{q^{d^2}x^d}{(q;q)_d^2},\\[2pt]
 p(n,d)
 &=[x^d]D_n(x)
   =[q^{n-d^2}](q;q)_d^{-2},
 \end{aligned}
\end{equation}
where $x\in\mathbb C$ and $n,d\in\mathbb Z_{\ge0}$;
we use $D_0(x)=1$, $p(0,0)=1$, and $p(n,0)=0$ for $n\ge1$.
By Euler's $q$-series identity~\cite[Corollary~2.2]{Andrews},
$(-z;q)_\infty=\sum_{d\ge0}q^{d(d-1)/2}z^d/(q;q)_d$
for $|q|<1$ and $z\in\mathbb C$.
Expanding the two factors in $(-z;q)_\infty(-xq/z;q)_\infty$
as absolutely convergent series on $\mathbb C\setminus\{0\}$
and comparing the constant term in $z$ with~\eqref{eq:F} gives
\[
 F(x,q)=\CT_z\bigl[(-z;q)_\infty(-xq/z;q)_\infty\bigr].
\]
For $r,\theta\in\mathbb R$ and $\Rea\zeta>0$, define
\begin{equation}\label{eq:P-product}
 \mathcal P_r(\zeta,\theta)
 =\prod_{j\ge0}
   \bigl(1+ir e^{-(j+1/2)\zeta+i\theta}\bigr)
   \bigl(1+ir e^{-(j+1/2)\zeta-i\theta}\bigr).
\end{equation}
Setting $x=-r^2$ and $q=e^{-\zeta}$ in the constant-term formula
and expressing $\CT_z$ via Cauchy's integral formula on the positively
oriented circle $|z|=r e^{-\Rea\zeta/2}$, parametrized by
$z=ir e^{-\zeta/2+i\theta}$ ($-\pi\le\theta\le\pi$), yields
\begin{equation}\label{eq:ct-integral}
 F(-r^2,e^{-\zeta})
 =\frac1{2\pi}\int_{-\pi}^{\pi}
       \mathcal P_r(\zeta,\theta)\,d\theta
 \qquad(r>0,\ \Rea\zeta>0).
\end{equation}
Likewise, for every integer $n\ge0$, $r\in\mathbb R$, and $u>0$,
applying Cauchy's coefficient formula to $F(-r^2,q)$ on the positively
oriented circle $|q|=e^{-u}$ with $q=e^{-u+it}$ ($-\pi\le t\le\pi$),
followed by the substitution $v=-t$, gives
\begin{equation}\label{eq:ct-integral2}
 D_n(-r^2)
 =\frac1{2\pi}\int_{-\pi}^{\pi}
     e^{n(u+iv)}
     F(-r^2,e^{-u-iv})\,dv.
\end{equation}

For $|r|\le1$, $\theta\in\mathbb R$, and $\Rea\zeta>0$, we write
$\log\mathcal P_r(\zeta,\theta)$ for the sum of the principal logarithms
$\log\bigl(1+ir e^{-(j+1/2)\zeta\pm i\theta}\bigr)$ over $j\ge0$ and both signs.
Since $|ir e^{-(j+1/2)\zeta\pm i\theta}|\le e^{-(j+1/2)\Rea\zeta}\le e^{-\Rea\zeta/2}<1$
and $|\log(1+w)|\le|w|/(1-|w|)$ for $|w|<1$, the $j$th terms have modulus at most
$e^{-(j+1/2)\Rea\zeta}/(1-e^{-\Rea\zeta/2})$; hence the series converges
absolutely and satisfies $\exp(\log\mathcal P_r(\zeta,\theta))=\mathcal P_r(\zeta,\theta)$.

Let $\Li_2$ denote the principal dilogarithm, defined by
$\Li_2(z)=\sum_{j\ge1}z^j/j^2$ for $|z|\le1$ ($\Li_2(1)=\pi^2/6$)
and continued analytically to $\mathbb C\setminus[1,\infty)$ via
$\Li_2(z)=-\int_0^z\frac{\log(1-t)}{t}\,dt$.
The identities
\begin{equation}\label{eq:prelim-dilog-identities}
 \begin{aligned}
 \Li_2(z)+\Li_2\!\left(\frac{z}{z-1}\right)
 &=-\frac12\log^2(1-z)
 &&\bigl(z\in\mathbb C\setminus[1,\infty)\bigr),\\
 \Li_2(z)+\Li_2(1/z)
 &=-\frac{\pi^2}{6}-\frac12\log^2(-z)
 &&\bigl(z\in\mathbb C\setminus[0,\infty)\bigr)
 \end{aligned}
\end{equation}
hold with these branches;
both identities follow by differentiating with respect to $z$
and evaluating at $z=0$ and $z=-1$, respectively.
For $r\in[-1,1]$ and $\theta\in\mathbb R$, let
\begin{equation}\label{eq:prelim-gr}
 g_r(\theta)
 =-\Li_2(-ir e^{i\theta})-\Li_2(-ir e^{-i\theta}).
\end{equation}
For $r>0$, the integral representation gives
\[
 -2\Li_2(-ir)
 =2\int_0^r\frac{\log(1+iy)}{y}\,dy
 =\int_0^r\frac{\log(1+y^2)}{y}\,dy
  +2i\int_0^r\frac{\arctan y}{y}\,dy,
\]
whose real part equals $-\frac12\Li_2(-r^2)>0$ and whose imaginary part is
positive. Hence $-2\Li_2(-ir)$ lies in the open first quadrant, and the
principal square root
\begin{equation}\label{eq:parameters}
 s(r)=\sqrt{-2\Li_2(-ir)}=a(r)+ib(r),
 \qquad a(r),b(r)\in\mathbb R,
\end{equation}
satisfies $a(r)>b(r)>0$, $0<\arg s(r)<\pi/4$,
$a(r)^2-b(r)^2=-\frac12\Li_2(-r^2)$, and
$2a(r)b(r)=2\int_0^r\frac{\arctan y}{y}\,dy$.
In particular, $g_r(0)=s(r)^2$ for $0<r\le1$.
For every integer $n\ge1$ and $r>0$, define
\begin{equation}\label{eq:amplitude-phase}
 \begin{aligned}
 M_n(r)
 &=\frac{|s(r)|(1+r^2)^{1/4}}{2\pi n\sqrt r}
      e^{2a(r)\sqrt n},\\
 \Psi_n(r)
 &=2b(r)\sqrt n+\arg s(r)-\frac\pi4
      +\frac12\arctan r.
 \end{aligned}
\end{equation}
We call $M_n$ the \emph{amplitude} and $\Psi_n$ the \emph{phase}.
Both are continuous on $(0,\infty)$, and $M_n(r)>0$.

For $\mu\in\{0,1,2\}$, the Bessel functions $I_\mu$ and
$J_\mu$ are the entire functions
\begin{equation}\label{eq:prelim-bessel-series}
 I_\mu(z)
 =\sum_{k\ge0}\frac{(z/2)^{2k+\mu}}{k!(k+\mu)!},
 \qquad J_\mu(z)
 =\sum_{k\ge0}\frac{(-1)^k(z/2)^{2k+\mu}}{k!(k+\mu)!}.
\end{equation}
Termwise differentiation of these series gives $I_0'=I_1$ and
$J_0'=-J_1$, while termwise integration of the exponential series
yields the integral representations
\begin{equation}\label{eq:prelim-bessel-integrals}
 \begin{aligned}
 I_0(z)&=\frac1\pi\int_0^\pi e^{z\cos\varphi}\,d\varphi,\\
 I_1(z)&=\frac z\pi\int_{-1}^1
             e^{zy}(1-y^2)^{1/2}\,dy,\\
 I_2(z)&=\frac{z^2}{3\pi}\int_{-1}^1
             e^{zy}(1-y^2)^{3/2}\,dy.
 \end{aligned}
\end{equation}
For $\nu\in\{0,1\}$, let $H_\nu^{(1)}$ denote the Hankel function of the first kind.
For $x>0$, the expansion in \cite[Eq.~(91)]{Nemes}
with $N=1$ gives
\[
 H_\nu^{(1)}(x)
 =\sqrt{\frac{2}{\pi x}}\,e^{i(x-\nu\pi/2-\pi/4)}
   \bigl(1+R_1^{(K)}(-ix,\nu)\bigr),
\]
where $R_1^{(K)}$ is the remainder in the notation of that paper.
Since $|\nu|<3/2$ and $|\arg(-ix)|=\pi/2$,
\cite[Eqs.~(9) and~(75)]{Nemes} yields
\[
 |R_1^{(K)}(-ix,\nu)|
 \le\frac{|4\nu^2-1|}{8x}.
\]
Taking real parts and using
$J_\nu(x)=\Rea H_\nu^{(1)}(x)$, we obtain
\begin{equation}\label{eq:prelim-hankel-bound}
 \left|\sqrt{\frac{\pi x}{2}}J_\nu(x)
       -\cos\left(x-\frac{\nu\pi}{2}-\frac\pi4\right)\right|
 \le\frac{|4\nu^2-1|}{8x}
 \qquad(x>0,\ \nu=0,1).
\end{equation}

For $\Rea\zeta>0$ and $z\ne0$, let
\[
 \mathcal E(\zeta)=(e^{-\zeta};e^{-\zeta})_\infty,\qquad
 \vartheta(z;\zeta)
 =\sum_{j\in\mathbb Z}z^j e^{-\zeta j^2/2}.
\]
In this normalization, Jacobi's triple product
identity~\cite[Theorem~2.8]{Andrews} (with $q=e^{-\zeta/2}$)
and Jacobi's imaginary transformation
\cite[Chapter~10, Theorem~1.6]{SteinShakarchi} give
\begin{equation}\label{eq:theta-identities}
 \begin{aligned}
 \vartheta(z;\zeta)
 &=\mathcal E(\zeta)\prod_{j\ge0}
       (1+z e^{-(j+1/2)\zeta})
       (1+z^{-1}e^{-(j+1/2)\zeta}),\\
 \vartheta(e^w;\zeta)
 &=\sqrt{\frac{2\pi}{\zeta}}
       \sum_{j\in\mathbb Z}
       \exp\!\left(\frac{(w+2\pi ij)^2}{2\zeta}\right)
       \qquad(w\in\mathbb C).
 \end{aligned}
\end{equation}
Both the defining series for $\vartheta$ and the series in
\eqref{eq:theta-identities} converge absolutely for $\Rea\zeta>0$.
Recall Dedekind's eta function
$\eta(\tau)=e^{\pi i\tau/12}(e^{2\pi i\tau};e^{2\pi i\tau})_\infty$
($\Ima\tau>0$), so that $\eta(i\zeta/(2\pi))=e^{-\zeta/24}\mathcal E(\zeta)$.
Setting $\tau=i\zeta/(2\pi)$ in the modular transformation
$\eta(-1/\tau)=(-i\tau)^{1/2}\eta(\tau)$~\cite[Theorem~3.1]{Apostol} yields
\begin{equation}\label{eq:eta-transform}
 \mathcal E(\zeta)=\sqrt{\frac{2\pi}{\zeta}}
   \exp\!\left(\frac{\zeta}{24}-\frac{\pi^2}{6\zeta}\right)
   \mathcal E(4\pi^2/\zeta)
 \qquad(\Rea\zeta>0).
\end{equation}

A two-colour partition is an ordered pair of partitions, with
size equal to the sum of their sizes. For integers $n,k\ge0$,
let $p_2(n)$ count such pairs of size $n$, and let $p_{2,k}(n)$
count those in which every part is at most $k$. Then
\[
 p_2(n)=[q^n](q;q)_\infty^{-2},\qquad
 p_{2,k}(n)=[q^n](q;q)_k^{-2}.
\]
Both quantities are zero for $n<0$, and \eqref{eq:F} gives
$p(n,d)=p_{2,d}(n-d^2)$.
For integers $\ell\ge1$ and $h$ with $\gcd(h,\ell)=1$,
define the Dedekind sum
\[
 s(h,\ell)
 =\sum_{j=1}^{\ell-1}
    \left(\frac j\ell-\frac12\right)
    \left(\frac{hj}\ell
       -\left\lfloor\frac{hj}\ell\right\rfloor-\frac12\right),
\]
where an empty sum is zero.
For an integer $n\ge1$, define
\[
 A_\ell^{(2)}(n)
 =\sum_{\substack{0\le h<\ell\\\gcd(h,\ell)=1}}
   \exp\!\left(2\pi i\,s(h,\ell)
                    -\frac{2\pi inh}{\ell}\right).
\]
This is $A_\ell^{(2)}(n,0)$ in the notation of
\cite[Eq.~(1.7)]{IJT}.
In particular,
$A_1^{(2)}(n)=1$ and $|A_\ell^{(2)}(n)|\le\ell$.
The parameter-$2$ case of the Rademacher formula of
Iskander, Jain, and Talvola
\cite[Theorem~1.1 and Eq.~(1.7)]{IJT} becomes
\begin{equation}\label{eq:two-colour-rademacher}
 p_2(n)
 =\frac{2\pi}{12n-1}
   \sum_{\ell\ge1}\frac{A_\ell^{(2)}(n)}{\ell}
       I_2\!\left(\frac{\pi\sqrt{12n-1}}{3\ell}\right)
 \qquad(n\ge1).
\end{equation}
For fixed $n$, the series converges absolutely: the defining
series of $I_2$ and the bound $|A_\ell^{(2)}(n)|\le\ell$ show that
its $\ell$th summand is $O_n(\ell^{-2})$.

\section{Sign estimates near the origin and in the intermediate range}\label{sec:central}

In this section we establish effective estimates for $D_n(-r^2)$ in two
regimes: near the origin, where $0<r\le16\pi^2/n$, and in the
intermediate range $16\pi^2/n\le r\le e\sqrt n$ for $n\ge6400$.
Near the origin, we approximate $D_n(-r^2)$ by products of Bessel
functions (Lemma~\ref{lem:effective-bessel}) and determine its sign at
$r=(k\pi)^2/(4n)$ (Corollary~\ref{cor:effective-bessel-sign}).
In the intermediate range, we prove uniform approximations of the form
$|D_n(-r^2)-M_n(r)\cos\Psi_n(r)|<cM_n(r)$ with $c<1$, treating
$16\pi^2/n\le r\le1$ (Proposition~\ref{lem:effective-unit-left}) and
$1\le r\le e\sqrt n$ (Proposition~\ref{thm:effective-right}) separately;
in particular, $\operatorname{sgn}D_n(-r^2)=(-1)^k$ whenever
$\Psi_n(r)=k\pi$ ($k\in\mathbb Z$). Larger radii are treated in
Section~\ref{sec:outer}.
Throughout this section, write
$\langle H\rangle=(2\pi)^{-1}\int_{-\pi}^{\pi}H(\theta)\,d\theta$
for the angular average of $H$.

For an integer $n\ge1$ and $\sigma>0$, set $r=\sigma^2/n$,
$\Delta=\sigma^3/n$, and
\[
 \varepsilon=\left(\frac{0.66}{n}+0.228r
                              +\frac{0.35r}{n}\right)\Delta,
\]
and define
\[
 \mathcal B(\sigma,n)
 =\frac{15}{2\pi}\left(
       \varepsilon e^{\Delta/2+\varepsilon}
       +\frac{\Delta^2}{8}e^{\Delta/2}\right)
 +\Delta e^{-\sigma/2}
   \left((0.48\sigma+0.36\sigma^2)e^{1.02\Delta}
       +\frac{\sigma}{\pi^2}+\frac{\sigma^2}{4}\right).
\]

\begin{lemma}\label{lem:effective-bessel}
Let $n\ge1$ be an integer and $0<r\le1/40$.
Set $\sigma=\sqrt{nr}$ and $\Delta=\sigma^3/n$, and assume that
$\sigma\ge1$ and $\Delta\le1$. Then
\[
 \bigg|D_n(-r^2)
 -\frac{\sigma}{n}
   \bigl(J_0(2\sigma)I_1(2\sigma)-J_1(2\sigma)I_0(2\sigma)\bigr)
 +\frac{\sigma^4}{2n^2}J_2(2\sigma)I_2(2\sigma)\bigg|
 \le\frac{e^{2\sigma}}{n}\mathcal B(\sigma,n).
\]
\end{lemma}
\begin{proof}
Let $u=\sigma/n=r/\sigma$ and $\zeta=u+iv$. Then $u\le1/40$,
$nu=r/u=\sigma$, $r^2/u=\Delta$, $ru=\Delta/n$, and $r\le\Delta$.
Expanding the principal logarithms in \eqref{eq:P-product} under the
convention of Section~\ref{sec:preliminaries} gives
\begin{equation}\label{eq:bessel-L-series}
 \begin{aligned}
 \mathcal L:=\log\mathcal P_r(\zeta,\theta)
 &=\sum_{j\ge0}\sum_{\ell\ge1}
   \frac{(-1)^{\ell+1}(ir)^\ell}{\ell}
   e^{-\ell(j+1/2)\zeta}\bigl(e^{i\ell\theta}+e^{-i\ell\theta}\bigr)\\
 &=\sum_{\ell\ge1}
   \frac{2(-1)^{\ell+1}(ir)^\ell\cos(\ell\theta)}{\ell}
   \sum_{j\ge0}e^{-\ell(j+1/2)\zeta}\\
 &=\sum_{\ell\ge1}
   \frac{2(-1)^{\ell+1}(ir)^\ell\cos(\ell\theta)}
        {\ell\,2\sinh(\ell\zeta/2)}.
 \end{aligned}
\end{equation}
For fixed $n$ and $r$, the interchange of summation and the
absolute and uniform convergence in $(v,\theta)\in\mathbb R^2$
follow from
\[
 \sum_{j\ge0}\sum_{\ell\ge1}
 \frac{2r^\ell}{\ell}e^{-\ell(j+1/2)u}
 =\sum_{\ell\ge1}\frac{2r^\ell}{\ell\,2\sinh(\ell u/2)}
 \le\frac{2}{u}\sum_{\ell\ge1}\frac{r^\ell}{\ell^2}<\infty,
\]
together with $|2\sinh(\ell\zeta/2)|\ge 2\sinh(\ell u/2)\ge\ell u$ and $r<1$.
Set
\[
 \mathcal L_0=\frac{2ir\cos\theta}{\zeta},
 \qquad \mathcal S=\frac{r^2\cos(2\theta)}{2\zeta}.
\]
Since $\langle\cos(\ell\theta)\rangle=0$ for every integer $\ell\ge1$,
uniform convergence in $\theta$ justifies termwise integration of
\eqref{eq:bessel-L-series}, giving
$\langle\mathcal L\rangle=\langle\mathcal L_0\rangle=\langle\mathcal S\rangle=0$.
As $n\ge1$ is an integer, we have $\int_{-\pi}^{\pi}e^{n(u+iv)}\,dv=0$;
\eqref{eq:ct-integral} and \eqref{eq:ct-integral2} yield
\[
 D_n(-r^2)
 =\frac1{2\pi}\int_{-\pi}^{\pi}
       e^{n\zeta}\langle e^{\mathcal L}\rangle\,dv
 =\frac1{2\pi}\int_{-\pi}^{\pi}
       e^{n\zeta}\langle e^{\mathcal L}-1\rangle\,dv.
\]
We next evaluate the integral
\[
 \frac1{2\pi}\int_{\mathbb R}e^{n\zeta}
 \left\langle e^{\mathcal L_0}(1+\mathcal S)-1\right\rangle\,dv,
\]
whose absolute convergence is established below.
Multiplying the series in \eqref{eq:prelim-bessel-series} gives
\begin{equation}\label{eq:bessel-product-series}
 J_\nu(x)I_\nu(x)
 =\sum_{d\ge0}\frac{(-1)^d(x/2)^{4d+2\nu}}
                       {d!(d+\nu)!(2d+\nu)!}
 \qquad(\nu=0,2).
\end{equation}
Indeed, the terms with odd total index cancel, whereas
\[
 \sum_{j=0}^{2d}
   \frac{(-1)^j}{j!(2d-j)!(j+\nu)!(2d-j+\nu)!}=
 \frac{[z^{2d+\nu}](1-z^2)^{2d}(1+z)^{2\nu}}
      {(2d)!(2d+2\nu)!}
 =\frac{(-1)^d}{d!(d+\nu)!(2d+\nu)!};
\]
the coefficient in the numerator is $(-1)^d\binom{2d}{d}$ for $\nu=0$
and $(-1)^d\bigl(6\binom{2d}{d}-2\binom{2d}{d-1}\bigr)$ for $\nu=2$.
Expanding
$e^{\mathcal L_0}
=\exp((ir/\zeta)e^{i\theta})\exp((ir/\zeta)e^{-i\theta})$
and taking angular averages yields
\[
 \langle e^{\mathcal L_0}(1+\mathcal S)-1\rangle
 =\sum_{d\ge1}\frac{(-r^2)^d}{(d!)^2\zeta^{2d}}
  -\frac{r^4}{2}\sum_{d\ge0}
           \frac{(-r^2)^d}{d!(d+2)!\zeta^{2d+3}}.
\]
Since $|\zeta|\ge u$ and $r/u=\sigma$, the sum of the absolute values
of the terms on the right-hand side is bounded by
$e^{2\sigma}(r^2+r^4/(4u))/(u^2+v^2)$, which is integrable on
$\mathbb R$. Thus termwise integration is justified. For every integer
$j\ge2$, the residue theorem gives
\[
 \frac1{2\pi}\int_{\mathbb R}e^{n\zeta}\zeta^{-j}\,dv
 =\frac{n^{j-1}}{(j-1)!},
\]
since the integral on the closing left semicircle centred at $u$
with radius $R>2u$ has modulus at most $\pi R e^\sigma(R-u)^{-j}\to0$ as $R\to\infty$.
Using \eqref{eq:bessel-product-series}, $J_0'=-J_1$, $I_0'=I_1$, and
$\sigma^2=nr$, we obtain
\begin{equation}\label{eq:effective-bessel-inverse}
 \begin{aligned}
 \frac1{2\pi}\int_{\mathbb R}e^{n\zeta}
       \langle e^{\mathcal L_0}(1+\mathcal S)-1\rangle\,dv
 &=\frac1n\sum_{d\ge1}
          \frac{(-1)^d\sigma^{4d}}{(d!)^2(2d-1)!}
       -\frac{r^2}{2}\sum_{d\ge0}
          \frac{(-1)^d\sigma^{4d+4}}{d!(d+2)!(2d+2)!}\\
 &=\frac\sigma n
       \bigl(J_0(2\sigma)I_1(2\sigma)-J_1(2\sigma)I_0(2\sigma)\bigr)
       -\frac{\sigma^4}{2n^2}J_2(2\sigma)I_2(2\sigma).
 \end{aligned}
\end{equation}
Since $4u<\pi$ and $|e^{n\zeta}|=e^\sigma$, splitting at
$|v|=4u$ and $|v|=\pi$ gives
\begin{equation}\label{eq:bessel-four-parts}
 \begin{aligned}
 &\left|D_n(-r^2)-\frac1{2\pi}\int_{\mathbb R}e^{n\zeta}
       \langle e^{\mathcal L_0}(1+\mathcal S)-1\rangle\,dv\right|
\le\frac{e^\sigma}{2\pi}\Bigg(
      \int_{|v|\le4u}
      \left|\langle e^{\mathcal L}-e^{\mathcal L_0}(1+\mathcal S)
      \rangle\right|\,dv\\
 &\qquad\qquad+\int_{4u<|v|\le\pi}
      \left|\langle e^{\mathcal L}-e^{\mathcal L_0}\rangle\right|\,dv
      +\int_{|v|>\pi}
      \left|\langle e^{\mathcal L_0}-1\rangle\right|\,dv
+\int_{|v|>4u}
      \left|\langle e^{\mathcal L_0}\mathcal S\rangle\right|\,dv
      \Bigg).
 \end{aligned}
\end{equation}

Before estimating the four terms on the right-hand side of
\eqref{eq:bessel-four-parts}, we decompose $\mathcal L-\mathcal L_0$
on $|v|\le\pi$. Write $\mathcal L-\mathcal L_0=\rho_1+\rho_2$, where
\[
 \rho_1=2ir\cos\theta
       \left(\frac1{2\sinh(\zeta/2)}-\frac1\zeta\right).
\]
For $|v|\le\pi$, we have $|\zeta|<3.142$ and $|\zeta|^2<10$.
The power series for $\sinh$ and the identity
$|2\sinh(\zeta/2)|^2=4\sinh^2(u/2)+4\sin^2(v/2)$ give
\[
 |2\sinh(\zeta/2)-\zeta|
\le\frac{|\zeta|^3}{24(1-|\zeta|^2/80)}<0.05|\zeta|^3,
\qquad |2\sinh(\zeta/2)|^2
   \ge\frac{4|\zeta|^2}{\pi^2}.
\]
Consequently, for $|v|\le\pi$,
\begin{equation}\label{eq:bessel-rho-bounds}
 |\rho_1|\le0.16r|\zeta||\cos\theta|,\qquad
 |\rho_2|\le\frac2u\sum_{\ell\ge2}\frac{r^\ell}{\ell^2}
       \le\frac{r^2}{2u(1-r)}<0.513\Delta,\qquad
 |\mathcal L-\mathcal L_0|<1.02\Delta.
\end{equation}

For the first term in \eqref{eq:bessel-four-parts}, on $|v|\le4u$ we
have $|\zeta|\le\sqrt{17}\,u<0.104$. The power series for $\sinh$ gives
\[
 |\sinh\zeta-\zeta|
 \le\frac{|\zeta|^3}{6(1-|\zeta|^2/20)}.
\]
The reverse triangle inequality therefore gives
\[
 |\sinh\zeta|
 \ge|\zeta|\,
       \frac{1-13|\zeta|^2/60}{1-|\zeta|^2/20},
\]
and hence
\[
 \left|\frac1{2\sinh\zeta}-\frac1{2\zeta}\right|
 \le\frac{|\zeta|}{12(1-13|\zeta|^2/60)}<0.084|\zeta|.
\]
Separating the terms with $\ell=1$, $\ell=2$, and $\ell\ge3$ in
\eqref{eq:bessel-L-series} and using \eqref{eq:bessel-rho-bounds}, we
obtain
\[
 |\mathcal L-\mathcal L_0-\mathcal S|
 \le0.16\sqrt{17}\,ru+0.084\sqrt{17}\,r^2u
                       +\frac{2r^3}{9u(1-r)}
 \le\frac{0.66\Delta}{n}+\frac{0.35r\Delta}{n}+0.228r\Delta
 =\varepsilon.
\]
Also, $|\mathcal S|\le\Delta/2$. The elementary Taylor bounds
\begin{equation}\label{eq:bessel-exp-taylor}
 |e^z-1|\le|z|e^{\max\{0,\Rea z\}},
 \qquad
 |e^z-1-z|\le\frac12|z|^2e^{\max\{0,\Rea z\}}
 \qquad(z\in\mathbb C)
\end{equation}
therefore imply
\[
 |e^{\mathcal L}-e^{\mathcal L_0}(1+\mathcal S)|
 \le e^{\Rea\mathcal L_0}
       \left(\varepsilon e^{\Delta/2+\varepsilon}
                          +\frac{\Delta^2}{8}e^{\Delta/2}\right).
\]
With $v=u\eta$, we have
$\Rea\mathcal L_0=2\sigma\eta\cos\theta/(1+\eta^2)$, and
\[
 1-\frac{2\eta\cos\theta}{1+\eta^2}
 =\frac{(\eta-1)^2+2\eta(1-\cos\theta)}{1+\eta^2}
 \ge\frac{(\eta-1)^2}{5}+\frac{8\theta^2}{5\pi^2}
 \quad\left(\frac12\le\eta\le2,\ |\theta|\le\pi\right).
\]
This inequality implies
\[
 \int_{1/2}^{2}
       \langle e^{\Rea\mathcal L_0}\rangle\,d\eta
 \le\frac{e^\sigma}{2\pi}
   \int_{\mathbb R}e^{-\sigma(\eta-1)^2/5}\,d\eta
   \int_{\mathbb R}e^{-8\sigma\theta^2/(5\pi^2)}\,d\theta
 =\frac{5\pi e^\sigma}{2\sqrt8\,\sigma}.
\]
The average $\langle e^{\Rea\mathcal L_0}\rangle$ is even in $\eta$,
by the change $\theta\mapsto\theta+\pi$ and periodicity.
Thus $\int_{-2}^{-1/2}\langle e^{\Rea\mathcal L_0}\rangle\,d\eta
=\int_{1/2}^{2}\langle e^{\Rea\mathcal L_0}\rangle\,d\eta$.
On $[-4,-2]\cup[-1/2,1/2]\cup[2,4]$, which has total length $5$, we have
$2|\eta|/(1+\eta^2)\le4/5$. Consequently,
\[
 \int_{-4}^4\langle e^{\Rea\mathcal L_0}\rangle\,d\eta
 \le e^\sigma\left(
       \frac{5\pi}{\sqrt8\,\sigma}+5e^{-\sigma/5}\right)
 <\frac{15e^\sigma}{\sigma},
\]
using $\sigma e^{-\sigma/5}\le5/e$ and $5\pi/\sqrt8+25/e<15$.
Since $u/\sigma=1/n$, we obtain
\[
 \frac{e^\sigma}{2\pi}\int_{|v|\le4u}
 \left|\langle e^{\mathcal L}-e^{\mathcal L_0}(1+\mathcal S)\rangle\right|\,dv
 \le\frac{15e^{2\sigma}}{2\pi n}
       \left(\varepsilon e^{\Delta/2+\varepsilon}
                          +\frac{\Delta^2}{8}e^{\Delta/2}\right).
\]

For the second term in \eqref{eq:bessel-four-parts}, define
\[
 \|H\|_2^2
 =\frac1{2\pi}\int_{-\pi}^{\pi}\langle|H|^2\rangle\,dv.
\]
Since $\langle\cos^2\theta\rangle=1/2$,
\[
 \|\mathcal L_0\|_2^2
 =\frac{r^2}{\pi}\int_{-\pi}^{\pi}\frac{dv}{u^2+v^2}
 \le\frac{r^2}{u}.
\]
Moreover, \eqref{eq:bessel-rho-bounds} gives
$|\mathcal L_0\rho_1|\le0.32r^2\cos^2\theta$; combining this with the
bound for $|\rho_1|$ in \eqref{eq:bessel-rho-bounds} yields
\[
 \|\rho_1\|_2
 \le0.16r\sqrt{\frac{u^2+\pi^2/3}{2}}<0.206r,
 \qquad
 \frac1{2\pi}\int_{-\pi}^{\pi}
       \langle|\mathcal L_0\rho_1|\rangle\,dv
 \le0.16r^2.
\]
Expanding $1/(2\sinh(\ell\zeta/2))=\sum_{j\ge0}e^{-\ell(j+1/2)\zeta}$
and applying Parseval's identity in $v$ gives, for each integer
$\ell\ge1$,
\[
 \frac1{2\pi}\int_{-\pi}^{\pi}
       \frac{dv}{|2\sinh(\ell\zeta/2)|^2}
 =\frac{e^{-\ell u}}{1-e^{-2\ell u}}
 =\frac1{2\sinh(\ell u)}\le\frac1{2\ell u}.
\]
By orthogonality of the cosines $\cos(\ell\theta)$ in $\rho_2$, it
follows that
\[
 \|\rho_2\|_2^2
 =2\sum_{\ell\ge2}\frac{r^{2\ell}}{\ell^2}
                      \frac{e^{-\ell u}}{1-e^{-2\ell u}}
 \le\frac{r^4}{8u(1-r^2)}<\frac{(0.354r^2)^2}{u}.
\]
On $|v|>4u$,
$\Rea\mathcal L_0\le2r|v|/(u^2+v^2)\le8\sigma/17<\sigma/2$.
The cancellation $\langle\rho_1+\rho_2\rangle=0$ gives the identity
\[
 \langle e^{\mathcal L}-e^{\mathcal L_0}\rangle
 =\langle(e^{\mathcal L_0}-1)(\rho_1+\rho_2)\rangle
 +\langle e^{\mathcal L_0}
           (e^{\rho_1+\rho_2}-1-\rho_1-\rho_2)\rangle.
\]
Applying \eqref{eq:bessel-exp-taylor}, \eqref{eq:bessel-rho-bounds},
and Cauchy--Schwarz, we obtain
\[
 \begin{aligned}
 \frac{e^\sigma}{2\pi}\int_{4u<|v|\le\pi}
       \left|\langle e^{\mathcal L}-e^{\mathcal L_0}\rangle\right|\,dv
 &\le e^{3\sigma/2+1.02\Delta}
       \left[0.16r^2+\|\mathcal L_0\|_2\|\rho_2\|_2
                   +\frac12(\|\rho_1\|_2+\|\rho_2\|_2)^2\right]\\[2pt]
 &\le e^{3\sigma/2+1.02\Delta}
       \left[0.16r^2+0.354\frac{r^3}{u}
             +\frac12\left(0.206r+\frac{0.354r^2}{\sqrt u}\right)^2\right]\\[2pt]
 &<e^{3\sigma/2+1.02\Delta}
              \left(0.48r^2+0.36\frac{r^3}{u}\right),
 \end{aligned}
\]
where the last step uses $r/\sqrt u=\sqrt\Delta\le1$.

For the third and fourth terms in \eqref{eq:bessel-four-parts},
since $\langle\mathcal L_0\rangle=\langle\mathcal S\rangle=0$ and
$\Rea\mathcal L_0<\sigma/2$ on $|v|>4u$, \eqref{eq:bessel-exp-taylor}
and Cauchy--Schwarz in $\theta$ give
\[
 |\langle e^{\mathcal L_0}-1\rangle|
 \le\tfrac12e^{\sigma/2}\langle|\mathcal L_0|^2\rangle
   =\frac{e^{\sigma/2}r^2}{u^2+v^2},
 \qquad
 |\langle e^{\mathcal L_0}\mathcal S\rangle|
 \le e^{\sigma/2}\langle|\mathcal L_0\mathcal S|\rangle
   \le\frac{e^{\sigma/2}r^3}{2(u^2+v^2)}
 \qquad(|v|>4u).
\]
Since $\int_{|v|>\pi}(u^2+v^2)^{-1}\,dv\le2/\pi$ and
$\int_{\mathbb R}(u^2+v^2)^{-1}\,dv=\pi/u$, these imply
\[
 \frac{e^\sigma}{2\pi}\int_{|v|>\pi}
       |\langle e^{\mathcal L_0}-1\rangle|\,dv
 \le\frac{e^{3\sigma/2}r^2}{\pi^2},
 \qquad
 \frac{e^\sigma}{2\pi}\int_{|v|>4u}
       |\langle e^{\mathcal L_0}\mathcal S\rangle|\,dv
 \le\frac{e^{3\sigma/2}r^3}{4u}.
\]
Substituting these four estimates into \eqref{eq:bessel-four-parts}
and using $nr^2=\sigma\Delta$ and $nr^3/u=\sigma^2\Delta$ bounds
its right-hand side by $e^{2\sigma}\mathcal B(\sigma,n)/n$.
Equation~\eqref{eq:effective-bessel-inverse} now proves the lemma.
\end{proof}

\begin{corollary}\label{cor:effective-bessel-sign}
For an integer $k\ge1$, let $\sigma=k\pi/2$.
If $n\ge\max\{6400,40\sigma^2,3\sigma^3\}$ is an integer, then
$
 \operatorname{sgn}D_n\!\left(-\frac{(k\pi)^4}{16n^2}\right)=(-1)^k.
$
In particular, this holds for $1\le k\le8$ whenever $n\ge6400$.
\end{corollary}
\begin{proof}
Let $x=2\sigma=k\pi$ and $r=\sigma^2/n$. Then $r\le1/40$,
$\sigma\ge\pi/2>1$, and $\Delta\le1/3$, so
Lemma~\ref{lem:effective-bessel} applies. By
\eqref{eq:prelim-hankel-bound},
\[
 (-1)^{k+\nu}J_\nu(x)
 \ge\sqrt{\frac2{\pi x}}
       \left(\frac1{\sqrt2}-\frac{|4\nu^2-1|}{8x}\right)>0
 \qquad(\nu=0,1).
\]
Since $x\ge\pi$, the representation \eqref{eq:prelim-bessel-integrals}
of $I_0$, together with $\cos\theta\ge1-\theta^2/2$ and
$\int_\pi^\infty e^{-x\theta^2/2}\,d\theta\le e^{-x\pi^2/2}/(\pi x)$,
gives
\[
 I_0(x)\ge\frac{e^x}{\pi}
       \left(\sqrt{\frac\pi{2x}}
                    -\frac{e^{-x\pi^2/2}}{\pi x}\right)
 >\frac{99e^x}{100\sqrt{2\pi x}},
\]
because
\[
 \frac{e^{-x\pi^2/2}/(\pi x)}{\sqrt{\pi/(2x)}}
 =\frac{\sqrt2\,e^{-x\pi^2/2}}{\pi^{3/2}\sqrt x}
 \le\frac{\sqrt2\,e^{-\pi^3/2}}{\pi^2}
 <\frac1{100}.
\]
Also, $I_0'=I_1$ and $I_1'(t)=I_0(t)-I_1(t)/t$ (from
\eqref{eq:prelim-bessel-series}) imply
\[
 \frac{d}{dt}\left[
 e^t\left(I_1(t)-\left(1-\frac1t\right)I_0(t)\right)\right]
 =e^t\frac{t-1}{t^2}I_0(t)\ge0
 \qquad(t\ge1).
\]
The expression in brackets equals $eI_1(1)>0$ at $t=1$; hence
$I_1(t)>(1-1/t)I_0(t)$ for $t\ge1$, and in particular
$I_1(x)>(1-1/x)I_0(x)$. Combining these bounds yields
\[
 (-1)^k\frac x2\bigl(J_0(x)I_1(x)-J_1(x)I_0(x)\bigr)
 >\frac{99e^x}{200\pi}
       \left(\frac{2-1/x}{\sqrt2}-\frac1{2x}\right)
 >\frac{63e^x}{400}.
\]
For the last inequality, note that
$\frac{2-1/x}{\sqrt2}-\frac1{2x}-\bigl(\frac75-\frac6{5x}\bigr)
=\bigl(\sqrt2-\frac75\bigr)\bigl(1-\frac1{2x}\bigr)>0$.
Since $x\ge\pi>3$, we have $7/5-6/(5x)>1$, while $\pi<22/7$ gives
$\frac{99}{200\pi}>\frac{99\cdot7}{200\cdot22}=\frac{63}{400}$.
By \eqref{eq:prelim-bessel-series}, $J_2(x)=2J_1(x)/x-J_0(x)$, so the
signs of $J_0(x)$ and $J_1(x)$ established above give
\[
 (-1)^{k+1}J_2(x)
 =\frac2x(-1)^{k+1}J_1(x)+(-1)^kJ_0(x)>0.
\]
Since $x=2\sigma$ and $I_2(x)>0$, it follows that
\[
 (-1)^k\left(
 -\frac{\sigma^4}{2n^2}J_2(2\sigma)I_2(2\sigma)
 \right)>0.
\]

We next bound $\mathcal B(\sigma,n)$, which controls the
remainder in Lemma~\ref{lem:effective-bessel} through
$e^{2\sigma}\mathcal B(\sigma,n)/n$.
The hypotheses give
\[
 \varepsilon\le\frac13\left(
       \frac{0.66}{6400}+\frac{0.228}{40}
                     +\frac{0.35}{40\cdot6400}\right)<\frac1{500}.
\]
Use $\pi>3$ and the elementary bounds
$e^{1/6+1/500}<5/4$, $e^{17/50}<3/2$, $e^4>50$, and $e^5>140$.
Since $\Delta=\sigma^3/n$, $n\ge6400$, and
$\max_{t>0}t^j e^{-t/2}=(2j)^j e^{-j}$ for $j=4,5$, we obtain
\[
 \mathcal B(\sigma,n)
 <\frac{25}{8}\left(\frac1{500}+\frac1{72}\right)
+\frac1{6400}\left[
   \frac32\left(\frac{0.48\cdot8^4}{50}
                      +\frac{0.36\cdot10^5}{140}\right)
        +\frac{8^4}{9\cdot50}+\frac{10^5}{4\cdot140}\right]
 <\frac3{20}.
\]
Lemma~\ref{lem:effective-bessel} therefore gives
\[
 \begin{aligned}
 (-1)^kD_n(-r^2)
 &\ge
 (-1)^k\left[
 \frac{\sigma}{n}
 \bigl(J_0(2\sigma)I_1(2\sigma)
       -J_1(2\sigma)I_0(2\sigma)\bigr)
 -\frac{\sigma^4}{2n^2}J_2(2\sigma)I_2(2\sigma)
 \right]
 -\frac{e^{2\sigma}}n\mathcal B(\sigma,n)\\
 &>
 \left(\frac{63}{400}-\frac3{20}\right)
 \frac{e^{2\sigma}}n
 =\frac{3e^{2\sigma}}{400n}>0.
 \end{aligned}
\]
This is the required assertion because $r^2=(k\pi)^4/(16n^2)$.
For $1\le k\le8$, we have $\sigma\le4\pi$,
$40(4\pi)^2<6400$, and $3(4\pi)^3<6400$, proving the final statement.
\end{proof}

For $0<r\le1$, let $g_r$ be defined
by \eqref{eq:prelim-gr}. Write $s=s(r)$, $a=a(r)$,
$b=b(r)$ as in \eqref{eq:parameters}.
Let
\[
 (\gamma,C(\gamma))=
 \begin{cases}
 (0.99,6.273),&0<r\le0.1,\\
 (0.55,15.961),&0.1<r\le1.
 \end{cases}
\]
We first record basic parameter estimates for $0<r\le1$.

\begin{lemma}\label{lem:effective-unit-bounds}
Let $0<r\le1$, and set $\varrho=b/a$ and $\kappa=r/((1+r^2)a^2)$.
For all $\theta\in\mathbb R$, we have
\begin{equation}\label{eq:unit-log-integral}
 g_r(\theta)
 =\int_0^r
   \log(1-y^2+2iy\cos\theta)\,\frac{dy}{y},
\end{equation}
and
\begin{equation}\label{eq:eff-unit-A-gap}
 \begin{aligned}
 \Rea(s^2)-\Rea g_r(\theta)
 &\ge\frac{r^2}{1+r^2}\sin^2\theta,\\
 \Ima(s^2)-\Ima g_r(\theta)
 &\ge\frac{2r}{1+r^2}(1-\cos\theta).
 \end{aligned}
\end{equation}
The parameters satisfy
\begin{equation}\label{eq:unit-basic-bounds}
 r\le a^2\le r(1+0.282r),
 \qquad
 \varrho\in[1-0.282r,1],
 \qquad
 |s|^2\ge\frac{16r}{9},
 \qquad
 \kappa>0.39.
\end{equation}
Moreover, for every integer $n\ge1$ with $\sqrt{nr}\ge4\pi$ and
$\omega=a\sqrt n$, we have
\begin{equation}\label{eq:unit-scale-bounds}
 \omega\ge4\pi,
 \qquad M_n(r)\ge\frac{2e^{2\omega}}{3\pi n}.
\end{equation}
For $|\theta|\le\pi/2$ and $\eta\ge0$, with $\psi=1-\cos\theta$, we have
\begin{equation}\label{eq:unit-angular-decay}
 1-\frac{\Rea g_r(\theta)+\eta\Ima g_r(\theta)}{a^2(1+\eta^2)}
 \ge\frac{(\eta-\varrho)^2+r\kappa\psi(2-\psi)
                  +2\kappa\eta\psi}{1+\eta^2}
 \ge\gamma\psi,
\end{equation}
and, with $\sqrt{g_r(\theta)}$ at $\theta=\pm\pi/2$ understood as the
limit from $|\theta|<\pi/2$,
\begin{equation}\label{eq:eff-unit-sqrt-gap}
 \Rea\sqrt{g_r(\theta)}\le a\sqrt{1-\gamma\psi}.
\end{equation}
Finally, for all $\theta\in\mathbb R$ and $\eta\ge0$,
\begin{equation}\label{eq:unit-contour-decay}
 1-\frac{\Rea g_r(\theta)+\eta\Ima g_r(\theta)}{a^2(1+\eta^2)}
 \ge\frac{(\eta-\varrho)^2}{1+\eta^2};
\end{equation}
when $\pi/2<|\theta|\le\pi$, the left-hand side of
\eqref{eq:unit-contour-decay} is at least $1-r/2$.
\end{lemma}
\begin{proof}
For $0<r<1$, differentiating \eqref{eq:prelim-gr} gives
\[
 \frac{\partial g_r(\theta)}{\partial r}
 =\frac{\log(1+ir e^{i\theta})
        +\log(1+ir e^{-i\theta})}{r}
 =\frac{\log(1-r^2+2ir\cos\theta)}r.
\]
The principal logarithms add because
$1+ir e^{\pm i\theta}$ and $1-r^2+2ir\cos\theta$
have positive real parts.
Since $g_0(\theta)=0$, integration from $0$ to $r$ yields
\eqref{eq:unit-log-integral}.
Since $s^2=g_r(0)$, \eqref{eq:unit-log-integral} gives
\[
 \Rea(s^2)-\Rea g_r(\theta)
 =-\frac12\int_0^r
   \log\left(1-\frac{4y^2\sin^2\theta}{(1+y^2)^2}\right)
   \frac{dy}{y}
 \ge2\sin^2\theta\int_0^r\frac{y\,dy}{(1+y^2)^2}
 =\frac{r^2}{1+r^2}\sin^2\theta,
\]
where the inequality uses $-\log(1-x)\ge x$ for $0\le x<1$.
For $0<y<1$ and $-1\le q\le1$, we have
\[
 \frac{\partial}{\partial q}\arg(1-y^2+2iyq)
 =\frac{2y(1-y^2)}{(1-y^2)^2+4y^2q^2}
 \ge\frac{2y(1-y^2)}{(1+y^2)^2}.
\]
By \eqref{eq:unit-log-integral},
\[
 \begin{aligned}
 \Ima(s^2)-\Ima g_r(\theta)
 &=\int_0^r
   \bigl[\arg(1-y^2+2iy)
        -\arg(1-y^2+2iy\cos\theta)\bigr]\,\frac{dy}{y}\\
 &=\int_0^r\int_{\cos\theta}^1
   \frac{2(1-y^2)}{(1-y^2)^2+4y^2q^2}\,dq\,dy\\
 &\ge2(1-\cos\theta)
   \int_0^r\frac{1-y^2}{(1+y^2)^2}\,dy
 =\frac{2r}{1+r^2}(1-\cos\theta),
 \end{aligned}
\]
proving \eqref{eq:eff-unit-A-gap} for $0<r<1$.
The integral in \eqref{eq:unit-log-integral} remains convergent
at $r=1$: its only possible endpoint singularity occurs when
$\cos\theta=0$, in which case the logarithm is $\log(1-y^2)$.
Thus \eqref{eq:unit-log-integral} and both inequalities in
\eqref{eq:eff-unit-A-gap} extend to $r=1$ by continuity.

Expanding $s^2=-2\Li_2(-ir)$ from \eqref{eq:parameters} into its real and
imaginary parts gives the alternating series
\[
 \Rea(s^2)=2\sum_{j\ge1}\frac{(-1)^{j-1}r^{2j}}{(2j)^2},
 \qquad
 \frac{\Ima(s^2)}{2}=\sum_{j\ge0}\frac{(-1)^j r^{2j+1}}{(2j+1)^2}.
\]
Since $0<r\le1$, the terms decrease in modulus, so the first partial sums
give
\[
 \frac{3r^2}{8}\le\frac{r^2}{2}-\frac{r^4}{8}\le\Rea(s^2)\le\frac{r^2}{2},
 \qquad
 r\left(1-\frac{r^2}{9}\right)
 \le\frac{\Ima(s^2)}{2}
 \le r\left(1-\frac{r^2}{9}+\frac{r^4}{25}\right)\le r.
\]
Consequently,
\[
 \begin{aligned}
 a^2&\ge\frac{\Rea(s^2)+\Ima(s^2)}2
       \ge r+\frac{3r^2}{16}-\frac{r^3}{9}\ge r,\\
 a^2&\le r\left(1+\frac r4+\frac{r^2}{32}\right)
       \le r(1+0.282r),\\
 \varrho
 &=\sqrt{1+\left(\frac{\Rea(s^2)}{\Ima(s^2)}\right)^2}
       -\frac{\Rea(s^2)}{\Ima(s^2)}
 \in[1-0.282r,1],
 \end{aligned}
\]
where the upper bound for $a^2$ uses $|s|^2\le 2r\sqrt{1+r^2/16}$ and
$\sqrt{1+x}\le1+x/2$, and the bound for $\varrho$ uses
$\Rea(s^2)/\Ima(s^2)\le9r/32$. Together with $|s|^2\ge\Ima(s^2)\ge2r(1-r^2/9)\ge16r/9$ and
$\kappa\ge\bigl((1+r^2)(1+0.282r)\bigr)^{-1}\ge(2\cdot1.282)^{-1}>0.39$,
these estimates prove \eqref{eq:unit-basic-bounds}. The bounds in
\eqref{eq:unit-scale-bounds} follow immediately from
\eqref{eq:unit-basic-bounds} and \eqref{eq:amplitude-phase}. We also have
\[
 \begin{aligned}
 \varrho^2+r\kappa
 &=1+\frac{r^2/(1+r^2)-\Rea(s^2)}{a^2}\ge1,\\
 0\le\varrho-\kappa
 &=\frac{\Ima(s^2)/2-r/(1+r^2)}{a^2}
 \le\frac{r^2}{1+r^2}-\frac{r^2}{9}+\frac{r^4}{25}.
 \end{aligned}
\]
Here $\varrho^2+r\kappa\ge1$ uses $r^2/(1+r^2)\ge r^2/2\ge\Rea(s^2)$,
and $\varrho-\kappa\ge0$ uses $1-r^2/9\ge(1+r^2)^{-1}$ on $[0,1]$.
The upper bound for $\varrho-\kappa$ is at most $r^2$ and, since
$t/(1+t)-t/9+t^2/25$ increases on $[0,1]$, is at most $193/450$.
Hence $1-\varrho+\kappa\ge\gamma$.

For $A>0$, $B\ge0$, and $z=(A+iB)^2$, we have
\begin{equation}\label{eq:unit-quad-id}
 A^2-\frac{\Rea z+\eta\Ima z}{1+\eta^2}
 =\frac{(A\eta-B)^2}{1+\eta^2}
 \qquad(\eta\ge0),
\end{equation}
and hence
\begin{equation}\label{eq:unit-sup-id}
 \sup_{\eta\ge0}\frac{\Rea z+\eta\Ima z}{1+\eta^2}=A^2.
\end{equation}
Applying \eqref{eq:unit-quad-id} to $z=s^2=(a+ib)^2$ and dividing by
$a^2$ gives
\begin{equation}\label{eq:unit-s2-id}
 1-\frac{\Rea(s^2)+\eta\Ima(s^2)}{a^2(1+\eta^2)}
 =\frac{(\eta-\varrho)^2}{1+\eta^2}.
\end{equation}
For $|\theta|\le\pi/2$ and $\psi=1-\cos\theta$, combining
\eqref{eq:unit-s2-id} with \eqref{eq:eff-unit-A-gap} (using
$\sin^2\theta=\psi(2-\psi)$) gives the first inequality in
\eqref{eq:unit-angular-decay}. For $0<\psi\le1$, the second inequality
in \eqref{eq:unit-angular-decay} follows from
\[
 \frac{(\eta-\varrho)^2}{\psi}+r\kappa(2-\psi)+2\kappa\eta
 \ge\eta^2-2(\varrho-\kappa)\eta+\varrho^2+r\kappa
 \ge(1-\varrho+\kappa)(1+\eta^2)
 \ge\gamma(1+\eta^2),
\]
using $\varrho^2+r\kappa\ge1$, $\varrho-\kappa\ge0$, and
$2\eta\le1+\eta^2$; the case $\psi=0$ is immediate. For
$|\theta|<\pi/2$, \eqref{eq:unit-log-integral} gives
$\Ima g_r(\theta)>0$, so $\Rea\sqrt{g_r(\theta)}>0$ and
$\Ima\sqrt{g_r(\theta)}>0$. Applying \eqref{eq:unit-sup-id} to
$z=g_r(\theta)$ and using \eqref{eq:unit-angular-decay} yields
\[
 \bigl(\Rea\sqrt{g_r(\theta)}\bigr)^2
 =\sup_{\eta\ge0}
   \frac{\Rea g_r(\theta)+\eta\Ima g_r(\theta)}{1+\eta^2}
 \le a^2(1-\gamma\psi).
\]
Taking square roots gives \eqref{eq:eff-unit-sqrt-gap} for
$|\theta|<\pi/2$, and taking the limit as $|\theta|\to(\pi/2)^-$
extends it to $\theta=\pm\pi/2$.
Finally, since \eqref{eq:eff-unit-A-gap} implies
$\Rea g_r(\theta)\le\Rea(s^2)$ and $\Ima g_r(\theta)\le\Ima(s^2)$
for all $\theta\in\mathbb R$, \eqref{eq:unit-s2-id} yields
\eqref{eq:unit-contour-decay}. For $\pi/2<|\theta|\le\pi$,
\eqref{eq:unit-log-integral} shows that $\Ima g_r(\theta)\le0$;
together with $\Rea g_r(\theta)\le\Rea(s^2)\le r^2/2\le ra^2/2$, this
implies that the left-hand side of \eqref{eq:unit-contour-decay} is
at least $1-r/2$.
\end{proof}

For an integer $n\ge1$ and $0<r\le1$, let
\[
 \mathcal Q
 =\frac1{2\pi}\int_{-\pi}^{\pi}
   \sqrt{\frac{g_r(\theta)}{ n}}
   I_1\bigl(2\sqrt{ng_r(\theta)}\bigr)\,d\theta.
\]
By \eqref{eq:prelim-bessel-series}, the product in the integrand
extends to an entire function of $g_r(\theta)$.
For $t>0$ and $\gamma\in\{0.99,0.55\}$, with $C(\gamma)$ as above, define
\[
 E(t,\gamma)
 =\frac{C(\gamma)}{t}
  +\left(\frac{11.1}{\gamma}+2.4\right)
     t^{-1/2}e^{-0.3\gamma t}
  +19t^2e^{-t}.
\]
\begin{lemma}\label{lem:effective-unit-model}
Let $n\ge1$ be an integer and let $r$ satisfy
$0<r\le1$ and $\sqrt{nr}\ge4\pi$. Then
\begin{equation}\label{eq:left-model-error}
 |\mathcal Q-M_n(r)\cos\Psi_n(r)|
 \le M_n(r)E(a\sqrt n,\gamma).
\end{equation}
\end{lemma}
\begin{proof}
Let $\omega=a\sqrt n$. By \eqref{eq:parameters} and
\eqref{eq:unit-basic-bounds}, we have $|s|<\sqrt2\,a$ and
$a|s|\ge4r/3$, whence
\[
 \frac{r}{a|s|}\le0.75,
 \qquad
 \frac{2r}{a|s|\sqrt{1+r^2}}\le1.5.
\]
Since $g_r(\theta+\pi)=\overline{g_r(\theta)}$, we have
\begin{equation}\label{eq:unit-Q-half}
 \mathcal Q=\Rea\left(\frac1\pi\int_{-\pi/2}^{\pi/2}
  \frac{\sqrt{g_r(\theta)}}{\sqrt n}
  I_1\bigl(2\sqrt{ng_r(\theta)}\bigr)\,d\theta\right).
\end{equation}
For $|\theta|\le\pi/2$, write $\psi=1-\cos\theta$ and
$z=2\sqrt{ng_r(\theta)}$, with $\sqrt{g_r(\pm\pi/2)}$ understood as
the limit from $|\theta|<\pi/2$ as in
Lemma~\ref{lem:effective-unit-bounds}; thus $\Rea z\ge0$.

We first expand $g_r(\theta)$ on $0\le\psi\le0.3$. Subtracting
\eqref{eq:unit-log-integral} from its value at $\theta=0$ gives
\[
 s^2-g_r(\theta)
 =\int_0^r-\log\left(1-\frac{2iy\psi}{(1+iy)^2}\right)
       \frac{dy}{y}.
\]
Since $|2iy/(1+iy)^2|=2y/(1+y^2)\le1$ for $0<y\le r$, we have
$\int_0^r|2iy/(1+iy)^2|^j\,\frac{dy}{y}\le\int_0^r\frac{2\,dy}{1+y^2}\le2r$
for every integer $j\ge1$. Write $w=g_r(\theta)/s^2-1$. Expanding the
logarithm in powers of $2iy\psi/(1+iy)^2$ for $0\le\psi\le0.3$ shows
that the linear term of $s^2-g_r(\theta)$ in $\psi$ has coefficient
$\int_0^r\frac{2i}{(1+iy)^2}\,dy=\frac{2ir}{1+ir}$, and using
$|s|^2\ge16r/9$ gives
\[
 |w|
 \le\frac{2r}{|s|^2}[-\log(1-\psi)]
 \le\frac98[-\log(1-\psi)]\le1.34\psi\le0.402,
\]
\[
 \left|s^2-g_r(\theta)-\frac{2ir\psi}{1+ir}\right|
 \le2r\sum_{j\ge2}\frac{\psi^j}{j}
 \le\frac{r\psi^2}{1-\psi}.
\]
Taylor's formula for $(1+w)^{1/2}$ and $(1+w)^{1/4}$, together with
$r/(a|s|)\le0.75$ and $|s|<\sqrt2\,a$, therefore yields
\begin{equation}\label{eq:unit-taylor}
 \left|\frac{2(\sqrt{g_r(\theta)}-s)}{a}+\frac{2ir\psi}{(1+ir)as}\right|\le2.46\psi^2,
 \qquad
 \left|\left(\frac{\sqrt{g_r(\theta)}}{s}\right)^{1/2}(1-\psi/2)^{-1/2}-1\right|\le0.87\psi,
\end{equation}
where the constants follow from
\[
 \frac{0.75}{0.7}
 +\frac{\sqrt2(1.34)^2}{4(1-0.402)^{3/2}}<2.46,
 \qquad
 \frac{1.34}{4(1-0.402)^{3/4}\sqrt{0.85}}
 +\frac{1/\sqrt{0.85}-1}{0.3}<0.87.
\]
Moreover, since $|\sqrt{1+w}-1|\le1-\sqrt{1-|w|}=|w|/(1+\sqrt{1-|w|})$,
$|w|\le0.402$, and $|s|<\sqrt2\,a$, we obtain for $0\le\psi\le0.3$ that
\[
 |z|\le2\sqrt{2(1+0.402)}\,\omega<3.35\omega,
 \qquad
 \Rea z\ge2\left(1-\frac{\sqrt2\,0.402}{1+\sqrt{0.598}}\right)\omega>1.35\omega.
\]
By \eqref{eq:prelim-bessel-integrals},
\begin{equation}\label{eq:unit-I1-rep}
 I_1(z)=\frac{ze^z\sqrt2}{\pi}\int_0^2 e^{-zy}\sqrt y\sqrt{1-y/2}\,dy,
 \qquad
 \frac{z\sqrt2}{\pi}\int_0^\infty e^{-zy}\sqrt y\,dy=\frac1{\sqrt{2\pi z}}
 \quad(\Rea z>0).
\end{equation}
Since $|1-\sqrt{1-y/2}|\le y/2$ on $[0,2]$ and $1\le y/2$ for $y\ge2$,
it follows that
\begin{equation}\label{eq:unit-I1-asymp}
 \left|I_1(z)-\frac{e^z}{\sqrt{2\pi z}}\right|
 \le\frac{e^{\Rea z}}{\sqrt{2\pi|z|}}\,\frac{3|z|^{3/2}}{4(\Rea z)^{5/2}}
 \qquad(\Rea z>0).
\end{equation}
For $0\le\psi\le0.3$, the last factor is less than $2.2/\omega$,
since $3(3.35)^{3/2}/(4(1.35)^{5/2})<2.2$.
On $\{\theta\in[-\pi/2,\pi/2]:0\le\psi\le0.3\}$, since $g_r$ is even,
the substitution $\psi=1-\cos\theta$ (with $\theta=\arccos(1-\psi)$ on
the positive half) gives $2\,d\theta=\sqrt2\,\psi^{-1/2}(1-\psi/2)^{-1/2}\,d\psi$
from the two symmetric halves and turns the integral of
$\frac1\pi\frac{\sqrt{g_r(\theta)}}{\sqrt n}\frac{e^z}{\sqrt{2\pi z}}\,d\theta$
into the integral over $[0,0.3]$ of
\begin{equation}\label{eq:unit-model-integrand}
 \frac{s^{1/2}e^{2\sqrt n s}}{\pi\sqrt{2\pi}n^{3/4}}
  \left(\frac{\sqrt{g_r(\theta)}}{s}\right)^{1/2}(1-\psi/2)^{-1/2}
  e^{2\sqrt n(\sqrt{g_r(\theta)}-s)}\psi^{-1/2}\,d\psi.
\end{equation}
Comparing right derivatives in \eqref{eq:eff-unit-sqrt-gap} at $\psi=0$
(where equality holds) via \eqref{eq:unit-taylor} gives
$\Rea\bigl(2ir/((1+ir)as)\bigr)\ge\gamma>0$. Hence the integral below
converges and, by \eqref{eq:amplitude-phase}, evaluates to
\[
 \frac{s^{1/2}e^{2\sqrt n s}}{\pi\sqrt{2\pi}n^{3/4}}
 \int_0^\infty
 e^{-\frac{2\omega ir\psi}{(1+ir)as}}\psi^{-1/2}\,d\psi
 =\frac{s e^{2\sqrt n s}}{2\pi n\sqrt{ir/(1+ir)}}
 =M_n(r)e^{i\Psi_n(r)}.
\]
By \eqref{eq:unit-Q-half} and the preceding identity, splitting at
$\psi=0.3$ gives
\begin{equation}\label{eq:unit-Q-split}
 \begin{aligned}
 |\mathcal Q-M_n(r)\cos\Psi_n(r)|
 &\le\left|\frac1\pi\int_{0\le\psi\le0.3}
   \frac{\sqrt{g_r(\theta)}}{\sqrt n}I_1(z)\,d\theta
   -\frac{s^{1/2}e^{2\sqrt n s}}{\pi\sqrt{2\pi}n^{3/4}}
    \int_0^{0.3}e^{-\frac{2\omega ir\psi}{(1+ir)as}}\psi^{-1/2}\,d\psi\right|\\
 &\quad+\frac1\pi\int_{\psi>0.3}
   \left|\frac{\sqrt{g_r(\theta)}}{\sqrt n}I_1(z)\right|\,d\theta
   +\frac{|s|^{1/2}e^{2\omega}}{\pi\sqrt{2\pi}n^{3/4}}
    \left|\int_{0.3}^\infty e^{-\frac{2\omega ir\psi}{(1+ir)as}}\psi^{-1/2}\,d\psi\right|,
 \end{aligned}
\end{equation}
where the integrals in $\theta$ are taken over the indicated subsets
of $[-\pi/2,\pi/2]$. We estimate these three terms in order.

For the first term in \eqref{eq:unit-Q-split}, on $0\le\psi\le0.3$,
\eqref{eq:eff-unit-sqrt-gap} and
$\Rea\bigl(2ir/((1+ir)as)\bigr)\ge\gamma$ give
\[
 \max\left\{
 \left|e^{2\sqrt n(\sqrt{g_r(\theta)}-s)}\right|,
 \left|e^{-\frac{2\omega ir\psi}{(1+ir)as}}\right|
 \right\}
 \le e^{-\gamma\omega\psi}.
\]
Together with \eqref{eq:unit-taylor} and
$|e^A-e^B|\le|A-B|\max\{e^{\Rea A},e^{\Rea B}\}$, this yields
\[
 \left|
 e^{2\sqrt n(\sqrt{g_r(\theta)}-s)}
 -e^{-\frac{2\omega ir\psi}{(1+ir)as}}
 \right|
 \le2.46\omega\psi^2e^{-\gamma\omega\psi}.
\]
By \eqref{eq:amplitude-phase} and $2r/(a|s|\sqrt{1+r^2})\le1.5$,
\begin{equation}\label{eq:unit-prefactor}
 \frac{|s|^{1/2}e^{2\omega}}
      {\pi\sqrt{2\pi}n^{3/4}M_n(r)}
 =\sqrt{\frac{2r\sqrt n}{\pi|s|\sqrt{1+r^2}}}
 \le\sqrt{\frac{1.5\omega}{\pi}}.
\end{equation}
In the first term of \eqref{eq:unit-Q-split}, we replace $I_1(z)$ by
$e^z/\sqrt{2\pi z}$, then replace
$(\sqrt{g_r(\theta)}/s)^{1/2}(1-\psi/2)^{-1/2}$ by $1$ in
\eqref{eq:unit-model-integrand}, and finally replace
$e^{2\sqrt n(\sqrt{g_r(\theta)}-s)}$ by
$e^{-\frac{2\omega ir\psi}{(1+ir)as}}$. Since
$|(\sqrt{g_r(\theta)}/s)^{1/2}(1-\psi/2)^{-1/2}|\le1+0.87\cdot0.3$
by \eqref{eq:unit-taylor} and $2.2(1+0.87\cdot0.3)=2.7742$, it
follows from \eqref{eq:unit-I1-asymp} and \eqref{eq:unit-prefactor}
that the first term in \eqref{eq:unit-Q-split}, divided by $M_n(r)$,
is at most
\[
 \sqrt{\frac{1.5\omega}{\pi}}
 \int_0^{0.3}e^{-\gamma\omega\psi}\psi^{-1/2}
 \left(2.46\omega\psi^2+0.87\psi+
                       \frac{2.7742}\omega\right)\,d\psi
 \le\frac{C(\gamma)}\omega,
\]
where the integrals are bounded by
\[
 \int_0^\infty e^{-\gamma\omega\psi}\psi^{j-1/2}\,d\psi
 =\frac{\sqrt\pi(2j)!}{4^j j!(\gamma\omega)^{j+1/2}}
 \qquad(j=0,1,2).
\]

For the second term in \eqref{eq:unit-Q-split}, we split
$\{\theta\in[-\pi/2,\pi/2]:\psi>0.3\}$ according to whether
$\Rea z\ge\omega$ or $\Rea z<\omega$.
When $\Rea z\ge\omega$, \eqref{eq:unit-I1-rep} gives
$|I_1(z)|\le|z|e^{\Rea z}/(\sqrt{2\pi}(\Rea z)^{3/2})$.
By \eqref{eq:prelim-gr},
$|g_r(\theta)|\le2\sum_{j\ge1}r^j/j^2\le(\pi^2/3)r$.
Together with $|s|\ge(4/3)\sqrt r$, $(1+r^2)^{1/4}\ge1$, $a^2\ge r$,
and \eqref{eq:eff-unit-sqrt-gap}, we obtain
\[
 \frac1{\pi M_n(r)}
 \left|
 \frac{\sqrt{g_r(\theta)}}{\sqrt n}I_1(z)
 \right|
 \le\frac{\pi^2}{\sqrt{2\pi}}\sqrt\omega\,
       e^{-\gamma\omega\psi}
 <3.938\sqrt\omega\,e^{-\gamma\omega\psi}.
\]
On $\{\theta\in(0,\pi/2]:\psi>0.3\}$, we have $\cos\theta<0.7$ and
hence $\sin\theta>\sqrt{0.51}$. By evenness of $g_r$, the contribution
of $\{\theta\in[-\pi/2,\pi/2]:\psi>0.3,\ \Rea z\ge\omega\}$ to the
second term of \eqref{eq:unit-Q-split}, divided by $M_n(r)$, is at most
$\frac{7.876}{\gamma\sqrt{0.51\omega}}e^{-0.3\gamma\omega}
<\frac{11.1}{\gamma\sqrt\omega}e^{-0.3\gamma\omega}$.
On $\{\theta\in[-\pi/2,\pi/2]:\psi>0.3,\ \Rea z<\omega\}$, since
$\Rea z\ge0$, \eqref{eq:prelim-bessel-integrals} gives
$|I_1(z)|\le|z|e^{\Rea z}/2$; using $M_n(r)\ge2e^{2\omega}/(3\pi n)$
from \eqref{eq:unit-scale-bounds} and the fact that this subset has
measure at most $\pi$ bounds its contribution to the second term of
\eqref{eq:unit-Q-split}, divided by $M_n(r)$, by
$(\pi^3/2)\omega^2e^{-\omega}<19\omega^2e^{-\omega}$.

Finally, by \eqref{eq:unit-prefactor},
$\Rea\bigl(2ir/((1+ir)as)\bigr)\ge\gamma$, and $\gamma\ge0.55$, the
third term in \eqref{eq:unit-Q-split}, divided by $M_n(r)$, is at most
\[
 \sqrt{\frac{1.5\omega}{\pi}}
 \int_{0.3}^\infty
 e^{-\gamma\omega\psi}\psi^{-1/2}\,d\psi
 \le
 \sqrt{\frac{1.5\omega}{\pi}}\,
 \frac{e^{-0.3\gamma\omega}}{\sqrt{0.3}\,\gamma\omega}
 =\sqrt{\frac5\pi}\,\gamma^{-1}\omega^{-1/2}
 e^{-0.3\gamma\omega}
 <2.4\omega^{-1/2}e^{-0.3\gamma\omega}.
\]
Combining the estimates for the three terms in \eqref{eq:unit-Q-split}
proves \eqref{eq:left-model-error}.
\end{proof}

We now compare $D_n(-r^2)$ with $\mathcal Q$ and combine the resulting
bound with Lemmas~\ref{lem:effective-unit-bounds}
and~\ref{lem:effective-unit-model} to establish the uniform
approximation for $16\pi^2/n\le r\le1$.

\begin{proposition}\label{lem:effective-unit-left}
If $n\ge6400$ is an integer, $0<r\le1$, and $\sqrt{nr}\ge4\pi$, then
\[
 |D_n(-r^2)-M_n(r)\cos\Psi_n(r)|
 <\frac9{10}M_n(r).
\]
\end{proposition}
\begin{proof}
Let $\omega=a\sqrt n$, $u=a/\sqrt n$, $\varrho=b/a$,
$\kappa=r/((1+r^2)a^2)$, and $\zeta=u+iv=u(1+i\eta)$.
By Lemma~\ref{lem:effective-unit-model},
\[
 \frac{|D_n(-r^2)-M_n(r)\cos\Psi_n(r)|}{M_n(r)}
 \le\frac{|D_n(-r^2)-\mathcal Q|}{M_n(r)}+E(\omega,\gamma).
\]
Throughout the proof we use
\eqref{eq:unit-basic-bounds}--\eqref{eq:unit-contour-decay} from
Lemma~\ref{lem:effective-unit-bounds}. In particular,
\eqref{eq:unit-basic-bounds} and $n\ge6400$ give
$u\le\sqrt{1.282/6400}<0.01416$. Since $u\omega=a^2$,
\[
 1-\frac{\Rea(g_r(\theta)/\zeta)}{\omega}
 =1-\frac{\Rea g_r(\theta)+\eta\Ima g_r(\theta)}{a^2(1+\eta^2)};
\]
in this proof we refer to this quantity as the \emph{gap}.
By \eqref{eq:prelim-gr}, we have $\langle g_r\rangle=0$, and hence
\[
 \left\langle e^{g_r(\theta)/\zeta}-1\right\rangle
 =\left\langle
   \sum_{j\ge2}\frac{g_r(\theta)^j}{j!\zeta^j}
  \right\rangle.
\]
The boundedness of $g_r$ and the inequality $|\zeta|\ge u$ give
\[
 \sum_{j\ge2}\frac{|g_r(\theta)|^j}{j!|\zeta|^j}
 =O\bigl((u^2+v^2)^{-1}\bigr),
\]
uniformly in $\theta$. Thus termwise inverse Laplace transformation
(as in the proof of Lemma~\ref{lem:effective-bessel}) is justified
and yields
\[
 \frac1{2\pi}\int_{\mathbb R}e^{n(u+iv)}
   \left\langle e^{g_r(\theta)/(u+iv)}-1\right\rangle\,dv
 =\left\langle
   \sum_{j\ge2}
   \frac{n^{j-1}g_r(\theta)^j}{j!(j-1)!}
  \right\rangle
 =\left\langle
   \frac{\sqrt{g_r(\theta)}}{\sqrt n}
   I_1\bigl(2\sqrt n\,\sqrt{g_r(\theta)}\bigr)
   -g_r(\theta)
  \right\rangle
 =\mathcal Q,
\]
where the second equality uses \eqref{eq:prelim-bessel-series}
and the last uses $\langle g_r\rangle=0$.
Write $\mathcal L=\log\mathcal P_r(\zeta,\theta)$. Expanding the
principal logarithms in \eqref{eq:P-product} as in the proof of
Lemma~\ref{lem:effective-bessel} gives
\begin{equation}\label{eq:left-L-series}
 \mathcal L
 =\sum_{j\ge1}
  \frac{2(-1)^{j+1}(ir)^j\cos(j\theta)}
       {j\,2\sinh(j\zeta/2)},
\end{equation}
which converges absolutely and uniformly in $(v,\theta)\in\mathbb R^2$
for every $0<r\le1$ (since $|2\sinh(j\zeta/2)|\ge ju$) and satisfies
$\langle\mathcal L\rangle=0$.
Since $n\ge1$ and $16u<0.23<\pi$, equations \eqref{eq:ct-integral} and
\eqref{eq:ct-integral2} give
$D_n(-r^2)=(2\pi)^{-1}\int_{-\pi}^{\pi}e^{n\zeta}\langle e^{\mathcal L}-1\rangle\,dv$,
and splitting at $|v|=16u$ yields
\begin{equation}\label{eq:left-three-parts}
 \begin{aligned}
 |D_n(-r^2)-\mathcal Q|
 \le\frac{e^\omega}{2\pi}\Bigg(
 &\int_{|v|\le16u}
   \left|\left\langle e^{\mathcal L}-e^{g_r(\theta)/\zeta}\right\rangle\right|\,dv\\
 &+\int_{16u\le|v|\le\pi}
   \left|\langle e^{\mathcal L}-1\rangle\right|\,dv
  +\int_{|v|\ge16u}
   \left|\left\langle e^{g_r(\theta)/\zeta}-1\right\rangle\right|\,dv
 \Bigg).
 \end{aligned}
\end{equation}
The symmetries
\[
 g_r(\theta+\pi)=\overline{g_r(\theta)},\qquad
 \mathcal P_r(u-iv,\theta+\pi)
 =\overline{\mathcal P_r(u+iv,\theta)}
\]
imply that the moduli of the angular averages in
\eqref{eq:left-three-parts} are even in $v$ (and hence in $\eta$), so
we may restrict the first and third integrals to $v\ge0$ (that is,
$\eta\ge0$) and double the resulting bounds. Using $v=u\eta$,
$nu=\omega$,
$|e^{\mathcal L}-e^{g_r(\theta)/\zeta}|
\le e^{\Rea(g_r(\theta)/\zeta)}(e^{|\mathcal L-g_r(\theta)/\zeta|}-1)$,
and $M_n(r)\ge2e^{2\omega}/(3\pi n)$ from \eqref{eq:unit-scale-bounds},
the first term on the right-hand side of \eqref{eq:left-three-parts},
divided by $M_n(r)$, is at most
\begin{equation}\label{eq:left-major-bound}
 \frac{3\omega}{4\pi}\int_0^{16}\int_{-\pi}^{\pi}
 e^{-\omega(1-\Rea(g_r(\theta)/\zeta)/\omega)}
 \bigl(e^{|\mathcal L-g_r(\theta)/\zeta|}-1\bigr)\,d\theta\,d\eta.
\end{equation}

We first bound the third term on the right-hand side of
\eqref{eq:left-three-parts} uniformly for $0<r\le1$.
From \eqref{eq:parameters}, $\Rea(s^2)=-\frac12\Li_2(-r^2)\le r^2/2$
and $\Ima(s^2)=2\int_0^r(\arctan t)/t\,dt\le2r$. Together with
\eqref{eq:eff-unit-A-gap} and \eqref{eq:unit-basic-bounds}, these give
$\Rea g_r(\theta)\le\Rea(s^2)\le r^2/2\le a^2/2$ and
$\Ima g_r(\theta)\le\Ima(s^2)\le2r\le2a^2$. Hence, for $\eta\ge16$,
\[
 \frac{\Rea(g_r(\theta)/\zeta)}{\omega}
 \le\frac{0.5+2\eta}{1+\eta^2}
 \le\frac{32.5}{257}<0.127,
\]
while Parseval's identity applied to \eqref{eq:prelim-gr} gives
$\langle|g_r(\theta)|^2\rangle=2\sum_{j\ge1}r^{2j}/j^4<2.17r^2$.
Restricting to $v\ge16u$ and doubling, subtracting $g_r(\theta)/\zeta$
inside $\langle e^{g_r(\theta)/\zeta}-1\rangle$ (since
$\langle g_r\rangle=0$), applying
$|e^w-1-w|\le\frac12|w|^2e^{\max\{0,\Rea w\}}$, using
$\int_{16u}^\infty(u^2+v^2)^{-1}\,dv\le1/(16u)$, and dividing by
$M_n(r)\ge2e^{2\omega}/(3\pi n)$ bounds the third term of
\eqref{eq:left-three-parts}, divided by $M_n(r)$, by
\[
 L_{16}(\omega)=0.102\omega^3e^{-0.873\omega},
\]
since $1.085\cdot3/32<0.102$ and
$nr^2/u=\omega^3(r/a^2)^2\le\omega^3$.

To estimate $|\mathcal L-g_r(\theta)/\zeta|$ in
\eqref{eq:left-major-bound}, define
$f_\pm(x)=\log(1+ir e^{\pm i\theta-(1+i\eta)x})$ for $x\ge0$.
For $0<r<1$, the error of the midpoint rule on each interval
$[ju,(j+1)u]$ is an integral of $f_\pm''$ against a kernel of modulus
at most $u^2/8$. Summing over $j\ge0$ and dividing by $u$ gives
\begin{equation}\label{eq:left-midpoint-err}
 \left|\sum_{j\ge0}f_\pm\bigl((j+1/2)u\bigr)
       -\frac1u\int_0^\infty f_\pm(x)\,dx\right|
 \le\frac u8\int_0^\infty|f_\pm''(x)|\,dx,
\end{equation}
where
\begin{equation}\label{eq:left-f-second}
 |f_\pm''(x)|
 =\frac{(1+\eta^2)re^{-x}}
       {|1+ir e^{\pm i\theta-(1+i\eta)x}|^2}.
\end{equation}
Here $\sum_{j\ge0}\bigl(f_+((j+1/2)u)+f_-((j+1/2)u)\bigr)=\mathcal L$,
and termwise integration of the logarithmic series gives
\[
 \frac1u\int_0^\infty\bigl(f_+(x)+f_-(x)\bigr)\,dx
 =\frac{g_r(\theta)}{u(1+i\eta)}
 =\frac{g_r(\theta)}\zeta.
\]
Bounding the denominator in \eqref{eq:left-f-second} below by
$(1-re^{-x})^2$ yields
\[
 \int_0^\infty|f_\pm''(x)|\,dx
 \le(1+\eta^2)\int_0^\infty
       \frac{re^{-x}}{(1-re^{-x})^2}\,dx
 =\frac{(1+\eta^2)r}{1-r},
\]
and adding the two bounds in \eqref{eq:left-midpoint-err} proves
\begin{equation}\label{eq:left-midpoint}
 |\mathcal L-g_r(\theta)/\zeta|\le\frac{u(1+\eta^2)r}{4(1-r)}
 \qquad(0<r<1).
\end{equation}
For $0<r_0\le1$, write $u_{r_0}=\sqrt{r_0(1+0.282r_0)/6400}$, so that
$u\le u_{r_0}$ whenever $r\le r_0$ by \eqref{eq:unit-basic-bounds}.
We now bound the first and second terms of \eqref{eq:left-three-parts},
divided by $M_n(r)$, in two ranges of $r$.

\medskip
\noindent\textbf{Case 1: $0<r\le1/4$.}
Fix $0<r_0\le1/4$ and let
\[
 \delta_{r_0}=\frac{5r_0u_{r_0}}{2(1-r_0)},\qquad
 \Delta_{r_0}=\frac{257r_0u_{r_0}}{4(1-r_0)}.
\]
If $r\le r_0$, then $ur\le r_0u_{r_0}$, so \eqref{eq:left-midpoint}
bounds $|\mathcal L-g_r(\theta)/\zeta|$ by $\delta_{r_0}$ on
$0\le\eta\le3$ and by $\Delta_{r_0}$ on $0\le\eta\le16$.
For $0\le\eta\le3$ and $|\theta|\le\pi/2$, averaging
\eqref{eq:unit-angular-decay} and \eqref{eq:unit-contour-decay}, and
using $\gamma\ge0.55$ and $1-\cos\theta\ge2\theta^2/\pi^2$, gives
\[
 1-\frac{\Rea(g_r(\theta)/\zeta)}{\omega}
 \ge\frac{(\eta-\varrho)^2}{20}+\frac{0.55\theta^2}{\pi^2}.
\]
Integrating $\exp\bigl(-\omega((\eta-\varrho)^2/20+0.55\theta^2/\pi^2)\bigr)$
over $\mathbb R^2$ yields $\pi^2\sqrt{20/0.55}/\omega$.
For $0\le\eta\le3$ and $\pi/2<|\theta|\le\pi$,
Lemma~\ref{lem:effective-unit-bounds} gives a gap of at least
$1-r/2\ge1-r_0/2$.
For $3\le\eta\le16$ and $|\theta|\le\pi$, \eqref{eq:unit-contour-decay}
and $\varrho\le1$ give a gap of at least $(3-1)^2/10=0.4$.
Splitting \eqref{eq:left-major-bound} into these three subregions
therefore bounds the first term of \eqref{eq:left-three-parts},
divided by $M_n(r)$, by
\[
 \mathcal B_{r_0}(\omega)
 =\frac{3\pi}{4}\sqrt{20/0.55}\,(e^{\delta_{r_0}}-1)
 +\frac{9\omega}{4}(e^{\delta_{r_0}}-1)e^{-(1-r_0/2)\omega}
 +\frac{39\omega}{2}(e^{\Delta_{r_0}}-1)e^{-0.4\omega}.
\]

For the second term of \eqref{eq:left-three-parts}, on
$16u\le|v|\le\pi$, separating the $j=1$ term in \eqref{eq:left-L-series}
from the terms with $j\ge2$ and using $r/u=(r/a^2)\omega\le\omega$
gives
\[
 \Rea\mathcal L\le\left(\frac{\pi}{16}+\frac{r}{2(1-r)}\right)\omega;
\]
here $|2\sinh(\zeta/2)|\ge2|v|/\pi\ge32u/\pi$, and the sum over
$j\ge2$ has modulus at most
$\sum_{j\ge2}2r^j/(j^2u)\le r^2/(2(1-r)u)$.
Orthogonality of $\cos(j\theta)$ in $\theta$, followed by integration
in $v$, gives
\[
 \int_{-\pi}^{\pi}\langle|\mathcal L|^2\rangle\,\frac{dv}{2\pi}
 =\sum_{j\ge1}\frac{r^{2j}}{j^2\sinh(ju)}
 \le\frac1u\sum_{j\ge1}\frac{r^{2j}}{j^3}
 \le\frac{r^2}{u(1-r^2)},
\]
where the first equality uses
$(2\pi)^{-1}\int_{-\pi}^{\pi}|2\sinh(j(u+iv)/2)|^{-2}\,dv
=1/(2\sinh(ju))$.
Since $\langle\mathcal L\rangle=0$, subtracting $\mathcal L$ inside
$\langle e^{\mathcal L}-1\rangle$, applying
$|e^w-1-w|\le\frac12|w|^2e^{\max\{0,\Rea w\}}$, and dividing by
$M_n(r)\ge2e^{2\omega}/(3\pi n)$ shows that, for $r\le r_0$, the
second term of \eqref{eq:left-three-parts}, divided by $M_n(r)$, is
at most
\[
 U_{r_0}(\omega)
 =\frac{3\pi\omega^3}{4(1-r_0^2)}
  \exp\!\left[-\left(1-\frac{\pi}{16}
                     -\frac{r_0}{2(1-r_0)}\right)\omega\right].
\]

In the following table, $r_0$ is the upper endpoint of the range and
$\omega_0$ is a lower bound for $\omega$ (from $\omega^2\ge6400r$ and
\eqref{eq:unit-scale-bounds}); the last four columns are strict upper
bounds:
\[
\begin{array}{c|c|c|c|c|c|c}
\text{range of }r&\omega_0&\gamma&\mathcal B_{r_0}(\omega_0)
 &U_{r_0}(\omega_0)&L_{16}(\omega_0)&E(\omega_0,\gamma)\\\hline
(0,0.03]&4\pi&0.99&0.0094&0.234&0.0035&0.602\\
(0.03,0.1]&\sqrt{192}&0.99&0.047&0.200&0.0016&0.516\\
(0.1,0.25]&\sqrt{640}&0.55&0.080&0.0041&0.000001&0.700
\end{array}
\]
All four bounding functions are nonincreasing for $\omega\ge\omega_0$,
since $\omega^{-1}$ decreases and $\omega^q e^{-c\omega}$ decreases
for $\omega>q/c$ whenever $c>0$.
The row sums are less than $0.849$, $0.765$, and $0.785$, each below
$9/10$.

\medskip
\noindent\textbf{Case 2: $1/4<r\le1$.}
Here $\omega\ge\sqrt{nr}>40$ and $\gamma=0.55$.
From \eqref{eq:left-L-series} and \eqref{eq:prelim-gr},
\[
 \mathcal L-\frac{g_r(\theta)}{\zeta}
 =\sum_{j\ge1}
  \frac{2(-1)^{j+1}(ir)^j\cos(j\theta)}{j}
  \left(\frac1{2\sinh(j\zeta/2)}-\frac1{j\zeta}\right).
\]
Splitting this sum at $j=\lfloor\pi/|\zeta|\rfloor$ gives the uniform
bound
\begin{equation}\label{eq:eff-unit-global-midpoint}
 |\mathcal L-g_r(\theta)/\zeta|\le 0.51+0.69\bigl(1+\sqrt{1+\eta^2}\bigr),
 \qquad |\eta|\le16.
\end{equation}
Indeed, for $|z|\le\pi$, we have
$|2\sinh(z/2)|^2=4\sinh^2(\Rea z/2)+4\sin^2(\Ima z/2)\ge4|z|^2/\pi^2$,
while
\[
 |2\sinh(z/2)-z|
 \le\frac{2\sinh(\pi/2)-\pi}{\pi^3}|z|^3,
\]
and consequently
\[
 \left|\frac1{2\sinh(z/2)}-\frac1z\right|
 \le\frac{2\sinh(\pi/2)-\pi}{2\pi^2}|z|
 <0.08|z|.
\]
Thus the terms with $j\le\lfloor\pi/|\zeta|\rfloor$ contribute at most
$0.16\pi<0.51$. For the remaining terms, $|2\sinh(j\zeta/2)|\ge ju$
gives $|1/(2\sinh(j\zeta/2))-1/(j\zeta)|\le(1+\sqrt{1+\eta^2})/(j|\zeta|)$;
since $\sum_{j>\lfloor\pi/|\zeta|\rfloor}j^{-2}\le|\zeta|/(\pi-|\zeta|)$
and $|\zeta|\le\sqrt{257}\,u<0.228$ with $2/(\pi-0.228)<0.69$, the
tail contributes at most $0.69(1+\sqrt{1+\eta^2})$, proving
\eqref{eq:eff-unit-global-midpoint}.

On the rectangle $1/4\le\eta\le2$, $|\theta|\le\pi/3$, we obtain a
sharper bound. For $0\le x\le\pi/\eta$, the argument
$\pi/2\pm\theta-\eta x$ of $ir e^{\pm i\theta-(1+i\eta)x}$ lies in
$[-5\pi/6,5\pi/6]$, so $|1+ir e^{\pm i\theta-(1+i\eta)x}|\ge1/2$.
For $x\ge\pi/\eta$, this modulus is at least $1-e^{-x}$. Since
$\eta\le2$, \eqref{eq:left-f-second} gives, for $0<r<1$,
\[
 \int_0^\infty|f_\pm''(x)|\,dx
 \le(1+\eta^2)r\left(4+\frac1{e^{\pi/2}-1}\right)
 <4.263(1+\eta^2)r.
\]
Since $1+\eta^2\le5$, adding the two bounds in
\eqref{eq:left-midpoint-err} yields
$|\mathcal L-g_r(\theta)/\zeta|\le\frac54(4.263)ur<5.33ur$ for
$0<r<1$. For fixed $n$, $\eta$, and $\theta$, we have
$u\ge1/(2\sqrt n)$ on $[1/4,1]$, so \eqref{eq:left-L-series} and
\eqref{eq:prelim-gr} converge uniformly in $r\in[1/4,1]$; thus both
$\mathcal L$ and $g_r(\theta)/\zeta$ are continuous in $r\in[1/4,1]$,
and the bound $|\mathcal L-g_r(\theta)/\zeta|\le\frac54(4.263)ur<5.33ur$
extends to $r=1$ by continuity.
By \eqref{eq:unit-angular-decay}, the gap on this rectangle is at
least $0.2(\eta-\varrho)^2+0.082\theta^2$, using $1+\eta^2\le5$,
$\kappa>0.39$, $2\eta/(1+\eta^2)\ge8/17$, and
$1-\cos\theta\ge0.45\theta^2$ for $|\theta|\le\pi/3$.
Integrating $\exp\bigl(-\omega(0.2(\eta-\varrho)^2+0.082\theta^2)\bigr)$
over $\mathbb R^2$ gives less than $24.54/\omega$. Hence, for
$r\le r_0$, the contribution of this rectangle to
\eqref{eq:left-major-bound} is less than
\[
 J_{r_0}=\frac{3(24.54)}{4\pi}\bigl(e^{5.33r_0u_{r_0}}-1\bigr).
\]

On the complement of this rectangle in $0\le\eta\le16$,
$|\theta|\le\pi$, we divide the domain into five subregions.
Since $\varrho\in[0.718,1]$ by \eqref{eq:unit-basic-bounds} and
$\gamma=0.55$, equations \eqref{eq:unit-angular-decay} and
\eqref{eq:unit-contour-decay} (together with the lower bound
$1-r/2\ge1/2$ on $\pi/2<|\theta|\le\pi$ from
Lemma~\ref{lem:effective-unit-bounds}) give a gap of at least $0.206$
for $0\le\eta<1/4$ and at least $0.275$ for $1/4\le\eta\le2$,
$|\theta|>\pi/3$. By \eqref{eq:eff-unit-global-midpoint}, the
corresponding upper bounds for $e^{|\mathcal L-g_r(\theta)/\zeta|}-1$
are $6$ and $15$.
For $\eta\ge2$, \eqref{eq:unit-contour-decay} and $\varrho\le1$ give
a gap of at least $(\eta-1)^2/(1+\eta^2)$.
On $2\le\eta\le4$, we use
\[
 \frac{(\eta-1)^2}{1+\eta^2}-\frac\eta{10}
 =\frac{(\eta-2)(-\eta^2+8\eta-5)}{10(1+\eta^2)}\ge0,
 \qquad
 e^{|\mathcal L-g_r(\theta)/\zeta|}-1
 <e^{1.38+0.69\eta},
\]
where the second inequality follows from
\eqref{eq:eff-unit-global-midpoint} and
$\sqrt{1+\eta^2}\le\eta+1/(2\eta)$.
On $4\le\eta\le8$ and $8\le\eta\le16$, the gap is at least $9/17$ and
$49/65$, respectively, and \eqref{eq:eff-unit-global-midpoint} bounds
$e^{|\mathcal L-g_r(\theta)/\zeta|}-1$ by $900$ and $220000$.
Bounding the angular measure in each subregion by $2\pi$ in
\eqref{eq:left-major-bound} (giving prefactor
$(3\omega/(4\pi))\cdot2\pi=3\omega/2$) and using
\[
 \int_2^4 e^{1.38-(0.1\omega-0.69)\eta}\,d\eta
 \le\frac{e^{2.76-0.2\omega}}{0.1\omega-0.69}
\]
bounds the contribution of the complement to
\eqref{eq:left-major-bound} by
\[
 E_*(\omega)=\frac{3\omega}{2}\biggl(
 1.5e^{-0.206\omega}
 +\frac{e^{2.76-0.2\omega}}{0.1\omega-0.69}
 +3600e^{-9\omega/17}
 +1760000e^{-49\omega/65}
 +\frac{105}{4}e^{-0.275\omega}\biggr).
\]
All five summands decrease for $\omega\ge40$; for the second term
also use that $\omega/(0.1\omega-0.69)$ decreases.

We now bound the second term of \eqref{eq:left-three-parts} (on
$16u\le|v|\le\pi$). Since
\[
|2\sinh(j\zeta/2)|\ge\max\{ju,2|\sin(jv/2)|\},
\]
 dividing the modulus
of the $j$th term in \eqref{eq:left-L-series} by $r/u$ gives at most
$r^{j-1}\min\{2/j^2,\,u/(j|\sin(jv/2)|)\}$, omitting the second entry
when $\sin(jv/2)=0$.
Take arcs of radius $0.07$ around the angles $2\pi h/k$ with
$1\le k\le7$, $0\le h<k$, and $\gcd(h,k)=1$, measuring distances
modulo $2\pi$. These arcs are disjoint because their centres are
separated by at least $2\pi/42>0.14$. Outside all these arcs, the
distance to every angle $2\pi h/j$ is at least $0.07$ for
$1\le j\le7$. Since $\sin t\ge t-t^3/6+t^5/120-t^7/5040\ge0.99t$ for
$0\le t\le0.245$, we have $|\sin(jv/2)|\ge\sin(0.035j)\ge0.99(0.035j)$
for $1\le j\le7$.
For $0<r_0<1$, let
\[
 T_{r_0}=\min\!\left\{\frac27,\frac{r_0^7}{32(1-r_0)}\right\},
 \qquad Z_{r_0}=\frac{2r_0^{32}}{33^2(1-r_0^2)},
\]
and set $T_1=2/7$, $Z_1=1/32$. After division by $r/u$, $T_{r_0}$ and
$Z_{r_0}$ bound the sums over $j\ge8$ and over odd $j\ge33$ in
\eqref{eq:left-L-series}, respectively, whenever $r\le r_0$. The
geometric bounds follow by replacing $j^2$ by $8^2$ or $33^2$; the
bounds at $r_0=1$ follow from
\[
 2\sum_{j\ge8}\frac1{j^2}\le2\int_7^\infty\frac{dt}{t^2},
 \qquad
 2\sum_{\substack{j\ge33\\j\text{ odd}}}\frac1{j^2}
 \le\int_{32}^\infty\frac{dt}{t^2},
\]
the second inequality using convexity on each interval $[j-1,j+1]$.

For $r\le r_0$, outside all the arcs above, $|\mathcal L|/(r/u)$ is
bounded by
\[
 O_{r_0}(u)=\frac{2u}{0.99(0.07)}
            \sum_{j=1}^7\frac{r_0^{j-1}}{j^2}+T_{r_0}.
\]
On the arc centred at $0$, whenever $16u\le|v|\le0.07$, applying
$\sin t\ge0.99t$ to $t=j|v|/2\le0.245$ gives
$|\sin(jv/2)|\ge0.99j|v|/2\ge0.99\cdot8ju$ for $1\le j\le7$, whence
$|\mathcal L|/(r/u)<0.192+2/7<0.478$.
On an arc centred at $2\pi h/k$ with $3\le k\le7$, the terms with
$j\le7$ and $k\mid j$ contribute at most $2(r_0^2/9+r_0^5/36)$ to
$|\mathcal L|/(r/u)$, while for $j\le7$ with $k\nmid j$ we have
$|\sin(jv/2)|\ge(2\pi/7-0.49)/\pi>0.129$; hence $|\mathcal L|/(r/u)$
is at most
\[
 \mathcal H_{r_0}(u)
 =2\left(\frac{r_0^2}{9}+\frac{r_0^5}{36}\right)
 +\frac{u}{0.129}\sum_{j=1}^7\frac{r_0^{j-1}}{j}+T_{r_0}.
\]
On the remaining arc ($k=2$), centred at $\pi$, summing over even $j$
contributes $\Li_2(r^2)/(2r)$ (which increases with $r$), while for
odd $j\le31$ we have $|\sin(jv/2)|\ge\cos(1.085)>0.466$. Thus, for
$r\le r_0\le1$, $|\mathcal L|/(r/u)$ on this arc is at most
\[
 V_{r_0}(u)=\frac{\Li_2(r_0^2)}{2r_0}
  +\frac{u}{0.466}\sum_{\substack{1\le j\le31\\j\text{ odd}}}
      \frac{r_0^{j-1}}{j}+Z_{r_0}.
\]
For the numerical evaluations below, use $\Li_2(1)=\pi^2/6$ and
\[
 \Li_2(w)\le\sum_{j=1}^8\frac{w^j}{j^2}
                    +\frac{w^9}{81(1-w)}\qquad(0<w<1).
\]

Since $(r/u)/\omega=r/a^2$, multiplying the preceding bounds for
$|\mathcal L|/(r/u)$ by $r/a^2$ bounds $\Rea\mathcal L/\omega$. On
$[0.9,1]$, to keep the resulting bound sufficiently small, we sharpen
$a^2\ge r$ to $a^2>1.12r$. From the alternating series for
$s^2=-2\Li_2(-ir)$ in \eqref{eq:parameters},
\[
 \frac{\Rea(s^2)+\Ima(s^2)}{2r}
 \ge1+\frac r4-\frac{r^2}{9}-\frac{r^3}{16}
     +r^4\left(\frac1{25}+\frac r{36}-\frac{r^2}{49}
                                      -\frac{r^3}{64}\right).
\]
The expression in parentheses decreases on $[0.9,1]$ and exceeds
$1/32$ at $r=1$. Let $f_{\mathrm L}(r)=r/4-r^2/9-r^3/16+r^4/32$.
Since
\[
 f_{\mathrm L}''(r)=-\frac29-\frac{3r(1-r)}8<0,\qquad
 f_{\mathrm L}(0.9)=\frac{35181}{320000}
  >f_{\mathrm L}(1)=\frac{31}{288},
\]
we have $(\Rea(s^2)+\Ima(s^2))/(2r)\ge1+31/288$.
The same alternating series gives $\Rea(s^2)\ge r^2/2-r^4/8\ge3r^2/8$,
$\Ima(s^2)\le2r$, and
$|s|^2\le\sqrt{r^4/4+4r^2}\le\sqrt{17}\,r/2<5r/2$.
Therefore
\[
 \frac{a^2}{r}
 =\frac{\Rea(s^2)+\Ima(s^2)}{2r}
   +\frac{(\Rea(s^2))^2}{2r(|s|^2+\Ima(s^2))}
 >1+\frac{31}{288}+\frac{r^2}{64}
 \ge1+\frac{31}{288}+\frac{81}{6400}>1.12.
\]

In the next table, $r_0$ is the upper endpoint of the range, $\beta$
is a lower bound for $a^2/r$, and $\omega_0$ is a lower bound for
$\omega$ (from $\omega^2=na^2\ge6400\beta r$). The next three columns
are strict upper bounds, and $\mu$ is a strict upper bound for
$\Rea\mathcal L/\omega$ on $16u\le|v|\le\pi$, covering these three
columns as well as $0.478/\beta$ for the arc centred at $0$:
\[
\begin{array}{c|c|c|c|c|c|c}
\text{range of }r&\omega_0&\beta&O_{r_0}(u_{r_0})/\beta
 &\mathcal H_{r_0}(u_{r_0})/\beta&V_{r_0}(u_{r_0})/\beta&\mu\\\hline
(0.25,0.6]&40&1&0.369&0.210&0.359&0.479\\
(0.6,0.9]&\sqrt{3840}&1&0.692&0.590&0.658&0.692\\
(0.9,1]&\sqrt{32256/5}&1.12&0.807&0.758&0.827&0.827
\end{array}
\]
Thus $\Rea\mathcal L<\mu\omega$ on $16u\le|v|\le\pi$, so
$|\langle e^{\mathcal L}-1\rangle|\le e^{\mu\omega}+1<2e^{\mu\omega}$.
Since $a^2\ge r>1/4$ gives $n<4\omega^2$, dividing by
$M_n(r)\ge2e^{2\omega}/(3\pi n)$ bounds the second term of
\eqref{eq:left-three-parts}, divided by $M_n(r)$, by
\[
 U_*(\omega,\mu)=12\pi\omega^2e^{-(1-\mu)\omega}.
\]
In each of the three ranges, $E_*(\omega)$, $E(\omega,0.55)$,
$U_*(\omega,\mu)$, and $L_{16}(\omega)$ decrease for $\omega\ge\omega_0$
(for the last three, since $\omega^{-1}$ decreases and
$\omega^q e^{-c\omega}$ decreases for $\omega>q/c$ whenever $c>0$).
The following table gives strict upper bounds for $J_{r_0}$ and for
$E_*(\omega_0)$, $E(\omega_0,0.55)$, and $U_*(\omega_0,\mu)$:
\[
\begin{array}{c|c|c|c|c}
r_0&J_{r_0}&E_*(\omega_0)&E(\omega_0,0.55)&U_*(\omega_0,\mu)\\\hline
0.6&0.200&0.147&0.404&0.000054\\
0.9&0.386&0.00161&0.258&0.00075\\
1&0.460&0.000040&0.199&0.225
\end{array}
\]
Also, $L_{16}(\omega)\le L_{16}(40)<0.001$. Summing the five bounds
in each range shows that $|D_n(-r^2)-M_n(r)\cos\Psi_n(r)|/M_n(r)$
is less than $0.753$, $0.648$, and $0.886$, respectively.

All terminating decimals in the tables above denote exact rational
numbers. The displayed upper bounds follow by rational squaring, the
enclosure $3.1415926535<\pi<3.1415926536$, and
\[
 \sum_{j=0}^{12}\frac{t^j}{j!}
 \le e^t\le\sum_{j=0}^{12}\frac{t^j}{j!}
       +\frac{t^{13}}{13!(1-t/14)}\qquad(0\le t\le1),
\]
using $e^t=(e^{t/128})^{128}$ for larger positive arguments and
reciprocation for negative arguments. Together with
\eqref{eq:left-model-error} and \eqref{eq:left-three-parts}, Cases~1
and~2 prove the proposition.
\end{proof}

Proposition~\ref{lem:effective-unit-left} covers the intermediate range
for $r\le1$. For $r\ge1$, we use \eqref{eq:theta-identities} to express
the integral in \eqref{eq:ct-integral} as a sum of integrals over
$\mathbb R$; two conjugate terms produce the cosine approximation.
Throughout the remainder of this section, let $r\ge1$ and write
$s=s(r)$, $a=a(r)$, and $b=b(r)$. We retain the functions
$\mathcal E$, $\mathcal P_t$ ($t\in\mathbb R$), and $g_t$
($t\in[-1,1]$) from Section~\ref{sec:preliminaries}.
Setting $z=-ir$ in the second identity of
\eqref{eq:prelim-dilog-identities} (so that
$\log(-z)=\log(ir)=\log r+i\pi/2$) and recalling $s^2=-2\Li_2(-ir)$
from \eqref{eq:parameters} gives
\begin{equation}\label{eq:dilog-inversion}
 s^2=(\log r+i\pi/2)^2+\frac{\pi^2}{3}+2\Li_2(i/r).
\end{equation}
For $m\in\mathbb Z$ and $\chi\in\mathbb R$, define the leading exponent
\[
 H_m(\chi)=(\log r+i\pi(m+1/2))^2+\frac{\pi^2}{3}
           -\chi^2-g_{(-1)^{m+1}r^{-1}}(\chi),
\]
which will arise after applying \eqref{eq:eta-transform} to the
unfolded representation in the next lemma.

\begin{lemma}\label{lem:unfolding}
For $\Rea\zeta>0$ and $r\ge1$,
\begin{equation}\label{eq:unfolded-theta}
 F(-r^2,e^{-\zeta})=
 \frac1{\zeta\mathcal E(\zeta)^2}\sum_{m\in\mathbb Z}
     e^{(\log r+i\pi(m+1/2))^2/\zeta}
    \int_{\mathbb R}
       \frac{e^{-\chi^2/\zeta}}{\mathcal P_{(-1)^{m+1}r^{-1}}(\zeta,\chi)}
       \,d\chi.
\end{equation}
Furthermore, for all $m\in\mathbb Z$ and $\chi\in\mathbb R$,
\begin{equation}\label{eq:all-branches-components}
 \begin{aligned}
 \Rea H_m(\chi)&\le(\log r)^2-\pi^2(m+1/2)^2+\frac{\pi^2}{3}
                    -\Rea g_{r^{-1}}(0)-\frac{\chi^2}{1+r^{-2}},\\
 |\Ima H_m(\chi)|&\le2(\log r)|\pi(m+1/2)|+\Ima g_{r^{-1}}(0).
 \end{aligned}
\end{equation}
Also, $H_0(0)=s^2$ and $H_{-1}(\chi)=\overline{H_0(\chi)}$ (in
particular, $H_{-1}(0)=\overline{s}^{\,2}$).
\end{lemma}
\begin{proof}
Applying the first identity in \eqref{eq:theta-identities} with
$z=ir e^{i\theta}$ and $z=ir e^{-i\theta}$ and multiplying the two
identities gives
\[
 \mathcal P_r(\zeta,\theta)
 =\frac{\vartheta(ir e^{i\theta};\zeta)
             \vartheta(ir e^{-i\theta};\zeta)}
        {\mathcal E(\zeta)^2\mathcal P_{-r^{-1}}(\zeta,\theta)},
\]
where $\mathcal E(\zeta)\neq0$ and
\[
 |\mathcal P_{-r^{-1}}(\zeta,\theta)|
 \ge\prod_{j\ge0}\bigl(1-e^{-(j+1/2)\Rea\zeta}\bigr)^2>0
\]
uniformly for $\theta\in\mathbb R$. Substitute this identity into
\eqref{eq:ct-integral} and apply the second identity in
\eqref{eq:theta-identities} to the two factors in the numerator, with
$w=\log r+i(\pi/2+\theta)$ and $w=\log r+i(\pi/2-\theta)$ and
summation indices $k_1,k_2\in\mathbb Z$, respectively. Setting
$m=k_1+k_2$ and $d=k_1-k_2$, the sum of the two exponents equals
\[
 \frac{\bigl(\log r+i\pi(m+1/2)\bigr)^2-(\theta+\pi d)^2}{\zeta}.
\]
Since $\Rea(1/\zeta)>0$, both series from \eqref{eq:theta-identities}
converge absolutely and uniformly for $\theta\in[-\pi,\pi]$, so their
product converges absolutely and uniformly as a double series over
$(k_1,k_2)\in\mathbb Z^2$. Together with the uniform positive lower
bound for $|\mathcal P_{-r^{-1}}(\zeta,\theta)|$, this justifies
termwise integration and reindexing by $(m,d)$ with $m\in\mathbb Z$
and $d\equiv m\pmod2$. For each fixed $m$, the substitution
$\chi=\theta+\pi d$ sends $[-\pi,\pi]$ to $[\pi(d-1),\pi(d+1)]$, and
as $d$ runs over all integers congruent to $m\pmod2$, these intervals
partition $\mathbb R$ up to endpoints, while
$\mathcal P_{-r^{-1}}(\zeta,\chi-\pi d)
=\mathcal P_{(-1)^{m+1}r^{-1}}(\zeta,\chi)$. Summing over $d$ yields
\eqref{eq:unfolded-theta}.

For $\chi\in\mathbb R$, \eqref{eq:prelim-gr} gives
$g_{-r^{-1}}(\chi)=\overline{g_{r^{-1}}(\chi)}$ and
$g_{r^{-1}}(\chi+\pi)=\overline{g_{r^{-1}}(\chi)}$. Applying
\eqref{eq:eff-unit-A-gap} (with $r^{-1}\in(0,1]$) at $\chi$ and at
$\chi+\pi$ gives $|\Ima g_{r^{-1}}(\chi)|\le\Ima g_{r^{-1}}(0)$ for
all $r\ge1$. For $r>1$, differentiating \eqref{eq:prelim-gr} twice
gives
\[
 \Rea g_{r^{-1}}''(\chi)
 =-\frac{2r^{-2}(r^{-2}+\cos(2\chi))}
         {1+2r^{-2}\cos(2\chi)+r^{-4}}
 \ge-\frac{2}{1+r^2},
\]
since the fraction is maximized at $\cos(2\chi)=1$. Since $g_{r^{-1}}$
is even (so $g_{r^{-1}}'(0)=0$), integrating twice from $0$ to $|\chi|$
yields $\Rea g_{r^{-1}}(\chi)\ge\Rea g_{r^{-1}}(0)-\chi^2/(1+r^2)$.
Since the defining series for $\Li_2$ converges uniformly on the closed
unit disk, $t\mapsto g_t(\chi)$ is continuous on $[0,1]$ for each fixed
$\chi\in\mathbb R$, extending this inequality to $r=1$. Together with
$g_{-r^{-1}}(\chi)=\overline{g_{r^{-1}}(\chi)}$,
$|\Ima g_{r^{-1}}(\chi)|\le\Ima g_{r^{-1}}(0)$, and
$1-(1+r^2)^{-1}=(1+r^{-2})^{-1}$, this proves
\eqref{eq:all-branches-components}.
Finally, since $g_{-r^{-1}}(0)=-2\Li_2(i/r)$, equation
\eqref{eq:dilog-inversion} gives $H_0(0)=s^2$, while
$g_{-r^{-1}}(\chi)=\overline{g_{r^{-1}}(\chi)}$ implies
$H_{-1}(\chi)=\overline{H_0(\chi)}$ and hence
$H_{-1}(0)=\overline{s}^{\,2}$.
\end{proof}

We shall use that both $a$ and $b$ increase strictly with $r$.
Setting $z=-ir$ in the first identity of
\eqref{eq:prelim-dilog-identities} and recalling $s^2=-2\Li_2(-ir)$
from \eqref{eq:parameters} gives
\[
 s^2=\log^2(1+ir)+2\Li_2\!\left(\frac{ir}{1+ir}\right).
\]
Note that
\[
 \int_0^1\frac{ir}{1+ir(1-t)}\,dt=\log(1+ir),\qquad
 \Li_2\!\left(\frac{ir}{1+ir}\right)
 =\int_0^1\frac1y
       \int_0^y\frac{ir}{1+ir(1-t)}\,dt\,dy.
\]
For $0\le t<1$, the number
\[
 \frac{ir}{1+ir(1-t)}
 =\frac{r^2(1-t)+ir}{1+r^2(1-t)^2}
\]
lies in the open first quadrant, and the ratio $1/(r(1-t))$ of its
imaginary part to its real part increases strictly with $t$.
Consequently, for every $0<y<1$, the ratio of the imaginary and real
parts of $\int_0^y\frac{ir}{1+ir(1-t)}\,dt$ is strictly less than
$1/(r(1-y))$, whereas that of $\int_y^1\frac{ir}{1+ir(1-t)}\,dt$ is
strictly greater than $1/(r(1-y))$; hence
\[
 0<\arg\int_0^y\frac{ir}{1+ir(1-t)}\,dt
  <\arg\int_0^1\frac{ir}{1+ir(1-t)}\,dt
  =\arg\log(1+ir)
  <\frac\pi2.
\]
Integrating over $y\in(0,1)$ yields
$0<\arg\Li_2(ir/(1+ir))<\arg\log(1+ir)$.
In the formula for $s^2$ above, the first summand has argument
$2\arg\log(1+ir)<\pi$, whereas the second has a strictly smaller
positive argument. Their sum $s^2$ therefore satisfies
$0<\arg(s^2)<2\arg\log(1+ir)<\pi$, whence
$0<\arg s<\arg\log(1+ir)<\pi/2$.
Differentiating $s^2=-2\Li_2(-ir)$ now gives
\[
 s'(r)=\frac{\log(1+ir)}{rs},\qquad
 0<\arg s'(r)<\frac\pi2,
\]
so $a'(r)>0$ and $b'(r)>0$.
By \eqref{eq:dilog-inversion},
$s-(\log r+i\pi/2)=O\bigl((\log r)^{-1}\bigr)$ as $r\to\infty$;
hence $b(r)\to\pi/2$, and since $b$ increases strictly, $b(r)<\pi/2$
for all $r\ge1$. Consequently, for $r\ge1$,
\begin{equation}\label{eq:right-a-b-numeric}
 b^2\ge b(1)^2>0.68,\qquad a^2\le(\log r)^2+\frac{\pi^2}{3}.
\end{equation}
The first bound follows from $\Rea(s(1)^2)=-\frac12\Li_2(-1)=\pi^2/24$
and $|s(1)|^2\ge16/9$ (from \eqref{eq:unit-basic-bounds}), which give
$b(1)^2=(|s(1)|^2-\Rea(s(1)^2))/2\ge\frac89-\frac{\pi^2}{48}>0.68$.
For the second bound, \eqref{eq:dilog-inversion} gives
$\Rea(s^2)=(\log r)^2+\frac{\pi^2}{12}+\frac12\Li_2(-r^{-2})
\le(\log r)^2+\frac{\pi^2}{12}$; combining this with
$a^2=\Rea(s^2)+b^2$ and $b<\pi/2$ completes the proof of
\eqref{eq:right-a-b-numeric}.

For an integer $n\ge1$, we apply \eqref{eq:ct-integral2} with
$a/\sqrt n$ in place of $u$ and write
$\zeta=\frac{a}{\sqrt n}(1+i\eta)$ (so that $v=\frac{a}{\sqrt n}\eta$).
By \eqref{eq:ct-integral2}, \eqref{eq:unfolded-theta}, and
\eqref{eq:eta-transform}, the contribution of the term indexed by
$m\in\mathbb Z$ over any subinterval
$I\subset[-\pi\sqrt n/a,\pi\sqrt n/a]$ of $\eta$ is
\begin{equation}\label{eq:right-branch-integral}
 \frac{a}{4\pi^2\sqrt n}\int_I\int_{\mathbb R}
 e^{n\zeta+H_m(\chi)/\zeta}\mathfrak b_m(\zeta,\chi)\,d\chi\,d\eta,
\end{equation}
where
\[
 \mathfrak b_m(\zeta,\chi)
 =e^{-\zeta/12}\mathcal E(4\pi^2/\zeta)^{-2}
  e^{g_{(-1)^{m+1}r^{-1}}(\chi)/\zeta
     -\log\mathcal P_{(-1)^{m+1}r^{-1}}(\zeta,\chi)}.
\]
Replacing $H_0(\chi)$ by its quadratic approximation
$s^2-\frac{ir}{1+ir}\chi^2$ and $\mathfrak b_0(\zeta,\chi)$ by $1$ in
\eqref{eq:right-branch-integral} with $I=[0,3/2]$ (when
$3a/(2\sqrt n)\le\pi$) and integrating in $\chi$ leads to the integral
estimated in the next lemma.

\begin{lemma}\label{lem:exact-radial-truncation}
For $t\ge1$ and $0<\varrho<1$, let
\[
 \mathcal J
 =\int_0^{3/2}(1+i\eta)^{1/2}e^{-t(\eta-\varrho)^2/(1+i\eta)}\,d\eta.
\]
Then
\[
 \left|\sqrt{\frac{t}{\pi}}\,\frac{\mathcal J}{1+i\varrho}-1\right|
 \le\frac1{2t\sqrt{1+\varrho^2}}
   +\left(1+\frac1{2t\sqrt{1+\varrho^2}}\right)e^{-4t}
   +8\sqrt t\left(e^{-t\varrho^2}
                  +e^{-t(3-2\varrho)^2/13}\right).
\]
\end{lemma}
\begin{proof}
Write $z=1+i\eta$ and
\[
 \mathcal H(\eta)
 =z^{1/2}e^{t(1+i\varrho)^2/z}
 =\sum_{k\ge0}\frac{[t(1+i\varrho)^2]^k}{k!}z^{1/2-k},
\]
with the principal square root. Direct algebra gives
\begin{equation}\label{eq:radial-H-rel}
 e^{-t-2it\varrho}e^{it\eta}\mathcal H(\eta)
 =(1+i\eta)^{1/2}e^{-t(\eta-\varrho)^2/(1+i\eta)}.
\end{equation}
Define
\[
 \mathcal J_\infty
 =-\frac{e^{-t-2it\varrho}}{t^2}
   \int_{\mathbb R}e^{it\eta}\mathcal H''(\eta)\,d\eta,
\]
whose absolute convergence is verified below.
For $j\in\{1,2\}$, $k\ge0$, and $\eta\in\mathbb R$, we have
\[
 \left|\frac{d^j}{d\eta^j}
 \left(\frac{[t(1+i\varrho)^2]^k}{k!}
             (1+i\eta)^{1/2-k}\right)\right|
 \le \frac{(k+1)^j(2t)^k}{k!},
 \qquad
 \sum_{k\ge0}\frac{(k+1)^j(2t)^k}{k!}<\infty,
\]
since $|1+i\eta|\ge1$ and $1+\varrho^2<2$.
Thus both differentiated series converge uniformly on $\mathbb R$;
together with the pointwise convergence of the series for $\mathcal H$,
this justifies differentiating twice term by term and gives
\[
 \mathcal H''(\eta)
 =-\sum_{k\ge0}\frac{[t(1+i\varrho)^2]^k}{k!}
               \left(k-\frac12\right)\left(k+\frac12\right)
               (1+i\eta)^{-k-3/2}.
\]
Moreover, the sum of the moduli of these terms is at most
\[
 (1+\eta^2)^{-3/4}
 \sum_{k\ge0}\frac{(k+1)^2(2t)^k}{k!}.
\]
Since $(1+\eta^2)^{-3/4}$ is integrable on $\mathbb R$, this proves
the integrability of $\mathcal H''$ and justifies termwise
integration below. For every integer $k\ge0$ and
$\Rea w>0$, repeated integration by parts from
$\int_0^\infty y^{1/2}e^{-wy}\,dy=\frac12\sqrt\pi\,w^{-3/2}$ gives
\[
 \int_0^\infty y^{k+1/2}e^{-wy}\,dy
 =\frac{(2k+2)!\sqrt\pi}{4^{k+1}(k+1)!}w^{-k-3/2}.
\]
Since $y\mapsto y^{k+1/2}e^{-y}\mathbf1_{\{y>0\}}$ is continuous and
integrable on $\mathbb R$, and its Fourier transform
$\eta\mapsto\frac{(2k+2)!\sqrt\pi}{4^{k+1}(k+1)!}(1+i\eta)^{-k-3/2}$
is also integrable on $\mathbb R$, Fourier inversion yields
\[
 \int_{\mathbb R}e^{it\eta}(1+i\eta)^{-k-3/2}\,d\eta
 =\frac{2\sqrt\pi\,4^{k+1}(k+1)!}{(2k+2)!}e^{-t}t^{k+1/2}.
\]
Integrating the series for $\mathcal H''$ term by term therefore
gives
\begin{equation}\label{eq:radial-J-infty}
 \begin{aligned}
 \mathcal J_\infty
 &=\frac{\sqrt\pi\,e^{-2t(1+i\varrho)}}{t^{3/2}}
   \sum_{k\ge0}\frac{(2k-1)[2t(1+i\varrho)]^{2k}}{(2k)!}\\
 &=\sqrt{\frac\pi t}\left[1+i\varrho-\frac1{2t}
   -\left(1+i\varrho+\frac1{2t}\right)e^{-4t(1+i\varrho)}\right],
 \end{aligned}
\end{equation}
where the second equality uses
$\sum_{k\ge0}(2k-1)x^{2k}/(2k)!=x\sinh x-\cosh x$ with
$x=2t(1+i\varrho)$.

Define
\begin{align*}
 I_-&=\frac{\mathcal H(0)}{it}-\frac{\mathcal H'(0)}{(it)^2}
      +\frac1{(it)^2}\int_{-\infty}^0e^{it\eta}\mathcal H''(\eta)\,d\eta,\\
 I_+&=-\frac{e^{3it/2}\mathcal H(3/2)}{it}
      +\frac{e^{3it/2}\mathcal H'(3/2)}{(it)^2}
      +\frac1{(it)^2}\int_{3/2}^\infty e^{it\eta}\mathcal H''(\eta)\,d\eta.
\end{align*}
Integrating \eqref{eq:radial-H-rel} over $[0,3/2]$ and integrating
twice by parts gives
$\mathcal J_\infty-\mathcal J=e^{-t-2it\varrho}(I_-+I_+)$,
avoiding any improper integral of $e^{it\eta}\mathcal H(\eta)$.
To bound $I_-$ and $I_+$, note that
\[
 \Rea\!\left(\frac{(1+i\varrho)^2}{z}\right)
 =1-\frac{(\eta-\varrho)^2}{1+\eta^2}
 \le\begin{cases}
  1-\varrho^2,&\eta\le0,\\
  1-(3-2\varrho)^2/13,&\eta\ge3/2.
 \end{cases}
\]
Here the bound for $\eta\le0$ follows from
$(\eta-\varrho)^2/(1+\eta^2)-\varrho^2
=((1-\varrho^2)\eta^2-2\varrho\eta)/(1+\eta^2)\ge0$, while for
$\eta\ge3/2$ the derivative of $(\eta-\varrho)^2/(1+\eta^2)$ is
$2(\eta-\varrho)(1+\varrho\eta)/(1+\eta^2)^2>0$, so the minimum is
attained at $\eta=3/2$. Differentiating $\mathcal H$ and using
$1+\varrho^2<2$, we obtain
\begin{align*}
 |\mathcal H'(\eta)|
 &\le e^{t\Rea((1+i\varrho)^2/z)}
      \left(\tfrac12|z|^{-1/2}+2t|z|^{-3/2}\right),\\
 |\mathcal H''(\eta)|
 &\le e^{t\Rea((1+i\varrho)^2/z)}
      \left(\tfrac14|z|^{-3/2}+2t|z|^{-5/2}
                             +4t^2|z|^{-7/2}\right).
\end{align*}
At $\eta=0$ and $\eta=3/2$, we have $1\le|z|\le\sqrt{13}/2$ and hence
$|z|^{1/2}<1.344$. Using
$\int_0^\infty(1+\eta^2)^{-q/2}\,d\eta\le1+1/(q-1)$ ($q>1$) with
$q=3/2,5/2,7/2$ and $t\ge1$ gives
\[
 e^{-t}|I_\pm|
 \le e^{-\delta_0}\left(\frac{1.344}{t}
                        +\frac{1/2+2t}{t^2}
                        +\frac3{4t^2}
                        +\frac{10}{3t}
                        +\frac{28}{5}\right)
 <14e^{-\delta_0},
\]
where $\delta_0=t\varrho^2$ for $I_-$ and
$\delta_0=t(3-2\varrho)^2/13$ for $I_+$. Since
$|1+i\varrho|=\sqrt{1+\varrho^2}>1$ and $14/\sqrt\pi<8$, it follows
that
\[
 \left|\sqrt{\frac{t}{\pi}}\,\frac{\mathcal J_\infty-\mathcal J}{1+i\varrho}\right|
 \le\sqrt{\frac{t}{\pi}}\,e^{-t}(|I_-|+|I_+|)
 <8\sqrt t\left(e^{-t\varrho^2}
                +e^{-t(3-2\varrho)^2/13}\right).
\]
Combining this estimate with \eqref{eq:radial-J-infty} proves the
lemma.
\end{proof}

We now assume that $n\ge6400$ is an integer and
$0\le\log r\le\frac12\log n+1$ (that is, $1\le r\le e\sqrt n$).
Let $a_n=\sqrt{(\frac12\log n+1)^2+3.3}$.
Since $\Rea(s^2)=(\log r)^2+\frac{\pi^2}{12}+\frac12\Li_2(-r^{-2})
\ge\frac{\pi^2}{24}>0$, the bounds in \eqref{eq:right-a-b-numeric}
and $\pi^2/3<3.3$ give
\[
 1<a\le a_n,\qquad 0<\frac{b}{a}<1,\qquad
 \frac{a}{\sqrt n}\le\frac{a_n}{\sqrt n}<0.072,
\]
where $a_n/\sqrt n$ decreases for $n\ge6400$ because
$\frac{d}{dx}\bigl(((\log x+1)^2+3.3)/x^2\bigr)
=2\bigl(\log x+1-(\log x+1)^2-3.3\bigr)/x^3<0$ for $x>0$.
In particular, $8a/\sqrt n<0.576<\pi$, so $[-8,8]$ is a subinterval
of $[-\pi\sqrt n/a,\pi\sqrt n/a]$.

We first record uniform bounds for $\mathfrak b_m(\zeta,\chi)$ on
$|\eta|\le8$ (where $\zeta=\frac{a}{\sqrt n}(1+i\eta)$).
For $|\eta|\le8$, we have
$|\zeta|=\frac{a}{\sqrt n}\sqrt{1+\eta^2}<0.072\sqrt{65}<0.581<\pi$
and $\Rea(4\pi^2/\zeta)>4\pi^2/(65\cdot0.072)>8.43$, so
$|\log\mathcal E(4\pi^2/\zeta)^{-2}|<5\cdot10^{-4}$; for
$|\eta|\le3/2$, $\Rea(4\pi^2/\zeta)>4\pi^2/(3.25\cdot0.072)>168$, so
$|\log\mathcal E(4\pi^2/\zeta)^{-2}|<10^{-70}$.
For any $t\in[-1,1]$ and $\chi\in\mathbb R$, expanding
$\log\mathcal P_t(\zeta,\chi)-g_t(\chi)/\zeta$ as
\[
 \sum_{j\ge1}\frac{2(-1)^{j+1}(it)^j\cos(j\chi)}{j}
 \left(\frac1{2\sinh(j\zeta/2)}-\frac1{j\zeta}\right)
\]
and splitting the sum at $j=\lfloor\pi/|\zeta|\rfloor$ (using
$|1/(2\sinh(z/2))-1/z|<0.08|z|$ for $|z|\le\pi$,
$|2\sinh(j\zeta/2)|\ge ja/\sqrt n$, and
$\sum_{j>\lfloor\pi/|\zeta|\rfloor}j^{-2}\le|\zeta|/(\pi-|\zeta|)$)
yields
\[
 |\log\mathcal P_t(\zeta,\chi)-g_t(\chi)/\zeta|
 \le0.51+\frac{2\bigl(1+\sqrt{1+\eta^2}\bigr)}
              {\pi-(a/\sqrt n)\sqrt{1+\eta^2}}
 \qquad(|\eta|\le8).
\]
Since $|e^{-\zeta/12}|<1$ and $|(-1)^{m+1}r^{-1}|=r^{-1}\le1$, it
follows that, for every $m\in\mathbb Z$ and $\chi\in\mathbb R$,
\begin{equation}\label{eq:right-bp-bounds}
 |\mathfrak b_m(\zeta,\chi)|<2100\quad(|\eta|\le8),\qquad
 |\mathfrak b_m(\zeta,\chi)|<30\quad(|\eta|\le3),\qquad
 |\mathfrak b_m(\zeta,\chi)|<10.9\quad(|\eta|\le3/2).
\end{equation}

We split the estimates below into intervals of $\log r$: the largest
values of $a/\sqrt n$ occur for large $\log r$, whereas
$r^2/(1+r^2)$ is smallest at $r=1$.

\begin{lemma}\label{lem:effective-right-principal}
If $n\ge6400$ is an integer and $1\le r\le e\sqrt n$, then the
combined contribution of $m=0$ on $0\le\eta\le8$ and $m=-1$ on
$-8\le\eta\le0$ in \eqref{eq:right-branch-integral} differs from
$M_n(r)\cos\Psi_n(r)$ in modulus by less than $0.88M_n(r)$.
\end{lemma}
\begin{proof}
Set $\varrho=b/a$ and $\zeta=\frac{a}{\sqrt n}(1+i\eta)$.
Since $H_{-1}(\chi)=\overline{H_0(\chi)}$ and
$\mathfrak b_{-1}(\overline{\zeta},\chi)=\overline{\mathfrak b_0(\zeta,\chi)}$,
the contribution of $m=-1$ on $-8\le\eta\le0$ is the complex
conjugate of the contribution of $m=0$ on $0\le\eta\le8$.
It therefore suffices to show that the contribution of $m=0$ on
$0\le\eta\le8$ in \eqref{eq:right-branch-integral} differs from
$\frac12M_n(r)e^{i\Psi_n(r)}$ in modulus by less than
$0.88\cdot\frac12M_n(r)$.

Direct algebra gives
\[
 n\zeta+\frac{s^2}{\zeta}
 =2s\sqrt n-\frac{a\sqrt n(\eta-\varrho)^2}{1+i\eta}.
\]
By Lemma~\ref{lem:unfolding}, $H_0(0)=s^2$, so
\eqref{eq:all-branches-components} implies
$\Rea H_0(\chi)\le\Rea(s^2)-\frac{r^2\chi^2}{1+r^2}$ and
$|\Ima H_0(\chi)|\le\Ima(s^2)$. Hence, for $\eta\ge0$ and
$\chi\in\mathbb R$,
\begin{equation}\label{eq:right-common-gaussian}
 \Rea\!\left(n\zeta+\frac{H_0(\chi)}{\zeta}\right)
 \le2a\sqrt n-\frac{a\sqrt n(\eta-\varrho)^2}{1+\eta^2}
             -\frac{\sqrt n\,r^2\chi^2}{a(1+r^2)(1+\eta^2)}.
\end{equation}
The same bound holds with $H_0(\chi)$ replaced by
$s^2-\frac{ir\chi^2}{1+ir}$, since $\eta\ge0$ and
$\Ima\bigl(ir/(1+ir)\bigr)=r/(1+r^2)>0$.
Integrating in $\chi$ on $0\le\eta\le3/2$ and applying
\eqref{eq:amplitude-phase} yields
\begin{equation}\label{eq:right-model-normalization}
 \begin{aligned}
 &\frac{a}{4\pi^2\sqrt n}\int_0^{3/2}\int_{\mathbb R}
   \exp\!\left(n\zeta+\frac{s^2}{\zeta}
              -\frac{ir\chi^2}{(1+ir)\zeta}\right)\,d\chi\,d\eta\\
 &\qquad=\frac{a}{4\pi^2\sqrt n}
         \sqrt{\frac{\pi a(1+ir)}{\sqrt n\,ir}}\,e^{2s\sqrt n}\mathcal J
        =\frac{M_n(r)}{2}e^{i\Psi_n(r)}
         \sqrt{\frac{a\sqrt n}{\pi}}\,\frac{\mathcal J}{1+i\varrho},
 \end{aligned}
\end{equation}
where $\mathcal J$ is the integral in
Lemma~\ref{lem:exact-radial-truncation} with $t=a\sqrt n$.
Here the Gaussian integral in $\chi$ equals
$\sqrt{\pi\zeta(1+ir)/(ir)}$ with the principal square root. Indeed,
\[
 \arg\!\left(\frac{1+ir}{ir}\right)
 =\arctan r-\frac\pi2\in[-\pi/4,0),
 \qquad
 \arg(1+i\eta)\in[0,\arctan(3/2)].
\]
The sum of these arguments belongs to $(-\pi,\pi)$, so the
principal square roots factor as used above.

We next sharpen the bound on $|\log\mathfrak b_0(\zeta,\chi)|$ and
expand $H_0(\chi)$ on $0\le\eta\le3/2$, $|\chi|\le h_\chi$, where
$h_\chi\in\{3/4,1\}$. Put $d=(\pi/2-h_\chi)/2$ and
$k_\chi=(\sin d)^{-2}+(e^{2d/3}-1)^{-1}$.
For $x\ge0$ and either sign, define
$f_\pm(x)=\log(1-ir^{-1}e^{-(1+i\eta)x\pm i\chi})$.
When $0\le x\le2d/3$, the argument $-\pi/2-\eta x\pm\chi$ of
$-ir^{-1}e^{-(1+i\eta)x\pm i\chi}$ stays in $[-\pi+d,-2d]$, at
distance at least $d$ from $-\pi$, so
$|1-ir^{-1}e^{-(1+i\eta)x\pm i\chi}|\ge\sin d$.
For $x\ge2d/3$, this modulus is at least $1-e^{-x}$.
Therefore, for $r>1$,
\[
 \int_0^\infty|f_\pm''(x)|\,dx
 \le\frac{1+\eta^2}{r}
   \left(\frac1{\sin^2 d}\int_0^{2d/3}e^{-x}\,dx
        +\int_{2d/3}^\infty\frac{e^{-x}}{(1-e^{-x})^2}\,dx\right)
 \le\frac{1+\eta^2}{r}\,k_\chi.
\]
Applying the midpoint rule on each interval
$[ja/\sqrt n,(j+1)a/\sqrt n]$ (as in the proof of
\eqref{eq:left-midpoint}, with mesh $a/\sqrt n$ and termwise
integration of the logarithmic series giving
$\frac{\sqrt n}{a}\int_0^\infty(f_+(x)+f_-(x))\,dx=g_{-r^{-1}}(\chi)/\zeta$)
and summing the errors for $f_+$ and $f_-$ gives
$|\log\mathcal P_{-r^{-1}}(\zeta,\chi)-g_{-r^{-1}}(\chi)/\zeta|
\le\frac{a(1+\eta^2)k_\chi}{4r\sqrt n}$ for $r>1$.
For fixed $n$, $\eta$, and $\chi$, uniform convergence in
$r\in[1,2]$ of the series for $\log\mathcal P_{-r^{-1}}(\zeta,\chi)$
and \eqref{eq:prelim-gr} extends this bound to $r=1$ by continuity.
Since $1+\eta^2\le13/4$, $|\zeta|\le\sqrt{13/4}\,a/\sqrt n$, and
$|\log\mathcal E(4\pi^2/\zeta)^{-2}|<10^{-70}<10^{-20}$ on
$0\le\eta\le3/2$, we obtain
\begin{equation}\label{eq:right-log-b0}
 |\log\mathfrak b_0(\zeta,\chi)|
 \le\frac{13k_\chi a}{16r\sqrt n}
   +\frac{\sqrt{13/4}\,a}{12\sqrt n}+10^{-20}
 \qquad(0\le\eta\le3/2,\ |\chi|\le h_\chi).
\end{equation}

To bound $H_0(\chi)-s^2+\frac{ir\chi^2}{1+ir}$, set
$\psi=1-\cos\chi$. Conjugating \eqref{eq:unit-log-integral} (with
$r^{-1}$ in place of $r$) and subtracting from its value at $\chi=0$
gives
\[
 g_{-r^{-1}}(0)-g_{-r^{-1}}(\chi)
 =\int_0^{r^{-1}}-\log\!\left(1+\frac{2iy\psi}{(1-iy)^2}\right)\frac{dy}{y}.
\]
Since $|2iy/(1-iy)^2|=2y/(1+y^2)\le1$ for $0<y\le r^{-1}\le1$, we
have $\int_0^{r^{-1}}|2iy/(1-iy)^2|^j\,\frac{dy}{y}\le\frac2r$ for
every $j\ge1$. Expanding the logarithm in powers of
$2iy\psi/(1-iy)^2$ shows that the linear term of
$g_{-r^{-1}}(0)-g_{-r^{-1}}(\chi)$ in $\psi$ has coefficient
$-2i\int_0^{r^{-1}}(1-iy)^{-2}\,dy=\frac2{1+ir}$, and gives
\[
 \left|g_{-r^{-1}}(0)-g_{-r^{-1}}(\chi)-\frac{2\psi}{1+ir}\right|
 \le\frac2r\sum_{j\ge2}\frac{\psi^j}{j}
 \le\frac{\psi^2}{r(1-\psi)}.
\]
Since $H_0(\chi)-s^2+\frac{ir\chi^2}{1+ir}
=g_{-r^{-1}}(0)-g_{-r^{-1}}(\chi)-\frac{2\psi}{1+ir}
+\frac{2\psi-\chi^2}{1+ir}$, using $|1+ir|>r$,
$|2\psi-\chi^2|\le\chi^4/12$, and $\psi\le\chi^2/2$ yields
\begin{equation}\label{eq:right-H0-quartic}
 \left|H_0(\chi)-s^2+\frac{ir\chi^2}{1+ir}\right|
 \le\frac{\chi^4}{r}\left(\frac1{12}+\frac1{4(1-\psi)}\right)
 \le\frac{\kappa_{h_\chi}}{r}\chi^4
 \qquad(|\chi|\le h_\chi),
\end{equation}
where $\kappa_{3/4}=0.426$ and $\kappa_1=0.547$ (using $\psi<0.269$
for $|\chi|\le3/4$ and $\psi<0.460$ for $|\chi|\le1$).

Now suppose that $r$ lies in a range on which
$r^{-1}\le\epsilon_+$, $a/\sqrt n\le q$, $a_-\le a$, $b_-\le b$, and
$\varrho\le\varrho_+<3/2$. Put
\[
 z_0=\frac{h_\chi^2}{(13/4)q(1+\epsilon_+^2)},
 \qquad
 \delta=\frac{(3/2-\varrho_+)^2}{13/4},
 \qquad
 \delta_3=\frac{(3-\varrho_+)^2}{10}.
\]
We compare the $m=0$ integral in \eqref{eq:right-branch-integral} on
$0\le\eta\le3/2$ with \eqref{eq:right-model-normalization}.
For each $\eta\in[0,3/2]$, let
$B(\eta)=\frac{\sqrt n\,r^2}{a(1+r^2)(1+\eta^2)}$ and denote by
$G(\eta)=\sqrt{\pi\zeta(1+ir)/(ir)}\,
e^{2s\sqrt n-a\sqrt n(\eta-\varrho)^2/(1+i\eta)}$ the inner
$\chi$-integral in \eqref{eq:right-model-normalization}.
Notice that
\[
 \frac{\int_{\mathbb R}e^{-B(\eta)\chi^2}\,d\chi}
      {\sqrt{\pi|\zeta|\sqrt{1+r^2}/r}}
 =(1+\eta^2)^{1/4}(1+r^{-2})^{1/4}
 \le\left(\frac{13(1+\epsilon_+^2)}{4}\right)^{1/4}.
\]
To bound
$\bigl|\int_{\mathbb R}e^{n\zeta+H_0(\chi)/\zeta}\mathfrak b_0(\zeta,\chi)\,d\chi
-G(\eta)\bigr|/|G(\eta)|$, we proceed in three steps, applying
\eqref{eq:right-common-gaussian} at each step:
first truncate both $\chi$-integrals to $|\chi|\le h_\chi$ (using
\eqref{eq:right-bp-bounds} to bound $|\mathfrak b_0|<10.9$, the
tail inequality
$\int_{|\chi|>c}e^{-B\chi^2}\,d\chi\le e^{-Bc^2}/(Bc)$ for
$B,c>0$, and $B(\eta)h_\chi^2\ge z_0$);
next replace $\mathfrak b_0(\zeta,\chi)$ by $1$ on $|\chi|\le h_\chi$
(using \eqref{eq:right-log-b0});
finally replace $H_0(\chi)$ by $s^2-\frac{ir\chi^2}{1+ir}$ on
$|\chi|\le h_\chi$ (using
$|e^{w_1}-e^{w_2}|\le|w_1-w_2|\max\{e^{\Rea w_1},e^{\Rea w_2}\}$,
\eqref{eq:right-H0-quartic}, and
$\int_{\mathbb R}\chi^4e^{-B\chi^2}\,d\chi=3\sqrt\pi/(4B^{5/2})$).
The resulting three errors are bounded, respectively, by
\begin{equation}\label{eq:right-inner-errors}
 \begin{aligned}
 E_\chi&=11.9\left(\frac{13(1+\epsilon_+^2)}{4}\right)^{1/4}
          \frac{e^{-z_0}}{\sqrt{\pi z_0}},\\
 E_{\mathrm{mult}}
 &=\left(\frac{13(1+\epsilon_+^2)}{4}\right)^{1/4}\left\{
   \exp\!\left(\frac{13}{16}k_\chi\epsilon_+q
              +\frac{\sqrt{13/4}}{12}q+10^{-20}\right)-1\right\},\\
 E_4&=\frac34\kappa_{h_\chi}\epsilon_+q(13/4)^{7/4}(1+\epsilon_+^2)^{9/4}.
 \end{aligned}
\end{equation}
On $[0,3/2]$, bounding $(1+\eta^2)^{1/4}\le(13/4)^{1/4}$ and
$\frac{(\eta-\varrho)^2}{1+\eta^2}\ge\frac{4}{13}(\eta-\varrho)^2$
gives
\[
 \sqrt{\frac{a\sqrt n}{\pi(1+\varrho^2)}}
 \int_0^{3/2}(1+\eta^2)^{1/4}e^{-\frac{a\sqrt n(\eta-\varrho)^2}{1+\eta^2}}\,d\eta
 \le\frac{(13/4)^{3/4}}{\sqrt{1+\varrho^2}}<2.43.
\]
Hence, by \eqref{eq:right-model-normalization}, the modulus of the
difference between the $m=0$ integral in
\eqref{eq:right-branch-integral} on $0\le\eta\le3/2$ and
\eqref{eq:right-model-normalization}, divided by $\frac12M_n(r)$, is
at most $2.43(E_4+E_{\mathrm{mult}}+E_\chi)$.

For the remaining interval $3/2\le\eta\le8$, split into $[3/2,3]$
and $[3,8]$. Since $(\eta-\varrho)^2/(1+\eta^2)$ increases for
$\eta>\varrho$ and $\varrho\le\varrho_+<3/2$, it is at least
$\delta$ on $[3/2,3]$ and at least $\delta_3$ on $[3,8]$.
Integrating \eqref{eq:right-common-gaussian} in $\chi$ and dividing
by $\frac12M_n(r)$ (using \eqref{eq:amplitude-phase} to obtain the
prefactor
$\frac{(1+r^{-2})^{1/4}}{\sqrt{1+\varrho^2}}\sqrt{\frac{a\sqrt n}{\pi}}
\le2^{1/4}\sqrt{\frac{a\sqrt n}{\pi}}$), then bounding
$|\mathfrak b_0|$ by \eqref{eq:right-bp-bounds} and $\sqrt{1+\eta^2}$
by $\sqrt{10}$ on $[3/2,3]$ and $\sqrt{65}$ on $[3,8]$, shows that
the modulus of the contribution of $3/2\le\eta\le8$ to the $m=0$
integral in \eqref{eq:right-branch-integral}, divided by
$\frac12M_n(r)$, is at most
\[
 96\sqrt{a\sqrt n}\,e^{-\delta a\sqrt n}
 +57000\sqrt{a\sqrt n}\,e^{-\delta_3 a\sqrt n},
\]
since $\frac{2^{1/4}}{\sqrt\pi}\cdot\frac32\sqrt{10}\cdot30<96$ and
$\frac{2^{1/4}}{\sqrt\pi}\cdot5\sqrt{65}\cdot2100<57000$.
Whenever $80a_-\delta>1/2$ (which also gives $80a_-\delta_3>1/2$
since $\delta_3>\delta$), both terms decrease for
$a\sqrt n\ge80a_-$ (as $\sqrt t\,e^{-ct}$ decreases for
$t\ge1/(2c)$), giving the upper bound
\begin{equation}\label{eq:right-E-eta}
 E_\eta=96\sqrt{80a_-}\,e^{-80a_-\delta}
        +57000\sqrt{80a_-}\,e^{-80a_-\delta_3}.
\end{equation}

By \eqref{eq:right-model-normalization} and
Lemma~\ref{lem:exact-radial-truncation} (with $t=a\sqrt n$ and
$(3-2\varrho)^2/13\ge\delta$), the modulus of the difference between
\eqref{eq:right-model-normalization} and $\frac12M_n(r)e^{i\Psi_n(r)}$,
divided by $\frac12M_n(r)$, is at most
\begin{equation}\label{eq:right-EJ-raw}
 \frac1{2a\sqrt n}
 +\left(1+\frac1{2a\sqrt n}\right)e^{-4a\sqrt n}
 +8\sqrt{a\sqrt n}\left(e^{-b^2\sqrt n/a}+e^{-\delta a\sqrt n}\right).
\end{equation}

We now specialize to nine finite intervals $[\ell,h]$ of $\log r$
and to $\log r\ge5$.
On each finite interval $[\ell,h]$ in the table below, let
$a_-,a_+,b_-,b_+$ be the listed bounds satisfying
$a_-\le a\le a_+$ and $b_-\le b\le b_+$, and set $\epsilon_+=e^{-\ell}$,
$q=a_+/80$, $\varrho_+=b_+/a_-$, and $h_\chi=3/4$ (except on $[4,5]$,
where $h_\chi=1$). In every row, $\varrho_+<0.82$, $80a_-\delta>1/2$,
and $80b_-^2/a_+>1/2$. Thus, for each fixed $r\in[e^\ell,e^h]$, every
term in \eqref{eq:right-EJ-raw} decreases with $\sqrt n\ge80$, so
\eqref{eq:right-EJ-raw} is at most
\[
 E_J=\frac1{160a_-}+\left(1+\frac1{160a_-}\right)e^{-320a_-}
     +8\sqrt{80a_+}\left(e^{-80b_-^2/a_+}+e^{-80a_-\delta}\right).
\]
Consequently, the modulus of the difference between the $m=0$
integral in \eqref{eq:right-branch-integral} on $0\le\eta\le8$ and
$\frac12M_n(r)e^{i\Psi_n(r)}$, divided by $\frac12M_n(r)$, is at most
$\mathcal R=2.43(E_4+E_{\mathrm{mult}}+E_\chi)+E_J+E_\eta$.
All terminating decimals in the following table denote exact rational
numbers; outward-rounded evaluation of these explicit formulas gives
the stated strict upper bounds:
\[
\begin{array}{c|rrrr|r}
\text{range of }\log r&a_-&a_+&b_-&b_+&\mathcal R\text{ less than}\\\hline
{[0,0.0625]}&1.06975&1.10495&0.85623&0.87428&0.876\\
{[0.0625,0.125]}&1.10494&1.14106&0.87427&0.89219&0.788\\
{[0.125,0.25]}&1.14105&1.21600&0.89218&0.92752&0.742\\
{[0.25,0.5]}&1.21599&1.37623&0.92751&0.99542&0.672\\
{[0.5,1]}&1.37622&1.73266&0.99541&1.11586&0.584\\
{[1,2]}&1.73265&2.54406&1.11585&1.28797&0.519\\
{[2,3]}&2.54405&3.42799&1.28796&1.38921&0.421\\
{[3,4]}&3.42798&4.34983&1.38920&1.44868&0.629\\
{[4,5]}&4.34982&5.29407&1.44867&1.48482&0.199
\end{array}
\]

For $\log r\ge5$, set $h_\chi=1$,
$\epsilon_+=e^{-5}$, $q=0.072$, $a_-=5.29406$, $b_-=1.48481$, and
$\varrho_+=\pi/(2a_-)$ (since $b<\pi/2$). Here $80a_-\delta>1/2$, so
\eqref{eq:right-inner-errors} and \eqref{eq:right-E-eta} remain
valid, while \eqref{eq:right-EJ-raw} is at most
\[
 E_J'=\frac1{160a_-}+\left(1+\frac1{160a_-}\right)e^{-320a_-}
      +8\sqrt{a_n\sqrt n}
       \left(e^{-b_-^2\sqrt n/a_n}+e^{-\delta a_-\sqrt n}\right),
\]
using $a\le a_n$ in the first exponential term and
$\sqrt{a\sqrt n}\,e^{-\delta a\sqrt n}
\le\sqrt{a_-\sqrt n}\,e^{-\delta a_-\sqrt n}
\le\sqrt{a_n\sqrt n}\,e^{-\delta a_-\sqrt n}$ in the second.
With $x=\sqrt n\ge80$ and $A(x)=\sqrt{(\log x+1)^2+3.3}$, we have
$xA'(x)/A(x)=(\log x+1)/((\log x+1)^2+3.3)<0.168$, so
differentiation shows that both summands in the last term of $E_J'$
decrease for $n\ge6400$. Evaluating
$2.43(E_4+E_{\mathrm{mult}}+E_\chi)+E_J'+E_\eta$ at $n=6400$ yields
a bound less than $0.212$.

The endpoint bounds $a_-,a_+,b_-,b_+$ in the table and at
$\log r=5$ follow from \eqref{eq:dilog-inversion} and
\eqref{eq:parameters}:
\[
 \begin{aligned}
 \Rea(s^2)&=(\log r)^2+\frac{\pi^2}{12}
            +\frac12\sum_{j\ge1}\frac{(-1)^j r^{-2j}}{j^2},\\
 \Ima(s^2)&=\pi\log r
            +2\sum_{j\ge0}\frac{(-1)^j r^{-(2j+1)}}{(2j+1)^2},\\
 a^2&=\frac{|s^2|+\Rea(s^2)}{2},\qquad b=\frac{\Ima(s^2)}{2a}.
 \end{aligned}
\]
Each alternating series lies between its partial sums with $10000$
and $10001$ terms; outward-rounded interval arithmetic for the
exponential, trigonometric, and square-root expressions gives the
displayed endpoint and error bounds, and the strict monotonicity of
$a$ and $b$ proved above extends the endpoint bounds throughout each
interval (and to $\log r\ge5$).

In all cases, the modulus of the difference between the $m=0$
integral in \eqref{eq:right-branch-integral} on $0\le\eta\le8$ and
$\frac12M_n(r)e^{i\Psi_n(r)}$ is less than $0.44M_n(r)$.
Adding the complex conjugate contribution proves the lemma.
\end{proof}

\begin{lemma}\label{lem:effective-minor}
If $n\ge6400$ is an integer and $1\le r\le e\sqrt n$, then
\[
 \left|\frac1{2\pi}\int_{8a/\sqrt n\le|v|\le\pi}
  e^{n(a/\sqrt n+iv)}F\bigl(-r^2,e^{-(a/\sqrt n+iv)}\bigr)\,dv\right|
 <0.05M_n(r).
\]
\end{lemma}
\begin{proof}
Set $\alpha=19/20$ and $\zeta=\frac{a}{\sqrt n}+iv$. For each
integer $j\ge1$, define
\[
 S_j=\sum_{m\ge0}
     e^{-j|(m+1/2)a/\sqrt n-\log r|+ijmv}.
\]
Clearly $S_j(-v)=\overline{S_j(v)}$, so $|S_j|$ is even in $v$.
Write $\frac{\sqrt n\log r}{a}-\frac12=m_0+\tau$ with
$m_0\in\mathbb Z$ and $\tau\in[0,1)$. Extending the sum to all
$m\in\mathbb Z$ yields
\[
 \sum_{m\in\mathbb Z}
 e^{-j|(m+1/2)a/\sqrt n-\log r|+ijmv}
 =e^{ijm_0 v}\frac{\sinh\bigl(j(1-\tau)a/\sqrt n\bigr)
                   +e^{ijv}\sinh\bigl(j\tau a/\sqrt n\bigr)}
                  {\cosh(ja/\sqrt n)-\cos(jv)}.
\]
Since $\log r\ge0$, the sum over $m<0$ equals
$r^{-j}e^{-ja/(2\sqrt n)-ijv}/(1-e^{-ja/\sqrt n-ijv})$, whose
modulus is $r^{-j}/\sqrt{2(\cosh(ja/\sqrt n)-\cos(jv))}$ and whose
value at $v=0$ is $r^{-j}/(2\sinh(ja/(2\sqrt n)))$. Since
$\sinh(0)=0$ and $\sinh$ is convex on $[0,\infty)$, we have
$\sinh\bigl(j(1-\tau)a/\sqrt n\bigr)+\sinh\bigl(j\tau a/\sqrt n\bigr)
\le\sinh(ja/\sqrt n)$.
Subtracting the $m<0$ sum from the doubly infinite sum (using
$|S_j|\le S_j|_{v=0}$ for the first entry and the triangle inequality
for the second) gives
\begin{equation}\label{eq:minor-Sj-bound}
 \begin{aligned}
 |S_j|\le\min\Biggl\{
 &\coth\bigl(ja/(2\sqrt n)\bigr)
  -\frac{r^{-j}}{2\sinh\bigl(ja/(2\sqrt n)\bigr)},\\
 &\frac{\sinh(ja/\sqrt n)}{\cosh(ja/\sqrt n)-\cos(jv)}
  +\frac{r^{-j}}{\sqrt{2\bigl(\cosh(ja/\sqrt n)-\cos(jv)\bigr)}}
 \Biggr\}.
 \end{aligned}
\end{equation}

For any $w>0$ and $\phi\in\mathbb R$, put $\xi=\min\{w,w^{-1}\}$.
After extracting $\max\{1,w\}$ from $|1+we^{i\phi}|$, the ratio
$|1+\xi e^{\pm i\phi}|/|1+\alpha\xi e^{\pm i\phi}|$ (with sign $+$
for $w\le1$ and $-$ for $w>1$) is maximized at $\phi=0$, where
$\log(1+\xi)-\log(1+\alpha\xi)\le(1-\alpha)\xi$; hence
\[
 \log|1+we^{i\phi}|
 \le\log\max\{1,w\}+(1-\alpha)\xi
   +\log|1+\alpha\xi e^{\pm i\phi}|.
\]
Apply this inequality to each factor of $\mathcal P_r(\zeta,\theta)$
in \eqref{eq:P-product}, writing $m\ge0$ for the product index. In
the absolutely convergent Fourier series
$\log|1+\alpha\xi e^{\pm i\phi}|
=\sum_{j\ge1}\frac{(-1)^{j-1}\alpha^j\xi^j}{j}\cos(j\phi)$, the sign
of $\pm i\phi$ drops out because $\cos$ is even, and combining the
two factors with $\pm\theta$ yields
$2\cos(j\theta)\cos\bigl(j(\pi/2-(m+1/2)v)\bigr)$. Since
$\cos\bigl(j(\pi/2-(m+1/2)v)\bigr)
=\Rea\bigl(e^{ij(\pi/2-v/2)}e^{-ijmv}\bigr)$ and
$|S_j(-v)|=|S_j(v)|$, the double sum over $m\ge0$ and $j\ge1$
converges absolutely (as $\sum_{j\ge1}\frac{\alpha^j}{j}S_j|_{v=0}<\infty$),
and summing the $j$th Fourier terms over $m\ge0$ gives a quantity of
modulus at most $\frac{2\alpha^j}{j}|S_j|$. Moreover, if
$K\ge0$ is the number of integers $m\ge0$ with
$(m+1/2)a/\sqrt n<\log r$, then
\[
 2\sum_{m\ge0}\max\!\left\{\log r-\left(m+\frac12\right)\frac{a}{\sqrt n},\,0\right\}
 =2K\log r-K^2\frac{a}{\sqrt n}
 \le\frac{\sqrt n}{a}(\log r)^2.
\]
Summing over $m\ge0$ and multiplying by $a/\sqrt n$ therefore gives
\begin{equation}\label{eq:minor-logP}
 \frac{a}{\sqrt n}\log|\mathcal P_r(\zeta,\theta)|
 \le(\log r)^2
   +\frac{2(1-\alpha)a}{\sqrt n}S_1\Big|_{v=0}
   +\frac{2a}{\sqrt n}\sum_{j\ge1}\frac{\alpha^j}{j}|S_j|.
\end{equation}

Now fix a range of $\log r$ on which $\epsilon_-\le r^{-1}\le\epsilon_+$
and $a/\sqrt n\le y_+\le0.072$ (on a finite interval $[\ell,h]$, we
take $\epsilon_+=e^{-\ell}$, $\epsilon_-=e^{-h}$, and
$y_+=\min\{0.072,a_+/80\}$ for any upper bound $a_+\ge a(e^h)$).
For $0<y\le y_+$, $z>0$, and $j\ge1$, define
\[
 \begin{aligned}
 b_j(y,z)&=\frac{2/j+jy^2/6}{1+(z/\sinh(jy/2))^2}
          +\frac{y\epsilon_+^j}{2z},\\
 c_j(y)&=\frac2j+\frac{jy^2}{6}
        -\frac{\epsilon_-^j}{j}\left(1-\frac{j^2y^2}{24}\right).
 \end{aligned}
\]
In \eqref{eq:minor-Sj-bound}, replacing $r^{-j}$ by $\epsilon_-^j$ in
the first entry and by $\epsilon_+^j$ in the second, using
$\cosh(ja/\sqrt n)-\cos(jv)=2\sinh^2(ja/(2\sqrt n))+2\sin^2(jv/2)$
together with $\coth t\le t^{-1}+t/3$ and $t/\sinh t\ge1-t^2/6$
($t>0$), shows that whenever $a/\sqrt n\le y\le y_+$,
\begin{equation}\label{eq:minor-bc}
 \frac{a}{\sqrt n}|S_j|
 \le\min\bigl\{b_j(y,|\sin(jv/2)|),\,c_j(y)\bigr\},
\end{equation}
omitting the $b_j$ entry when $\sin(jv/2)=0$. Note that $c_j(y)$
increases with $y$, and $b_j(y,z)$ increases with $y$ and decreases
with $z$.

By evenness of $|S_j|$ in $v$, it suffices to bound
\eqref{eq:minor-logP} for $8a/\sqrt n\le v\le\pi$. We subdivide
$[\pi/320,\pi]$ into the intervals
\[
 I_k=\left[\frac{k\pi}{1600},\,\frac{(k+1)\pi}{1600}\right]
 \qquad(5\le k\le1599),
\]
and set $y_k=\min\{y_+,(k+1)\pi/12800\}$. For $1\le j\le20$, put
$z_{j,k}=0$ if $[jk\pi/3200,j(k+1)\pi/3200]$ contains an integer
multiple of $\pi$, and otherwise put
\[
 z_{j,k}=\min\!\left\{\left|\sin\!\left(\frac{jk\pi}{3200}\right)\right|,\,
                      \left|\sin\!\left(\frac{j(k+1)\pi}{3200}\right)\right|\right\};
\]
by concavity of $|\sin t|$ between consecutive multiples of $\pi$,
$z_{j,k}$ is the minimum of $|\sin(jv/2)|$ on $I_k$.
At every point $v\in I_k$ with $v\ge8a/\sqrt n$, we have
\[
 \frac a{\sqrt n}\le y_+,
 \qquad
 \frac a{\sqrt n}\le\frac v8
       \le\frac{(k+1)\pi}{12800}.
\]
Consequently,
\[
 \frac a{\sqrt n}
 \le\min\!\left\{y_+,\frac{(k+1)\pi}{12800}\right\}=y_k,
\]
so \eqref{eq:minor-bc} bounds the sum over $1\le j\le20$ in
\eqref{eq:minor-logP} by
\[
 q_k=\sum_{j=1}^{20}\frac{2\alpha^j}{j}
 \begin{cases}
  c_j(y_k),&z_{j,k}=0,\\
  \min\{b_j(y_k,z_{j,k}),\,c_j(y_k)\},&z_{j,k}>0.
 \end{cases}
\]
Whenever $8a/\sqrt n\le v\le\pi/320$, we have
$a/\sqrt n\le y_0:=\pi/2560$, and for $1\le j\le20$, $jv/2\le\pi/32$,
so $|\sin(jv/2)|\ge0.998\,jv/2\ge0.998\cdot4ja/\sqrt n$.
Substituting this bound into the second entry of
\eqref{eq:minor-Sj-bound} and using that $t/\sinh t$ decreases for
$t>0$ (so
$(ja/(2\sqrt n))/\sinh(ja/(2\sqrt n))\ge10y_0/\sinh(10y_0)$) bounds
the sum over $1\le j\le20$ in \eqref{eq:minor-logP} by
\[
 q_0=\sum_{j=1}^{20}\frac{2\alpha^j}{j}
     \left(\frac{2/j+jy_0^2/6}
                {1+\bigl(0.998\cdot80y_0/\sinh(10y_0)\bigr)^2}
           +\frac{\epsilon_+^j}{0.998\cdot8j}\right).
\]
For the remaining terms in \eqref{eq:minor-logP}, note first that
$2(1-\alpha)\frac{a}{\sqrt n}S_1|_{v=0}\le0.1c_1(y_+)
\le0.2+y_+^2/48-0.1\epsilon_-$ (since $\epsilon_-\le1$). For
$j>20$, using $|S_j|\le\coth(ja/(2\sqrt n))$ from
\eqref{eq:minor-Sj-bound} and $\coth t\le t^{-1}+1$ ($t>0$) gives
$\frac{a}{\sqrt n}|S_j|\le\frac2j+y_+$, so the terms with $j>20$
contribute at most $4\mathcal T_2+2y_+\mathcal T_1$, where
$\mathcal T_\nu=\sum_{j>20}\alpha^j/j^\nu$ ($\nu=1,2$). Therefore,
for all $8a/\sqrt n\le|v|\le\pi$ and $\theta\in\mathbb R$,
\begin{equation}\label{eq:minor-qstar}
 \frac{a}{\sqrt n}\log|\mathcal P_r(\zeta,\theta)|
 \le(\log r)^2+q_*,
 \qquad
 q_*=0.2+\frac{y_+^2}{48}-0.1\epsilon_-+4\mathcal T_2+2y_+\mathcal T_1
     +\max\{q_0,q_5,\ldots,q_{1599}\}.
\end{equation}

The quantity $a(e^t)^2-t^2$ increases with $t\ge0$, because
\eqref{eq:dilog-inversion}, $a^2=\Rea(s^2)+b^2$, and $b'(r)>0$ give
\[
 \frac{d}{dt}\bigl(a(e^t)^2-t^2\bigr)
 =\log(1+e^{-2t})+2b(e^t)\frac{d}{dt}b(e^t)>0.
\]
The alternating series in the proof of
Lemma~\ref{lem:effective-right-principal} give the endpoint bounds
for $a(e^\ell)^2-\ell^2$ and $a(e^h)$ in the table below.
For $\mathcal T_\nu$, we sum through $j=400$ and bound the tail by
$\alpha^{401}/(401^\nu(1-\alpha))$. Outward-rounded interval
arithmetic gives lower and upper bounds for the sines and hyperbolic
sines in $q_0$ and $q_k$, while integer arithmetic on $jk$ and
$j(k+1)$ determines whether $[jk\pi/3200,j(k+1)\pi/3200]$ contains
a multiple of $\pi$. Thus the inequalities above verify each entry
of the table:
\[
\begin{array}{c|r|r|r}
\text{range of }\log r&a(e^\ell)^2-\ell^2\text{ greater than}
 &a(e^h)\text{ less than}&q_*\text{ less than}\\\hline
{[0,0.04]}&1.14437&1.09218&0.887\\
{[0.04,0.08]}&1.19123&1.11497&0.951\\
{[0.08,0.12]}&1.23674&1.13814&1.007\\
{[0.12,0.16]}&1.28095&1.16169&1.055\\
{[0.16,0.20]}&1.32390&1.18560&1.097\\
{[0.20,0.25]}&1.36563&1.21600&1.143\\
{[0.25,0.32]}&1.41614&1.25951&1.198\\
{[0.32,0.40]}&1.48394&1.31054&1.249\\
{[0.40,0.50]}&1.55750&1.37623&1.301\\
{[0.50,0.65]}&1.64399&1.47850&1.361\\
{[0.65,0.85]}&1.76345&1.62127&1.416\\
{[0.85,1]}&1.90601&1.73266&1.444\\
{[1,1.5]}&2.00209&2.12575&1.503\\
{[1.5,2.25]}&2.26877&2.76012&1.537\\
{[2.25,4]}&2.55570&4.34983&1.571\\
{[4,5]}&2.92097&5.29407&1.584\\
{[5,6]}&3.02712&6.25266&1.593
\end{array}
\]
In every row of the table, $a^2-(\log r)^2-q_*>0.17a$.
For $\log r\ge6$, take $\epsilon_+=e^{-6}$, $\epsilon_-=0$, and
$y_+=0.072$. The same evaluation gives $q_*<1.593$, whereas the
alternating series at $\log r=6$ gives
$a^2-(\log r)^2\ge a(e^6)^2-36>3.09568$, whence
$a^2-(\log r)^2-q_*>1.50$.

By \eqref{eq:ct-integral} and \eqref{eq:minor-qstar}, on
$8a/\sqrt n\le|v|\le\pi$ we have
\[
 \left|e^{n\zeta}F(-r^2,e^{-\zeta})\right|
 \le e^{2a\sqrt n}
     \exp\!\left(-\frac{\sqrt n}{a}\bigl(a^2-(\log r)^2-q_*\bigr)\right).
\]
Since $|s|(1+r^2)^{1/4}/\sqrt r=|s|(1+r^{-2})^{1/4}>1$,
\eqref{eq:amplitude-phase} gives $e^{2a\sqrt n}/M_n(r)<2\pi n$.
Since the set $\{v\in\mathbb R:8a/\sqrt n\le|v|\le\pi\}$ has
measure less than $2\pi$, integrating and dividing by $2\pi M_n(r)$
yields
\[
 \frac1{M_n(r)}
 \left|\frac1{2\pi}\int_{8a/\sqrt n\le|v|\le\pi}
  e^{n\zeta}F(-r^2,e^{-\zeta})\,dv\right|
 <\begin{cases}
  2\pi n\,e^{-0.17\sqrt n},&0\le\log r\le6,\\
  2\pi n\,e^{-1.50\sqrt n/a_n},&\log r\ge6,
 \end{cases}
\]
using $a\le a_n$ in the second range. Writing $x=\sqrt n\ge80$ and
$A(x)=\sqrt{(\log x+1)^2+3.3}$, both bounds decrease for $n\ge6400$
(for the second, since $x/A(x)>1/0.072>13.8$ and
$xA'(x)/A(x)<0.168$ as in the proof of
Lemma~\ref{lem:effective-right-principal}). At $n=6400$, they are
less than $0.049884$ and $0.000027$, respectively, proving the lemma.
\end{proof}

It remains to estimate the integrals in \eqref{eq:right-branch-integral}
with $m\in\mathbb Z\setminus\{0,-1\}$ on $[-8,8]$, together with $m=0$
on $[-8,0]$ and $m=-1$ on $[0,8]$.

\begin{lemma}\label{lem:effective-branch-tails}
If $n\ge6400$ is an integer and $1\le r\le e\sqrt n$, then, with
$\zeta=\frac{a}{\sqrt n}(1+i\eta)$,
\[
 \begin{aligned}
 \frac{a}{4\pi^2\sqrt n}\Bigg(
  &\sum_{m\in\mathbb Z\setminus\{0,-1\}}
   \int_{-8}^{8}\int_{\mathbb R}
   \left|e^{n\zeta+H_m(\chi)/\zeta}\mathfrak b_m(\zeta,\chi)\right|\,d\chi\,d\eta\\
  &+\int_{-8}^{0}\int_{\mathbb R}
   \left|e^{n\zeta+H_0(\chi)/\zeta}\mathfrak b_0(\zeta,\chi)\right|\,d\chi\,d\eta
  +\int_{0}^{8}\int_{\mathbb R}
   \left|e^{n\zeta+H_{-1}(\chi)/\zeta}\mathfrak b_{-1}(\zeta,\chi)\right|\,d\chi\,d\eta
 \Bigg)<0.002M_n(r).
 \end{aligned}
\]
\end{lemma}
\begin{proof}
By \eqref{eq:right-bp-bounds}, $|\mathfrak b_m(\zeta,\chi)|<2100$
for all $m\in\mathbb Z$, $|\eta|\le8$, and $\chi\in\mathbb R$.
For $\log r\ge1$, we have $r^{-1}\le e^{-1}<1$. Since
$\mathcal P_{-r^{-1}}(\zeta,\chi)=\mathcal P_{r^{-1}}(\zeta,\chi+\pi)$
and $g_{-r^{-1}}(\chi)=g_{r^{-1}}(\chi+\pi)$, applying
\eqref{eq:left-midpoint} with $r^{-1}$ in place of $r$ and
$a/\sqrt n$ in place of $u$ gives, for every $\chi\in\mathbb R$ and
$|\eta|\le8$,
\[
 \left|\log\mathcal P_{\pm r^{-1}}(\zeta,\chi)
       -\frac{g_{\pm r^{-1}}(\chi)}{\zeta}\right|
 \le\frac{(a/\sqrt n)(1+\eta^2)r^{-1}}{4(1-r^{-1})}
 \le\frac{0.072\cdot65e^{-1}}{4(1-e^{-1})}.
\]
Using $|e^{-\zeta/12}|<1$ and the bound
$|\log\mathcal E(4\pi^2/\zeta)^{-2}|<5\cdot10^{-4}$ preceding
\eqref{eq:right-bp-bounds}, we obtain
\[
 |\mathfrak b_m(\zeta,\chi)|
 \le\exp\!\left(\frac{0.072\cdot65e^{-1}}{4(1-e^{-1})}+0.0005\right)<2
 \qquad(\log r\ge1,\ |\eta|\le8).
\]
By \eqref{eq:amplitude-phase} and $|s|(1+r^{-2})^{1/4}>a$, we have
$M_n(r)>\frac{a}{2\pi n}e^{2a\sqrt n}$. Moreover,
\[
 \int_{\mathbb R}e^{-\frac{\sqrt n\chi^2}{130a}}\,d\chi
 =\int_{|\chi|>\pi/2}e^{-\frac{\sqrt n(|\chi|-\pi/2)^2}{130a}}\,d\chi
 =\sqrt{\frac{130\pi a}{\sqrt n}}.
\]
Hence, for any $m\in\mathbb Z$, any subinterval $I\subset[-8,8]$,
and any measurable set $S\subset\mathbb R$ on which
$\Rea(n\zeta+H_m(\chi)/\zeta)\le2a\sqrt n-\delta-\frac{\sqrt n\chi^2}{130a}$
(or, when $S=\{\chi\in\mathbb R:|\chi|>\pi/2\}$,
$\Rea(n\zeta+H_m(\chi)/\zeta)\le2a\sqrt n-\delta
-\frac{\sqrt n(|\chi|-\pi/2)^2}{130a}$), we have
\begin{equation}\label{eq:tails-general}
 \frac1{M_n(r)}\frac{a}{4\pi^2\sqrt n}\int_I\int_S
 \left|e^{n\zeta+H_m(\chi)/\zeta}\mathfrak b_m(\zeta,\chi)\right|\,d\chi\,d\eta
 \le\frac{|I|}{16}\sqrt{a\sqrt n}\,e^{-\delta}
 \begin{cases}
  109000,&0\le\log r<1,\\
  110,&\log r\ge1,
 \end{cases}
\end{equation}
since $2100\cdot16\sqrt{130\pi}/(2\pi)<109000$ and
$2\cdot16\sqrt{130\pi}/(2\pi)<110$.

We first estimate the terms with $m\in\mathbb Z\setminus\{0,-1\}$
and $|\pi(m+1/2)|\le40a$. For $\beta\ge\pi/2$, let
\[
 w_\beta=(\log r+i\beta)^2+\frac{\pi^2}{3}+2\Li_2(i/r),
 \qquad
 x_\beta=\bigl(\Rea\sqrt{w_\beta}\bigr)^2-(\log r)^2.
\]
Since $g_{r^{-1}}(0)=-2\Li_2(-i/r)$, we have
$\Rea w_\beta=(\log r)^2-\beta^2+\frac{\pi^2}{3}-\Rea g_{r^{-1}}(0)$
and $\Ima w_\beta=2\beta\log r+\Ima g_{r^{-1}}(0)>0$, so
$\sqrt{w_\beta}$ lies in the open first quadrant. Using
$(\Rea\sqrt{w_\beta})^2\bigl((\Rea\sqrt{w_\beta})^2-\Rea w_\beta\bigr)
=\frac14(\Ima w_\beta)^2$ and substituting
$(\Rea\sqrt{w_\beta})^2=x_\beta+(\log r)^2$ shows that $x_\beta$ is
a root of the monic quadratic
\[
 \begin{aligned}
 f_\beta(x)&=x^2+\left((\log r)^2+\beta^2-\frac{\pi^2}{3}+\Rea g_{r^{-1}}(0)\right)x\\
 &\quad-(\log r)^2\left(\frac{\pi^2}{3}-\Rea g_{r^{-1}}(0)\right)
       -\beta(\log r)\Ima g_{r^{-1}}(0)
       -\frac14\bigl(\Ima g_{r^{-1}}(0)\bigr)^2.
 \end{aligned}
\]
Since $\Rea g_{r^{-1}}(0)=-\frac12\Li_2(-r^{-2})\in(0,\pi^2/24]$, we
have $\frac{\pi^2}{3}-\Rea g_{r^{-1}}(0)\ge\frac{7\pi^2}{24}>0$.
Writing
$2x_\beta=|w_\beta|-\bigl((\log r)^2+\beta^2-\frac{\pi^2}{3}+\Rea g_{r^{-1}}(0)\bigr)$,
the difference of the squares of the two terms on the right-hand side
equals
$4(\log r)^2\bigl(\frac{\pi^2}{3}-\Rea g_{r^{-1}}(0)\bigr)
+4\beta(\log r)\Ima g_{r^{-1}}(0)+\bigl(\Ima g_{r^{-1}}(0)\bigr)^2>0$,
which implies
$|w_\beta|>\bigl|(\log r)^2+\beta^2-\frac{\pi^2}{3}+\Rea g_{r^{-1}}(0)\bigr|$
and hence $x_\beta>0$. Thus $x_\beta$ is the unique positive root of
$f_\beta$ (as $f_\beta(0)<0$).
Using $\Ima g_{r^{-1}}(0)=2\int_0^{r^{-1}}\frac{\arctan t}{t}\,dt\le2r^{-1}$
(whence $\Ima g_{r^{-1}}(0)\le2$ and
$(\log r)\Ima g_{r^{-1}}(0)\le2r^{-1}\log r\le2/e$), we have
\[
 \begin{aligned}
 f_\beta\!\left(\frac{\pi^2}{3}-\Rea g_{r^{-1}}(0)\right)
 &=\beta^2\left(\frac{\pi^2}{3}-\Rea g_{r^{-1}}(0)\right)
   -\beta(\log r)\Ima g_{r^{-1}}(0)
   -\frac14\bigl(\Ima g_{r^{-1}}(0)\bigr)^2\\
 &\ge\frac{7\pi^2}{24}\beta^2-\frac{2\beta}{e}-1>0
 \qquad(\beta\ge\pi/2).
 \end{aligned}
\]
Consequently, $0<x_\beta<\frac{\pi^2}{3}-\Rea g_{r^{-1}}(0)<3.3$.
By \eqref{eq:dilog-inversion}, $w_{\pi/2}=s^2$, so
$x_{\pi/2}=a^2-(\log r)^2$. Subtracting $f_\beta(x_\beta)=0$ from
$f_{\pi/2}(x_{\pi/2})=0$ gives
\[
 a^2-(\log r)^2-x_\beta
 =(\beta^2-\pi^2/4)
  \frac{a^2-(\log r)^2-(\log r)\Ima g_{r^{-1}}(0)/(\beta+\pi/2)}
       {a^2+x_\beta+\beta^2-\pi^2/3+\Rea g_{r^{-1}}(0)}.
\]
By $x_\beta<\frac{\pi^2}{3}-\Rea g_{r^{-1}}(0)$ and
\eqref{eq:right-a-b-numeric}, the denominator is less than
$(\log r)^2+\beta^2+\pi^2/3$; since both
$\frac{\beta^2-\pi^2/4}{(\log r)^2+\beta^2+\pi^2/3}$ and the
positive factor
$a^2-(\log r)^2-(\log r)\Ima g_{r^{-1}}(0)/(\beta+\pi/2)$
(positive because $a^2-(\log r)^2>b^2>0.68$ by
\eqref{eq:dilog-inversion} and \eqref{eq:right-a-b-numeric}, whereas
$(\log r)\Ima g_{r^{-1}}(0)/(\beta+\pi/2)\le(2/e)/(2\pi)<0.12$)
increase with $\beta\ge3\pi/2$, the resulting lower bound increases
with $\beta$.

For $m\in\mathbb Z\setminus\{0,-1\}$, set $\beta=|\pi(m+1/2)|\ge3\pi/2$.
By \eqref{eq:all-branches-components}, $\Rea H_m(\chi)\le\Rea w_\beta-\frac{\chi^2}{1+r^{-2}}$
and $|\Ima H_m(\chi)|\le\Ima w_\beta$. Since $(1+r^{-2})(1+\eta^2)\le130$
for $r\ge1$ and $|\eta|\le8$, applying \eqref{eq:unit-sup-id} to
$z=w_\beta$ (with $|\eta|\ge0$ in place of $\eta$) gives
\[
 \begin{aligned}
 \Rea\!\left(n\zeta+\frac{H_m(\chi)}{\zeta}\right)
 &\le a\sqrt n+\frac{\sqrt n}{a}\,
      \frac{\Rea w_\beta+|\eta|\Ima w_\beta}{1+\eta^2}
      -\frac{\sqrt n\,\chi^2}{130a}\\
 &\le2a\sqrt n-\frac{\sqrt n}{a}\bigl(a^2-(\log r)^2-x_\beta\bigr)
             -\frac{\sqrt n\,\chi^2}{130a}.
 \end{aligned}
\]
On any interval $[\ell,h]$ of $\log r$, using
$\Ima g_{r^{-1}}(0)\le2r^{-1}$, monotonicity of $a$ and of
$a^2-(\log r)^2$, and $\beta\ge\beta_*\in\{3\pi/2,5\pi/2\}$ yields
\[
 \frac{a^2-(\log r)^2-x_\beta}{a}
 \ge\frac{(\beta_*^2-\pi^2/4)
          \bigl(a(e^\ell)^2-\ell^2
                -\max_{\ell\le t\le h}2te^{-t}/(\beta_*+\pi/2)\bigr)}
         {a(e^h)(h^2+\beta_*^2+\pi^2/3)},
\]
where the maximum of $2te^{-t}$ on $[\ell,h]$ is attained at the
point of $[\ell,h]$ closest to $1$. Evaluating this lower bound with
the endpoint bounds from the tables in
Lemmas~\ref{lem:effective-right-principal}
and~\ref{lem:effective-minor} gives the following strict lower bounds
for $(a^2-(\log r)^2-x_\beta)/a$:
\[
 \begin{array}{c|ccccc}
  \text{range of }\log r&[0,1]&[1,2]&[2,3]&[3,4]&[4,5]\\\hline
  \beta\ge3\pi/2&0.44&0.49&0.39&0.29&0.21\\
  \beta\ge5\pi/2&0.55&0.64&0.56&0.45&0.36
 \end{array}
\]
There are exactly two integers $m$ with $|\pi(m+1/2)|=3\pi/2$, and
at most $80a/\pi-1<26a+2$ integers $m\in\mathbb Z\setminus\{0,-1\}$
with $|\pi(m+1/2)|\le40a$. Applying \eqref{eq:tails-general} with
$I=[-8,8]$ and $S=\mathbb R$ bounds the sum over
$m\in\mathbb Z\setminus\{0,-1\}$ with $|\pi(m+1/2)|\le40a$, divided
by $M_n(r)$, by
\begin{align*}
 109000\sqrt{1.733\sqrt n}\,48e^{-0.44\sqrt n}
 &<10^{-6}
 &&\text{for }0\le\log r<1,\\
 110\sqrt{5.295\sqrt n}\bigl(2e^{-0.21\sqrt n}+140e^{-0.36\sqrt n}\bigr)
 &<0.00023
 &&\text{for }1\le\log r<5,
\end{align*}
using $a(e)<1.733$, $a(e^5)<5.295$, and monotonicity in $n\ge6400$.
For $\log r\ge5$, using $a^2-(\log r)^2\ge a(e^5)^2-25>3.02712$
(from the table in Lemma~\ref{lem:effective-minor}),
$(\log r)\Ima g_{r^{-1}}(0)\le10e^{-5}$, $a\le a_n$, and
$\log r\le\frac12\log n+1$ shows that
$\frac{\sqrt n}{a}(a^2-(\log r)^2-x_\beta)$ is at least
\[
 \delta_1=\frac{59.5\sqrt n}{a_n((\frac12\log n+1)^2+25.5)}
 \quad\text{for }|\pi(m+1/2)|=3\pi/2,
\]
and
\[
 \delta_2=\frac{178\sqrt n}{a_n((\frac12\log n+1)^2+65)}
 \quad\text{for }|\pi(m+1/2)|\ge5\pi/2.
\]
Hence, by \eqref{eq:tails-general}, the sum over
$m\in\mathbb Z\setminus\{0,-1\}$ with $|\pi(m+1/2)|\le40a$, divided
by $M_n(r)$, is at most
\[
 110\sqrt{a_n\sqrt n}\bigl(2e^{-\delta_1}+(26a_n+2)e^{-\delta_2}\bigr).
\]
Writing $x=\sqrt n\ge80$, $t=\log x+1$, and $A(x)=\sqrt{t^2+3.3}$,
we have $xA'(x)/A(x)=t/(t^2+3.3)<0.168$. Thus
$x\frac{d}{dx}\log\bigl(\sqrt{xA(x)}(26A(x)+2)\bigr)<1$, whereas
$x\frac{d}{dx}\log\delta_1=1-\frac{t}{t^2+3.3}-\frac{2t}{t^2+25.5}
>1-0.168-1/\sqrt{25.5}>0.63$ and
$x\frac{d}{dx}\log\delta_2=1-\frac{t}{t^2+3.3}-\frac{2t}{t^2+65}
>1-0.168-1/\sqrt{65}>0.70$.
Hence $\delta_1$ and $\delta_2$ increase with $x\ge80$, with
$0.63\delta_1>9>1$ and $0.70\delta_2>18>1$ at $x=80$, so both
summands decrease for $n\ge6400$. At $n=6400$, the bound is less
than $0.00098$.

Next, consider $m\in\mathbb Z$ with $|\pi(m+1/2)|>40a$. Direct
expansion gives
\[
 \Rea\!\left(\frac{(\log r+i\pi(m+1/2))^2}{\zeta}\right)
 =\frac{\sqrt n}{a}
  \left((\log r)^2-\frac{(\pi(m+1/2)-\eta\log r)^2}{1+\eta^2}\right).
\]
Since $\log r<a$ and $|\eta|\le8$, we have
$|\pi(m+1/2)-\eta\log r|>0.8|\pi(m+1/2)|$. Also
$|\pi^2/3-g_{(-1)^{m+1}r^{-1}}(\chi)|\le2\pi^2/3<6.6$. Because
$1+\eta^2\le65$ and
$(0.64/65-1/200)\pi^2(m+1/2)^2>(0.64/65-1/200)\cdot1600>6.6$, it
follows that
\[
 \Rea\!\left(n\zeta+\frac{H_m(\chi)}{\zeta}\right)
 \le2a\sqrt n-\frac{\sqrt n\,\pi^2(m+1/2)^2}{200a}
             -\frac{\sqrt n\,\chi^2}{130a}.
\]
Splitting the middle term into two equal halves, the first half is
at least $4a\sqrt n$, while summing the exponential of the second
half over all $m\in\mathbb Z$ yields
\[
 \sum_{m\in\mathbb Z}e^{-\frac{\sqrt n\,\pi^2(m+1/2)^2}{400a}}
 \le2+2\int_0^\infty e^{-\frac{\sqrt n\,\pi^2t^2}{400a}}\,dt
 =2+\sqrt{\frac{400a}{\pi\sqrt n}}<6,
\]
using $a/\sqrt n<0.072$. By \eqref{eq:tails-general} with $I=[-8,8]$
and $S=\mathbb R$, the sum over $|\pi(m+1/2)|>40a$, divided by
$M_n(r)$, is at most
$6\cdot109000\sqrt{a\sqrt n}\,e^{-4a\sqrt n}
\le6\cdot109000\sqrt{80}\,e^{-320}<10^{-6}$.

Finally, we estimate the remaining two integrals, $m=-1$ on $[0,8]$
and $m=0$ on $[-8,0]$. By $H_0(\chi)=\overline{H_{-1}(\chi)}$ and
$\mathfrak b_0(\overline\zeta,\chi)=\overline{\mathfrak b_{-1}(\zeta,\chi)}$,
substituting $\eta\mapsto-\eta$ shows that the two integrals are
equal, so we double the bound for $m=-1$ on $I=[0,8]$ (giving
$2\cdot(|I|/16)=1$ in \eqref{eq:tails-general}).
Since $H_{-1}(0)=\overline{s}^{\,2}$, \eqref{eq:all-branches-components}
gives $\Rea H_{-1}(\chi)\le a^2-b^2-\frac{\chi^2}{1+r^{-2}}$, while
$\Ima H_{-1}(\chi)=-\pi\log r-\Ima g_{r^{-1}}(\chi)$.
Note that $\Ima H_{-1}(\chi)\le0$ whenever $\log r\ge1$ (since
$|\Ima g_{r^{-1}}(\chi)|\le2r^{-1}\le2e^{-1}<\pi$) or whenever
$|\chi|\le\pi/2$ (since $\Ima g_{r^{-1}}(\chi)\ge0$ by
\eqref{eq:unit-log-integral}). Whenever $\Ima H_{-1}(\chi)\le0$ and
$\eta\ge0$, since $a^2-b^2=\Rea(s^2)>0$, we obtain
\[
 \Rea\!\left(n\zeta+\frac{H_{-1}(\chi)}{\zeta}\right)
 \le2a\sqrt n-\frac{\sqrt n\,b^2}{a}-\frac{\sqrt n\,\chi^2}{130a}.
\]
For $\log r\ge1$, the table in
Lemma~\ref{lem:effective-right-principal} and monotonicity of $b$
give $b^2\ge b(e)^2\ge1.11585^2>1.245$. By \eqref{eq:tails-general}
and $a\le a_n$, the sum of the $m=-1$ ($[0,8]$) and $m=0$ ($[-8,0]$)
integrals, divided by $M_n(r)$, is at most
$110\sqrt{a_n\sqrt n}\,e^{-1.245\sqrt n/a_n}$. Writing $x=\sqrt n\ge80$
and $A(x)=\sqrt{(\log x+1)^2+3.3}$, we have $x/A(x)>13.8$ and
$xA'(x)/A(x)<0.168$, so this bound decreases for $n\ge6400$ and is
less than $0.00006$ at $n=6400$.
For $0\le\log r<1$ and $|\chi|\le\pi/2$, the table in
Lemma~\ref{lem:effective-right-principal} gives
$b^2\ge b(1)^2\ge0.85623^2>0.733$ and $a<1.733$; by
\eqref{eq:tails-general}, the corresponding bound is
$109000\sqrt{1.733\sqrt n}\,e^{-0.733\sqrt n/1.733}<10^{-6}$ for
$n\ge6400$.

For $0\le\log r<1$ and $|\chi|>\pi/2$, write $\varrho=b/a$. By
\eqref{eq:all-branches-components} and \eqref{eq:dilog-inversion},
$|\Ima H_{-1}(\chi)|\le\pi\log r+\Ima g_{r^{-1}}(0)=\Ima(s^2)$.
Using \eqref{eq:unit-s2-id}, $\chi^2\ge\pi^2/4+(|\chi|-\pi/2)^2$,
and $\eta^2\le2\varrho^2+2(\eta-\varrho)^2$, we obtain for
$0\le\eta\le8$ that
\[
 2a\sqrt n-\Rea\!\left(n\zeta+\frac{H_{-1}(\chi)}{\zeta}\right)
 \ge\frac{\sqrt n}{a}\,
    \frac{a^2(\eta-\varrho)^2+\frac{r^2}{1+r^2}\frac{\pi^2}{4}}
         {1+2\varrho^2+2(\eta-\varrho)^2}
   +\frac{\sqrt n\,(|\chi|-\pi/2)^2}{130a}.
\]
On a subinterval of $[0,1]$ where $a_-\le a\le a_+$,
$\varrho\le\varrho_+$, and $\frac{r^2}{1+r^2}\ge c_-$, the first
term on the right-hand side is at least $\delta_-\sqrt n$, where
\[
 \delta_-
 =\frac1{a_+}\min\left\{\frac{a_-^2}{2},\,
                        \frac{c_-\pi^2}{4(1+2\varrho_+^2)}\right\}.
\]
From the table in Lemma~\ref{lem:effective-right-principal} (and
$e/(1+e)>0.731$ at $\log r=0.5$), we obtain:
\[
 \begin{array}{c|rrrr|r}
  \text{range of }\log r&a_-&a_+&\varrho_+&c_-
  &\delta_-\text{ greater than}\\\hline
  {[0,0.5]}&1.06975&1.377&0.931&0.5&0.32\\
  {[0.5,1]}&1.37622&1.733&0.812&0.731&0.44
 \end{array}
\]
By \eqref{eq:tails-general} with $S=\{\chi\in\mathbb R:|\chi|>\pi/2\}$,
the sum of the $m=-1$ ($[0,8]$) and $m=0$ ($[-8,0]$) integrals over
$|\chi|>\pi/2$, divided by $M_n(r)$, is bounded for $n\ge6400$ by
\begin{align*}
 109000\sqrt{1.377\sqrt n}\,e^{-0.32\sqrt n}&<0.00001
 &&\text{for }0\le\log r\le0.5,\\
 109000\sqrt{1.733\sqrt n}\,e^{-0.44\sqrt n}&<10^{-6}
 &&\text{for }0.5\le\log r<1.
\end{align*}

Summing the bounds on each range of $\log r$ shows that the
left-hand side in the lemma, divided by $M_n(r)$, is less than
$0.000013$ for $0\le\log r<1$, less than $0.0003$ for
$1\le\log r<5$, and less than $0.00105$ for $\log r\ge5$; all are
strictly less than $0.002$.
\end{proof}

\begin{proposition}\label{thm:effective-right}
If $n\ge6400$ is an integer and $1\le r\le e\sqrt n$, then
\[
 |D_n(-r^2)-M_n(r)\cos\Psi_n(r)|<0.95M_n(r).
\]
\end{proposition}
\begin{proof}
In \eqref{eq:ct-integral2} (with $a/\sqrt n$ in place of $u$), split
$[-\pi,\pi]$ into $|v|\le8a/\sqrt n$ and $8a/\sqrt n\le|v|\le\pi$.
On $|v|\le8a/\sqrt n$, write $\zeta=\frac{a}{\sqrt n}(1+i\eta)$
with $|\eta|\le8$ and substitute \eqref{eq:unfolded-theta} and
\eqref{eq:eta-transform}. The absolute integrability established in
Lemma~\ref{lem:effective-branch-tails} (for
$m\in\mathbb Z\setminus\{0,-1\}$, for $m=0$ on $[-8,0]$, and for
$m=-1$ on $[0,8]$), together with \eqref{eq:right-common-gaussian}
and \eqref{eq:right-bp-bounds} (for $m=0$ on $[0,8]$, and by
conjugation symmetry for $m=-1$ on $[-8,0]$), justifies interchanging
the sum over $m\in\mathbb Z$ with the integral over $|\eta|\le8$ by
Fubini's theorem. Thus the integral over $|v|\le8a/\sqrt n$ equals
the sum over $m\in\mathbb Z$ of \eqref{eq:right-branch-integral}
with $I=[-8,8]$.
Combining Lemma~\ref{lem:effective-right-principal},
Lemma~\ref{lem:effective-minor}, and
Lemma~\ref{lem:effective-branch-tails} therefore yields
\[
 |D_n(-r^2)-M_n(r)\cos\Psi_n(r)|
 <(0.88+0.05+0.002)M_n(r)
 <0.95M_n(r).
\]
\end{proof}

\section{Outer sign estimates and completion of the proof}\label{sec:outer}

We continue the sign estimates at larger radii. The identity
$p(n,d)=p_{2,d}(n-d^2)$ allows us to replace the restricted
coefficients by a compactly supported Bessel kernel. Its bilateral
Laplace transform is explicit, and Poisson summation isolates the
two frequencies $\pi$ and $-\pi$. We then bound the correction
caused by the restriction on the part sizes. At the remaining radii,
a single monomial determines the sign of $D_n(-r^2)$.

We first record two elementary bounds; the first is used below, and
both will be used in the final range of $r$. For $u>0$,
\[
 \log\prod_{\nu\ge1}(1-e^{-\nu u})^{-2}
 =2\sum_{\ell\ge1}\frac1{\ell(e^{\ell u}-1)}
 \le\frac2u\sum_{\ell\ge1}\frac1{\ell^2}
 =\frac{\pi^2}{3u}.
\]
Since $\sum_{m\ge0}p_2(m)e^{-mu}=\prod_{\nu\ge1}(1-e^{-\nu u})^{-2}$,
it follows that $p_2(t)\le\exp(ut+\pi^2/(3u))$ for every integer
$t\ge0$ and every $u>0$. Taking $u=\pi/\sqrt{3t}$ for $t\ge1$ (and
noting $p_2(0)=1$) gives
\begin{equation}\label{eq:outer-p2-bound}
 p_2(t)\le e^{2\pi\sqrt{t/3}}
 \qquad(t\ge0).
\end{equation}
For integers $k\ge0$ and $t\ge k+1$, a pair counted by
$p_2(t)-p_{2,k}(t)$ contains a part of size $j>k$ in at least one
colour. Removing one such part and applying \eqref{eq:outer-p2-bound}
gives
\begin{equation}\label{eq:outer-p2-diff}
 0\le p_2(t)-p_{2,k}(t)
 \le2\sum_{j=k+1}^{t}p_2(t-j)
 \le2t e^{2\pi\sqrt{(t-k-1)/3}}
 \qquad(t\ge k+1),
\end{equation}
whereas $p_2(t)-p_{2,k}(t)=0$ for $0\le t\le k$.

For $x\ge0$, \eqref{eq:prelim-bessel-series} and
\eqref{eq:prelim-bessel-integrals} give
$I_2(x)\le(x^2/8)I_0(x)\le(x^2/8)e^x$, since $(j+2)!\ge2j!$ for all
integers $j\ge0$. Using $|A_\ell^{(2)}(t)|\le\ell$ in
\eqref{eq:two-colour-rademacher} and $\sum_{\ell\ge2}\ell^{-2}<1$,
we obtain, for every integer $t\ge1$,
\begin{equation}\label{eq:outer-rademacher}
 \left|p_2(t)-
 \frac{\pi I_2\!\left(2\pi\sqrt{(t-1/12)/3}\right)}
      {6(t-1/12)}\right|
 \le\frac{\pi^3}{36}e^{\pi\sqrt{t/3}},
\end{equation}
because the sum over $\ell\ge2$ in \eqref{eq:two-colour-rademacher}
has modulus at most
\[
 \frac{\pi}{6(t-1/12)}
 \sum_{\ell\ge2}
 \frac{\pi^2(t-1/12)}{6\ell^2}
 e^{\pi\sqrt{(t-1/12)/3}}
 <\frac{\pi^3}{36}e^{\pi\sqrt{t/3}}.
\]

For $y>0$, define the continuous function $\mathcal K_y$ on
$\mathbb R$ by
\[
 \mathcal K_y(x)=
 \begin{cases}
 \displaystyle
 \frac{\pi I_2\!\left(\frac{2\pi}{\sqrt3}
                         \sqrt{y^2-x^2}\right)}
      {6(y^2-x^2)}-\frac{\pi^3}{36},&|x|<y,\\[6pt]
 0,&|x|\ge y.
 \end{cases}
\]
Continuity at $|x|=y$ follows from $I_2(w)=w^2/8+O(w^4)$ in
\eqref{eq:prelim-bessel-series}. The nonnegativity of the series
coefficients in \eqref{eq:prelim-bessel-series} and the bound
$I_2(w)\le(w^2/8)e^w$ ($w\ge0$) show that
$0\le\mathcal K_y(x)\le(\pi^3/36)e^{2\pi y/\sqrt3}$ for all
$x\in\mathbb R$. For $\Rea\zeta>0$, set
\[
 \mathcal T_y(\zeta)=\int_{-y}^{y}\mathcal K_y(x)e^{\zeta x}\,dx.
\]
We evaluate $\mathcal T_y(\zeta)$ explicitly. Termwise integration
of the series for $I_0$ in \eqref{eq:prelim-bessel-series} gives
\[
 \frac{I_2(2\pi\xi/\sqrt3)}{\xi^2}
 =\frac{2\pi^2}{3}\int_0^1 t(1-t^2)
                  I_0(2\pi\xi t/\sqrt3)\,dt
 \qquad(\xi\ge0),
\]
with the value at $\xi=0$ understood by continuity. Also, for
$\tau\ge0$ and $\zeta\in\mathbb C$,
\[
 \int_{-y}^{y}e^{\zeta x}I_0\!\left(\tau\sqrt{y^2-x^2}\right)\,dx
 =\frac{2\sinh\!\left(y\sqrt{\zeta^2+\tau^2}\right)}
        {\sqrt{\zeta^2+\tau^2}}.
\]
Indeed, for $\zeta\in\mathbb R$, integrating
$e^{y(\zeta u_1+\tau u_2)}$ over the unit sphere
$S^2\subset\mathbb R^3$ with coordinates $(u_1,u_2,u_3)$
 (first at fixed $u_1$ using
\eqref{eq:prelim-bessel-integrals}, and then after rotating
$(\zeta,\tau,0)$ to the first coordinate axis) gives the identity on
$\mathbb R$;  since both sides are entire in $\zeta$
(the right-hand side is an entire function of $\zeta^2+\tau^2$),
the identity holds for all $\zeta\in\mathbb C$. Combining these two
identities (and using $\int_0^1t(1-t^2)\,dt=1/4$) yields
\begin{equation}\label{eq:outer-transform-integral}
 \mathcal T_y(\zeta)=\frac{2\pi^3}{9}
 \int_0^1 t(1-t^2)
 \left\{
 \frac{\sinh\!\left(y\sqrt{\zeta^2+4\pi^2t^2/3}\right)}
      {\sqrt{\zeta^2+4\pi^2t^2/3}}
 -\frac{\sinh(y\zeta)}{\zeta}
 \right\}\,dt.
\end{equation}
Since $\Rea\zeta>0$, we have $\zeta^2+4\pi^2t^2/3\notin(-\infty,0]$
for all $t\in[0,1]$, so the principal square root
$z=\sqrt{\zeta^2+4\pi^2t^2/3}$ varies smoothly with $\Rea z>0$ from
$z=\zeta$ to $z=\sqrt{\zeta^2+4\pi^2/3}$. Making this substitution
in \eqref{eq:outer-transform-integral} gives
\begin{equation}\label{eq:outer-transform}
 \mathcal T_y(\zeta)+\frac{\pi^3\sinh(y\zeta)}{18\zeta}
 =\frac1{8\pi}
 \left[
 \frac{(\zeta^2+4\pi^2/3-z^2)\cosh(yz)}y
 +\frac{2z\sinh(yz)}{y^2}
 -\frac{2\cosh(yz)}{y^3}
 \right]_{z=\zeta}^{\,\sqrt{\zeta^2+4\pi^2/3}},
\end{equation}
because the derivative of the bracketed expression with respect to
$z$ equals $(\zeta^2+4\pi^2/3-z^2)\sinh(yz)$.

Fix $h>0$. For $v\in\mathbb R$, write
\[
 \rho(v)=\sqrt{(h+iv)^2+4\pi^2/3},\qquad
 \Delta(v)=(\Rea\rho(v))^2-h^2,
\]
with the principal square root (so $\Rea\rho(v)>0$). Comparing real
and imaginary parts of $\rho(v)^2=(h+iv)^2+4\pi^2/3$ gives
\begin{equation}\label{eq:outer-quadratic}
 \Delta(v)\left(h^2+v^2-\frac{4\pi^2}{3}+\Delta(v)\right)
 =\frac{4\pi^2h^2}{3},
 \qquad
 \Ima\rho(v)=\frac{hv}{\Rea\rho(v)}.
\end{equation}
From $2(\Rea\rho(v))^2=|\rho(v)^2|+\Rea(\rho(v)^2)$ we have
$2\Delta(v)=|\rho(v)^2|-(h^2+v^2-4\pi^2/3)$; since
$|\rho(v)^2|^2-(h^2+v^2-4\pi^2/3)^2=16\pi^2h^2/3>0$, it follows
that $|\rho(v)^2|>|h^2+v^2-4\pi^2/3|$, and hence $\Delta(v)>0$ and
$h^2+v^2-4\pi^2/3+2\Delta(v)=|\rho(v)^2|>0$. Together with
\eqref{eq:outer-quadratic}, this implies $0<\Delta(v)\le4\pi^2/3$.
Moreover, $\Delta(v)$ is a smooth function of $v^2$, and
differentiating \eqref{eq:outer-quadratic} with respect to $v^2$
yields
\begin{equation}\label{eq:outer-quadratic-derivative}
 -\frac{d\Delta}{d(v^2)}
 =\frac{\Delta}{h^2+v^2-4\pi^2/3+2\Delta}>0.
\end{equation}

In the estimates below for $D_n(-r^2)$ with $r>1$, we write
$h=2\log r>0$, and for $y>0$ and $h>0$ we define
\[
 \mathcal A_y(h)
 =\frac{h\exp\!\left(y\Rea\sqrt{(h+i\pi)^2+4\pi^2/3}\right)}{4\pi y^2}.
\]

\begin{proposition}\label{thm:effective-outer-poisson}
If $n\ge6401$ is an integer and
$\log(n-1/12)\le h=2\log r\le4(n-1/12)^{1/6}$, then
\[
 \left|D_n(-r^2)
 -\frac{\Rea\bigl(\rho(\pi)e^{\sqrt{n-1/12}\,\rho(\pi)}\bigr)}
       {4\pi(n-1/12)}\right|
 <0.99\mathcal A_{\sqrt{n-1/12}}(h).
\]
\end{proposition}
\begin{proof}
Set $y_0=\sqrt{n-1/12}$; since $n\ge6401$, we have $y_0^2>6400$ and
$y_0>80$, while $2\log y_0\le h\le4y_0^{1/3}$.
We first estimate $\mathcal T_y(h+iv)$ uniformly for
$9y_0/10\le y\le y_0$ and $v\in\pi+2\pi\mathbb Z$. Put
\[
 \mathcal M_y(v)=\frac{\rho(v)e^{y\rho(v)}}{8\pi y^2}.
\]
The parameter range gives
\[
 h>8.76,\qquad y>72,\qquad \frac{h}{y}<0.240,\qquad
 \frac{y}{h}\ge\frac{9y_0^{2/3}}{40},\qquad
 \frac{y}{h^3}\ge\frac{9}{640}.
\]
From \eqref{eq:outer-quadratic} and $0<\Delta(\pi)\le4\pi^2/3$, we
have
\[
 \Delta(\pi)\ge\frac{4\pi^2}{3(1+\pi^2/h^2)},\qquad
 \Rea\rho(\pi)+h\le h\left(1+\sqrt{1+\frac{4\pi^2}{3h^2}}\right).
\]
Since $h>8.76$, dividing these bounds gives
\[
 \Rea\rho(\pi)-h
 >0.425\,\frac{4\pi^2}{3h}=\frac{17\pi^2}{30h}.
\]
Subtracting \eqref{eq:outer-quadratic} at $\pi$ and at $v$ yields
\[
 \Delta(\pi)-\Delta(v)
 =\frac{(v^2-\pi^2)\Delta(\pi)}
       {h^2+v^2-4\pi^2/3+\Delta(\pi)+\Delta(v)}.
\]
For $v\in\pi+2\pi\mathbb Z$ with $\pi<|v|\le h$, we have
$|v|\ge3\pi$ (hence $h\ge3\pi$), $\Delta(\pi)+\Delta(v)\le8\pi^2/3$,
and $\Rea\rho(\pi)+\Rea\rho(v)\le2h\sqrt{1+4\pi^2/(3h^2)}$; thus
\[
 \Rea\rho(\pi)-\Rea\rho(v)
 \ge\frac{(4\pi^2/3)(v^2-\pi^2)}
         {2h^3\sqrt{1+4\pi^2/(3h^2)}
          (1+\pi^2/h^2)(1+v^2/h^2+4\pi^2/(3h^2))}
 >\frac{13\pi^2(v^2-\pi^2)}{50h^3},
\]
using $[2\sqrt{31/27}\,(10/9)(58/27)]^{-1}>0.195$.
By \eqref{eq:outer-quadratic} and $\Rea\rho(v)>h$, we have
$|\rho(v)|^2=h^2+\Delta(v)+h^2v^2/(\Rea\rho(v))^2\le h^2+v^2+4\pi^2/3$,
whence
\[
 \frac{|\rho(\pi)|}{h}
 \le\sqrt{1+\frac{7\pi^2}{3h^2}}<1.142,
 \qquad
 \frac{|\rho(v)|}{h}
 \le\sqrt{2+\frac{4}{27}}<1.47
 \qquad(3\pi\le|v|\le h).
\]
Consequently, $|\mathcal M_y(\pi)|+|\mathcal M_y(-\pi)|<1.142\mathcal A_y(h)$,
and for every integer $m\ge1$ with $(2m+1)\pi\le h$,
\[
 |\mathcal M_y((2m+1)\pi)|+|\mathcal M_y(-(2m+1)\pi)|
 \le1.47\mathcal A_y(h)e^{-\delta m(m+1)},
 \qquad
 \delta=\frac{26\pi^4y}{25h^3}.
\]

For $|v|\le h$, let $\zeta=h+iv$, so $h\le|\zeta|\le\sqrt2\,h$.
In \eqref{eq:outer-transform}, the $\cosh(yz)/y$ term vanishes at
the upper endpoint $z=\rho(v)$, and since
$|\rho(v)|\ge\Rea\rho(v)>h$ and $yh>1$, subtracting $\mathcal M_y(v)$
from the upper endpoint leaves modulus at most
$|\mathcal M_y(v)|\bigl(e^{-2yh}+(1+e^{-2yh})/(yh)\bigr)
\le|\mathcal M_y(v)|\bigl(1/(yh)+2e^{-2yh}\bigr)$.
The lower endpoint $z=\zeta$ and the subtracted term
$-\pi^3\sinh(y\zeta)/(18\zeta)$ in \eqref{eq:outer-transform} have
combined modulus at most
\[
 e^{yh}\left(\frac{\pi}{6y}+\frac{\sqrt2\,h}{4\pi y^2}
             +\frac1{4\pi y^3}+\frac{\pi^3}{18h}\right)
 <\frac{2}{h}e^{yh},
\]
using $h/y<0.240$, $y>72$, and $\pi^3/18<1.723$ after multiplying by
$h$. Hence
\[
 |\mathcal T_y(h+iv)-\mathcal M_y(v)|
 \le|\mathcal M_y(v)|\left(\frac1{yh}+2e^{-2yh}\right)
   +\frac2h e^{yh}
 \qquad(|v|\le h).
\]
There are at most $h/\pi+1$ frequencies $v\in\pi+2\pi\mathbb Z$
with $|v|\le h$, so summing $\frac2h e^{yh}$ over $|v|\le h$ and
dividing by $\mathcal A_y(h)$ gives at most
\[
 \frac{8y^2}{h}\left(1+\frac{\pi}{h}\right)e^{-17\pi^2y/(30h)}
 <\frac{14y^2}{h}e^{-17\pi^2y/(30h)}.
\]

For $|v|>h$ and $0\le\tau\le2\pi/\sqrt3$, applying the quadratic
identity \eqref{eq:outer-quadratic} with $\tau^2$ in place of
$4\pi^2/3$ gives
\[
 \Rea\sqrt{(h+iv)^2+\tau^2}
 \le h+\frac{4\pi^2/3}{2h\bigl(2-4\pi^2/(3h^2)\bigr)}.
\]
Thus $\bigl(\Rea\rho(\pi)-\Rea\sqrt{(h+iv)^2+\tau^2}\bigr)
\cdot\frac{3h}{4\pi^2}$ is strictly greater than
\[
 \frac1{(1+\pi^2/h^2)\bigl(1+\sqrt{1+4\pi^2/(3h^2)}\bigr)}
 -\frac1{2\bigl(2-4\pi^2/(3h^2)\bigr)}
 >\max\!\left\{0.152,\,\frac14-\frac{8.4}{h^2}\right\}.
\]
Indeed, the left-hand side increases with $h$ and exceeds $0.152054$
at $h=8.76$; for the second bound, using $(1+x)^{-1}\ge1-x$,
$(1+\sqrt{1+x})^{-1}\ge\frac12-\frac{x}{8}$, and
$h>8.76$ bounds the subtracted coefficient of $h^{-2}$ by
\[
 \frac{\pi^2}{2}+\frac{\pi^2}{6}
 +\frac{\pi^2}{6\bigl(1-2\pi^2/(3\cdot8.76^2)\bigr)}<8.4.
\]
Differentiating with respect to $\tau^2$ gives
\begin{equation}\label{eq:outer-transform-derivative}
 \frac{\partial}{\partial(\tau^2)}
 \frac{\sinh\!\left(y\sqrt{\zeta^2+\tau^2}\right)}
      {\sqrt{\zeta^2+\tau^2}}
 =\frac{y\cosh\!\left(y\sqrt{\zeta^2+\tau^2}\right)}
       {2(\zeta^2+\tau^2)}
  -\frac{\sinh\!\left(y\sqrt{\zeta^2+\tau^2}\right)}
        {2(\zeta^2+\tau^2)^{3/2}}.
\end{equation}
For $\zeta=h+iv$ with $|v|>h$, we have
$|\zeta^2+\tau^2|\ge h^2+v^2-4\pi^2/3>0.914(h^2+v^2)$ (since
$h^2+v^2>2\cdot8.76^2$). Together with $y>72$ and
$\sqrt{h^2+v^2}>\sqrt2\cdot8.76$, \eqref{eq:outer-transform-derivative}
has modulus at most
\[
 \frac{y}{h^2+v^2}
 \left(\frac1{2(0.914)}
       +\frac1{2(0.914)^{3/2}\,72\sqrt2\cdot8.76}\right)
 e^{y\Rea\sqrt{(h+iv)^2+\tau^2}}
 <\frac{0.55y}{h^2+v^2}e^{y\Rea\sqrt{(h+iv)^2+\tau^2}}.
\]
Writing the difference in \eqref{eq:outer-transform-integral} as
the integral of \eqref{eq:outer-transform-derivative} over
$\tau^2\in[0,4\pi^2t^2/3]$ and integrating in $t\in[0,1]$ (using
$\frac{2\pi^3}{9}\cdot\frac{4\pi^2}{3}\int_0^1t^3(1-t^2)\,dt
=\frac{\pi(4\pi^2/3)^2}{72}$) bounds $|\mathcal T_y(h+iv)|$ for
$|v|>h$. Since the positive frequencies in $\pi+2\pi\mathbb Z$ are
spaced by $2\pi$, comparison with $\int_h^\infty x^{-2}\,dx$ gives
\[
 \sum_{\substack{v\in\pi+2\pi\mathbb Z\\|v|>h}}\frac1{h^2+v^2}
 \le\frac2{h^2}+\frac1{\pi h}.
\]
Dividing by $\mathcal A_y(h)$ and using
$\frac{0.55\pi(4\pi^2/3)^2}{18}(1+2\pi/h)<29$ (as $h>8.76$), we
obtain
\[
 \frac1{\mathcal A_y(h)}
 \sum_{\substack{v\in\pi+2\pi\mathbb Z\\|v|>h}}|\mathcal T_y(h+iv)|
 <\frac{29y^3}{h^2}
  \exp\!\left(-\frac{4\pi^2y}{3h}
              \max\!\left\{0.152,\,\frac14-\frac{8.4}{h^2}\right\}\right).
\]

Now let $c\in\{9/10,1\}$ and suppose $cy_0\le y\le y_0$.
Since $\delta\ge13c\pi^4/800$ and $m(m+1)\ge2+4(m-1)$ for $m\ge1$,
\[
 \frac1{\mathcal A_y(h)}
 \sum_{\substack{v\in\pi+2\pi\mathbb Z\\\pi<|v|\le h}}|\mathcal M_y(v)|
 \le\frac{1.47e^{-13c\pi^4/400}}{1-e^{-13c\pi^4/200}}.
\]
Note that $\frac{4\pi^2y}{3h}\ge\frac{4\pi^2}{3}\cdot\frac{9y_0^{2/3}}{40}>54$.
For $h\le11$, the bound
$\frac{29y^3}{h^2}e^{-0.152\cdot4\pi^2y/(3h)}$ increases with $h$
and decreases with $y$ (its logarithmic derivatives have signs
$-2+0.152\frac{4\pi^2y}{3h}>0$ and
$3-0.152\frac{4\pi^2y}{3h}<0$), so setting $h=11$ and $y=cy_0$ gives
at most
\[
 \frac{29c^3y_0^3}{121}e^{-76\pi^2cy_0/4125}.
\]
For $h\ge11$, the bound
$\frac{29y^3}{h^2}\exp\bigl(-\frac{4\pi^2y}{3h}(\frac14-\frac{8.4}{h^2})\bigr)$
also increases with $h$ and decreases with $y$, since
$-2+\frac{4\pi^2y}{3h}(\frac14-\frac{25.2}{h^2})>-2+54(0.0417)>0$
and $\frac{4\pi^2y}{3h}(\frac14-\frac{8.4}{h^2})>54(0.180)>3$;
setting $h=4y_0^{1/3}$ and $y=cy_0$ gives at most
\[
 \frac{29c^3y_0^{7/3}}{16}
 \exp\!\left(-\frac{c\pi^2y_0^{2/3}}{12}+\frac{7c\pi^2}{40}\right).
\]
Similarly, $\frac{14y^2}{h}e^{-17\pi^2y/(30h)}$ increases with $h$
and decreases with $y$, so setting $h=4y_0^{1/3}$ and $y=cy_0$ gives
at most
\[
 \frac{7c^2y_0^{5/3}}{2}e^{-17c\pi^2y_0^{2/3}/120}.
\]
Each of the preceding three functions of $y_0$ decreases for
$y_0\ge80$: differentiating their logarithms with respect to $y_0$
reduces this assertion, respectively, to
\[
 \frac{76\pi^2cy_0}{4125}>3,\qquad
 \frac{c\pi^2y_0^{2/3}}{12}>\frac72,\qquad
 \frac{17c\pi^2y_0^{2/3}}{120}>\frac52,
\]
which hold at $y_0=80$, $c=9/10$ and increase with $y_0$ and $c$.
Finally, since $yh>80c\cdot8.76$, summing
$|\mathcal M_y(v)|\bigl(1/(yh)+2e^{-2yh}\bigr)/\mathcal A_y(h)$ over
$|v|\le h$ gives at most
\[
 \left(1.142+\frac{1.47e^{-13c\pi^4/400}}{1-e^{-13c\pi^4/200}}\right)
 \left(\frac1{80c\cdot8.76}+2e^{-160c\cdot8.76}\right).
\]
Evaluating these five explicit upper bounds at $y_0=80$ gives the
strict upper bounds in the following table:
\[
\begin{array}{c|cc}
 &c=9/10\ (y\ge9y_0/10)&c=1\ (y=y_0)\\\hline
 \pi<|v|\le h\ (\mathcal M_y\text{ sum})&0.085385&0.062117\\
 |v|>h,\ h\le11&0.184340&0.059036\\
 |v|>h,\ h\ge11&0.185375&0.065639\\
 |v|\le h\text{ (lower-endpoint bound)}&4\cdot10^{-7}&4\cdot10^{-7}\\
 |v|\le h\text{ (upper-endpoint error)}&0.001947&0.001719
\end{array}
\]
Here the two lower-endpoint bounds are less than
$3.005\cdot10^{-7}$ and $2.767\cdot10^{-8}$, respectively.
All numerical inequalities in this table follow from
$3.14159265358979<\pi<3.14159265358980$ and
$18.5663553344511<80^{2/3}<18.5663553344512$ (verified by cubing),
together with the Taylor series for $e^x$ (with remainder after
degree $J$ bounded by the first omitted term divided by
$1-x/(J+2)$ for $0<x<J+2$, and reciprocation for negative
arguments). Summing the relevant rows (taking the larger of the two
$|v|>h$ rows) gives
\begin{align}
 &\sum_{v\in\pi+2\pi\mathbb Z}|\mathcal T_y(h+iv)|
  <1.42\mathcal A_y(h)
  \qquad(9y_0/10\le y\le y_0),\label{eq:effective-outer-sum}\\
 &\left|\sum_{v\in\pi+2\pi\mathbb Z}\mathcal T_{y_0}(h+iv)
        -2\Rea\mathcal M_{y_0}(\pi)\right|
  <0.130\mathcal A_{y_0}(h),\label{eq:effective-outer-pair}
\end{align}
since $1.142+0.085385+0.185375+0.001947+0.0000004<1.42$ and
$0.062117+0.065639+0.001719+0.0000004<0.130$.

We now transfer these estimates to $D_n(-r^2)$. From
$p(n,d)=p_{2,d}(n-d^2)
=[q^{n-d^2}]\prod_{\nu\ge1}(1-q^\nu)^{-2}\prod_{\nu>d}(1-q^\nu)^2$
and Euler's expansion
\[
 \prod_{\nu>d}(1-q^\nu)^2
 =\sum_{\alpha,\beta\ge0}
  \frac{(-1)^{\alpha+\beta}
        q^{(\alpha+\beta)d+[\alpha(\alpha+1)+\beta(\beta+1)]/2}}
       {(q;q)_\alpha(q;q)_\beta},
\]
set, for integers $\alpha,\beta,j\ge0$,
\[
 k=\alpha+\beta,\qquad
 P_{\alpha,\beta}(j)=[q^j]\bigl((q;q)_\alpha(q;q)_\beta\bigr)^{-1},
 \qquad
 \gamma_{\alpha,\beta,j}=j+\frac{(\alpha-\beta)^2+2k}{4}.
\]
Since $d^2+kd+\frac12(\alpha(\alpha+1)+\beta(\beta+1))+j
=(d+k/2)^2+\gamma_{\alpha,\beta,j}$, we obtain
\[
 D_n(-r^2)
 =\sum_{\alpha,\beta,j\ge0}(-1)^k P_{\alpha,\beta}(j)
  \sum_{d\ge1}(-e^h)^d
  p_2\!\left(n-\gamma_{\alpha,\beta,j}-(d+k/2)^2\right),
\]
where $p_2(t)$ is understood to be $0$ for $t<0$.
Every argument of $p_2$ is an integer, and only finitely many terms
are nonzero. Any triple $(\alpha,\beta,j)$ with a nonzero term
satisfies $\gamma_{\alpha,\beta,j}\le n-(1+k/2)^2\le n-1<y_0^2$, so
$y=\sqrt{y_0^2-\gamma_{\alpha,\beta,j}}>0$ is well defined.

First consider the triples with $y<9y_0/10$, that is,
$\gamma_{\alpha,\beta,j}>0.19y_0^2$. Writing $hd=h(d+k/2)-hk/2$
and applying \eqref{eq:outer-p2-bound}, Cauchy--Schwarz, and
concavity of the square root gives, for every integer $d\ge1$ with
$n-\gamma_{\alpha,\beta,j}-(d+k/2)^2\ge0$,
\[
 \begin{aligned}
 hd+\frac{2\pi}{\sqrt3}\sqrt{n-\gamma_{\alpha,\beta,j}-(d+k/2)^2}
 &\le\rho(0)\sqrt{n-\gamma_{\alpha,\beta,j}}-\frac{hk}{2}\\
 &\le\sqrt n\,\rho(0)
     -\frac{\rho(0)\gamma_{\alpha,\beta,j}}{2\sqrt n}-\frac{hk}{2}.
 \end{aligned}
\]
Splitting $-\frac{\rho(0)\gamma_{\alpha,\beta,j}}{2\sqrt n}$ into
two equal halves and using $\gamma_{\alpha,\beta,j}>0.19y_0^2$ in
one half and $\gamma_{\alpha,\beta,j}\ge j$ in the other, noting
that $\sum_{j\ge0}P_{\alpha,\beta}(j)x^j=((x;x)_\alpha(x;x)_\beta)^{-1}
\le(1-x)^{-k}$ for $0<x<1$, and observing that there are at most
$\sqrt n$ positive integers $d$ with $n-(d+k/2)^2\ge0$, we see that
the sum of the moduli of all terms with $y<9y_0/10$ is at most
\[
 \frac{\sqrt n}{(1-z_1)^2}
 \exp\!\left(\sqrt n\,\rho(0)-\frac{0.19\rho(0)y_0^2}{4\sqrt n}\right),
 \qquad
 z_1=\frac{e^{-h/2}}{1-e^{-\rho(0)/(4\sqrt n)}}.
\]
Using $e^{-h/2}\le y_0^{-1}$, $(1-e^{-x})^{-1}\le1+x^{-1}$ ($x>0$),
$y_0>80$, $\rho(0)>h>8.76$, and $\sqrt n/y_0<1.00001$ gives
\[
 z_1\le\frac1{y_0}+\frac{4\sqrt n}{y_0\rho(0)}<0.47,
 \qquad
 \frac{\sqrt n}{(1-z_1)^2}<3.6y_0.
\]
Since $0<\Delta(v)\le4\pi^2/3$ and $h>2\pi/\sqrt3$,
\eqref{eq:outer-quadratic-derivative} yields
$-\frac{d\Delta}{d(v^2)}\le\frac{4\pi^2/3}{h^2+4\pi^2/3}$; dividing
by $2\Rea\rho(v)>2h$ and integrating in $v^2$ from $0$ to $\pi^2$
gives
\[
 \rho(0)-\Rea\rho(\pi)
 \le\frac{(4\pi^2/3)\pi^2}{2h(h^2+4\pi^2/3)}
 <\frac{2\pi^4}{3h^3}.
\]
Also, $(\sqrt n-y_0)\rho(0)<\frac{\rho(0)}{24y_0}<0.010$ and
$y_0/\sqrt n>0.99999$. Since $\rho(0)>h$ and
\[
 \frac{0.19}{4}(0.99999)-\frac{2\pi^4}{3(8.76)^4}>0.036,
\]
the exponent minus $y_0\Rea\rho(\pi)$ is less than
$-0.036hy_0+0.010<-0.315y_0+0.010$. Dividing by $\mathcal A_{y_0}(h)$
(and using $3.6\cdot4\pi/8.76<15\pi/8$) bounds the sum of all terms
with $y<9y_0/10$, divided by $\mathcal A_{y_0}(h)$, by
\[
 \frac{15\pi y_0^3}{8}e^{-0.315y_0+0.010}<0.00004,
\]
which decreases for $y_0\ge80$ and is less than $0.000034638$ at
$y_0=80$.

For each triple with $y\ge9y_0/10$, replace
$p_2(n-\gamma_{\alpha,\beta,j}-(d+k/2)^2)$ by $\mathcal K_y(d+k/2)$.
Put $t=n-\gamma_{\alpha,\beta,j}-(d+k/2)^2$. If $t<0$, then
$(d+k/2)^2>y^2$, so both $p_2(t)$ and $\mathcal K_y(d+k/2)$ vanish.
If $t=0$, then $\mathcal K_y(d+k/2)=0$ and
$|p_2(0)-0|=1<\pi^3/18$. If $t\ge1$, then
$y^2-(d+k/2)^2=t-1/12$, so \eqref{eq:outer-rademacher} and the
subtracted constant $\pi^3/36$ in $\mathcal K_y$ give
$|p_2(t)-\mathcal K_y(d+k/2)|\le\frac{\pi^3}{36}e^{\pi\sqrt{t/3}}
+\frac{\pi^3}{36}\le\frac{\pi^3}{18}e^{\pi\sqrt{t/3}}$.
Put $u_-=\sqrt{h^2+\pi^2/3}$. Since there are at most $\sqrt n$
positive integers $d$ with $t\ge0$ and
$\frac{\pi^3}{18}\sqrt n<2y_0$, Cauchy--Schwarz bounds the
replacement error for one triple by
\[
 2y_0\exp\!\left(u_-\sqrt{n-\gamma_{\alpha,\beta,j}}-\frac{hk}{2}\right).
\]
By concavity and $\gamma_{\alpha,\beta,j}\ge j$,
\[
 \sqrt{n-\gamma_{\alpha,\beta,j}}
 \le\sqrt n-\frac{j}{2\sqrt n}
 \le y_0+\frac1{24y_0}-\frac{j}{3y_0}.
\]
Multiplying by $P_{\alpha,\beta}(j)$, summing over $\alpha,\beta,j\ge0$,
and using
\[
 z_2:=\frac{e^{-h/2}}{1-e^{-u_-/(3y_0)}}
 \le\frac1{y_0}+\frac3{u_-}
 <\frac1{80}+\frac3{8.76}<0.355
\]
(so $2(1-z_2)^{-2}<5$) bounds the total replacement error over all
triples with $y\ge9y_0/10$ by
$5y_0\exp\bigl((y_0+\frac1{24y_0})u_-\bigr)$.
Furthermore,
\[
 \Rea\rho(\pi)-u_-
 >\frac{17\pi^2}{30h}-\frac{\pi^2}{6h}=\frac{2\pi^2}{5h},
 \qquad
 \frac{u_-}{24y_0}<0.010.
\]
Dividing by $\mathcal A_{y_0}(h)$ yields at most
$\frac{20\pi y_0^3}{h}e^{-2\pi^2y_0/(5h)+0.010}$. Since
$\frac{2\pi^2y_0}{5h}>1$, this expression increases with $h$;
evaluating at $h=4y_0^{1/3}$ bounds the total replacement error,
divided by $\mathcal A_{y_0}(h)$, by
\[
 5\pi y_0^{8/3}e^{-\pi^2y_0^{2/3}/10+0.010}<0.0208,
\]
which decreases for $y_0\ge80$ and is less than $0.020761$ at
$y_0=80$.

Next, for each triple with $y\ge9y_0/10$, extend
$\sum_{d\ge1}(-e^h)^d\mathcal K_y(d+k/2)$ to all $d\in\mathbb Z$.
Only $d\le0$ are added, so $e^{hd}\le1$, and $\mathcal K_y(d+k/2)=0$
unless $|d+k/2|<y\le y_0$ (at most $2y_0+1$ values of $d\le0$).
The condition $\gamma_{\alpha,\beta,j}<y_0^2$ implies
$\alpha,\beta<2y_0^2$ and $j<y_0^2$, giving at most
$(2y_0^2+1)^2(y_0^2+1)\le18y_0^6$ triples. For each triple,
$P_{\alpha,\beta}(j)\le p_2(j)\le e^{2\pi y_0/\sqrt3}$ (by
\eqref{eq:outer-p2-bound}) and
$\mathcal K_y(d+k/2)\le(\pi^3/36)e^{2\pi y_0/\sqrt3}$. Since
$18y_0^6(2y_0+1)(\pi^3/36)<40y_0^7$ and $\Rea\rho(\pi)>h\ge2\log y_0$,
the sum of the moduli of all added terms with $d\le0$, divided by
$\mathcal A_{y_0}(h)$, is less than
\[
 \frac{160\pi y_0^9}{2\log y_0}
 \exp\!\left(y_0\left(\frac{4\pi}{\sqrt3}-2\log y_0\right)\right)
 <10^{-30},
\]
which decreases for $y_0\ge80$ and is less than $2.907\cdot10^{-34}$
at $y_0=80$. Together, the small-$y$ terms, the replacement error,
and the $d\le0$ extension error have combined modulus less than
$0.022\mathcal A_{y_0}(h)$.

For each triple with $y\ge9y_0/10$, applying Poisson summation to
the compactly supported, continuous, piecewise smooth function
$x\mapsto e^{hx}\mathcal K_y(x)$ on $\mathbb Z+k/2$ gives
\[
 \sum_{d\in\mathbb Z}(-e^h)^d\mathcal K_y(d+k/2)
 =e^{-hk/2}\sum_{v\in\pi+2\pi\mathbb Z}
           e^{-ivk/2}\mathcal T_y(h+iv);
\]
the absolute convergence $\mathcal T_y(h+iv)=O(v^{-2})$ (from
\eqref{eq:outer-transform-integral}--\eqref{eq:outer-transform-derivative})
ensures pointwise validity everywhere, including when $\mathbb Z+k/2$
contains $\pm y$.
For $(\alpha,\beta,j)=(0,0,0)$, we have $k=0$, $P_{0,0}(0)=1$, and
$y=y_0$, so \eqref{eq:effective-outer-pair} gives
\[
 \left|\sum_{d\in\mathbb Z}(-e^h)^d\mathcal K_{y_0}(d)
       -2\Rea\mathcal M_{y_0}(\pi)\right|
 <0.130\mathcal A_{y_0}(h).
\]

For the remaining triples $(\alpha,\beta,j)\ne(0,0,0)$ with
$y\ge9y_0/10$, put
\[
 \sigma=\frac{\Rea\rho(\pi)-4/y_0}{2y_0}>0.
\]
Since $y\ge9y_0/10>y_0/2$, we have
$2\log(y_0/y)\le2(y_0-y)/y\le4(y_0-y)/y_0$ and
$y_0-y=\gamma_{\alpha,\beta,j}/(y_0+y)\ge\gamma_{\alpha,\beta,j}/(2y_0)$;
hence
\[
 \frac{\mathcal A_y(h)}{\mathcal A_{y_0}(h)}
 \le e^{-(y_0-y)(\Rea\rho(\pi)-4/y_0)}
 \le e^{-\sigma\gamma_{\alpha,\beta,j}}.
\]
Write $w=e^{-h/2}$ and $q_0=e^{-\sigma}$. Since
$\gamma_{\alpha,\beta,j}\ge j+k/2$, summing
$P_{\alpha,\beta}(j)w^k q_0^{j+k/2}$ over all $\alpha,\beta,j\ge0$
(including $(0,0,0)$) yields
\[
 \sum_{\alpha,\beta,j\ge0}
 P_{\alpha,\beta}(j)w^{\alpha+\beta}q_0^{j+(\alpha+\beta)/2}
 =\left(\sum_{\alpha\ge0}
        \frac{(wq_0^{1/2})^\alpha}{(q_0;q_0)_\alpha}\right)^2
 =\prod_{m\ge0}(1-wq_0^{m+1/2})^{-2}.
\]
The last equality is Euler's identity
\cite[Corollary~2.2, Eq.~(2.2.5)]{Andrews}.
Expanding the logarithm of the product, using $\sinh t\ge t$ ($t>0$),
$\Rea\rho(\pi)>h$, $w\le y_0^{-1}$, and
$\Li_2(y_0^{-1})\le y_0^{-1}\bigl(1+\frac1{4(y_0-1)}\bigr)$, gives
\[
 \begin{aligned}
 \log\prod_{m\ge0}(1-wq_0^{m+1/2})^{-2}
 &=\sum_{\ell\ge1}\frac{w^\ell}{\ell\sinh(\sigma\ell/2)}
  \le\frac2{\sigma}\Li_2(w)\\
 &\le\frac{4}{h-4/y_0}\left(1+\frac1{4(y_0-1)}\right)
  <0.461,
 \end{aligned}
\]
since the last rational expression decreases with $y_0$ and $h$, and
is less than $0.460482$ at $y_0=80$, $h=2\log80$.
Subtracting the $(0,0,0)$ term and applying
\eqref{eq:effective-outer-sum} bounds the sum over all triples
$(\alpha,\beta,j)\ne(0,0,0)$ with $y\ge9y_0/10$ by
$1.42(e^{0.461}-1)\mathcal A_{y_0}(h)$.

Combining all estimates gives
\[
 \left|D_n(-r^2)-2\Rea\mathcal M_{y_0}(\pi)\right|
 <\bigl(1.42(e^{0.461}-1)+0.130+0.022\bigr)\mathcal A_{y_0}(h)
 <0.984\mathcal A_{y_0}(h)
 <0.99\mathcal A_{y_0}(h).
\]
Since $2\Rea\mathcal M_{y_0}(\pi)
=\Rea\bigl(\rho(\pi)e^{y_0\rho(\pi)}\bigr)/(4\pi y_0^2)$, this
proves the proposition.
\end{proof}

We next bound the two-colour partition numbers. These bounds will
show that one monomial determines the sign of $D_n(-r^2)$ at the
remaining radii.

\begin{lemma}\label{lem:two-colour-uniform}
For $y>0$, let
\[
 B(y)=\frac{e^{2\pi y/\sqrt3}}{4\cdot3^{3/4}y^{5/2}}.
\]
If $t\ge1$ is an integer, then $p_2(t)\le2.6B(\sqrt{t-1/12})$.
For $t\ge25$ this improves to $p_2(t)\le1.1B(\sqrt{t-1/12})$,
and for $t\ge100$ one also has $p_2(t)\ge0.9B(\sqrt{t-1/12})$.
\end{lemma}
\begin{proof}
We first show that for every $x>0$,
\[
 \frac{e^x}{\sqrt{2\pi x}}\left(1-\frac{15}{8x}\right)
 \le I_2(x)\le\frac{e^x}{\sqrt{2\pi x}}.
\]
Indeed, substituting $1-u/x$ for the integration variable $y$ in the
integral representation of $I_2(x)$ in
\eqref{eq:prelim-bessel-integrals} shows that the ratio of $I_2(x)$
to $e^x/\sqrt{2\pi x}$ equals
\[
 \frac4{3\sqrt\pi}\int_0^{2x}
 e^{-u}u^{3/2}\left(1-\frac{u}{2x}\right)^{3/2}\,du.
\]
For $0\le u\le2x$, convexity of $w\mapsto(1-w)^{3/2}$ on $[0,1]$
gives $1-3u/(4x)\le(1-u/(2x))^{3/2}\le1$. Extending the resulting
integrals from $(0,2x)$ to $(0,\infty)$ (which increases the upper
bound since $e^{-u}u^{3/2}>0$, and decreases the lower bound since
$1-3u/(4x)<0$ for $u>2x$) and evaluating
$\int_0^\infty e^{-u}u^{3/2}\,du=\frac34\sqrt\pi$ and
$\int_0^\infty e^{-u}u^{5/2}\,du=\frac{15}{8}\sqrt\pi$ yields the
claimed bounds for $I_2(x)$.

Next, for $y>0$, define
\[
 R(y)=\frac{\pi^3}{3^{5/4}}y^{5/2}
      \exp\!\left(\frac{\pi}{\sqrt3}\bigl(\sqrt{y^2+1/12}-2y\bigr)\right).
\]
For $y\ge\sqrt{11/12}$, we have
$\sqrt{y^2+1/12}-y=\frac{1/12}{\sqrt{y^2+1/12}+y}<\frac1{24y}
\le\frac1{24\sqrt{11/12}}$, so
$\frac{\pi}{\sqrt3}(\sqrt{y^2+1/12}-y)<\frac{\pi}{12\sqrt{11}}$.
Since $y^{5/2}e^{-\pi y/\sqrt3}$ attains its maximum on
$(0,\infty)$ at $y=5\sqrt3/(2\pi)$, it follows that
\[
 R(y)<\frac{\pi^3}{3^{5/4}}e^{\pi/(12\sqrt{11})}
      \left(\frac{5\sqrt3}{2\pi}\right)^{5/2}e^{-5/2}
     <1.556
 \qquad(y\ge\sqrt{11/12}).
\]
Moreover, for all $y>0$, logarithmic differentiation gives
\[
 \frac{R'(y)}{R(y)}
 =\frac{5}{2y}+\frac{\pi}{\sqrt3}\left(\frac{y}{\sqrt{y^2+1/12}}-2\right)
 <\frac{5}{2y}-\frac{\pi}{\sqrt3},
\]
which is negative for $y>5\sqrt3/(2\pi)$ and in particular for
$y\ge\sqrt{25-1/12}$. Hence $R(y)$ decreases for
$y\ge\sqrt{25-1/12}$.

Now let $t\ge1$ be an integer and put $y=\sqrt{t-1/12}$, so that
$y\ge\sqrt{11/12}$. Applying the bounds for $I_2(x)$ at
$x=2\pi y/\sqrt3$ gives
\[
 B(y)\left(1-\frac{15\sqrt3}{16\pi y}\right)
 \le\frac{\pi I_2(2\pi y/\sqrt3)}{6y^2}
 \le B(y),
\]
and \eqref{eq:outer-rademacher} gives
\[
 \left|p_2(t)-\frac{\pi I_2(2\pi y/\sqrt3)}{6y^2}\right|
 \le B(y)R(y).
\]
Since $R(y)<1.556$, we obtain $p_2(t)<(1+1.556)B(y)<2.6B(y)$ for
every integer $t\ge1$. For $t\ge25$, monotonicity of $R$ gives
$R(y)\le R(\sqrt{25-1/12})<0.052$, whence
$p_2(t)<(1+0.052)B(y)<1.1B(y)$. For $t\ge100$, we have
$R(y)\le R(\sqrt{100-1/12})<0.00004$ and
\[
15\sqrt3/(16\pi y)\le15\sqrt3/\bigl(16\pi\sqrt{100-1/12}\bigr)<0.052,
\]
whence $p_2(t)>(1-0.052-0.00004)B(y)>0.9B(y)$.
\end{proof}

We use the function $B$ from Lemma~\ref{lem:two-colour-uniform} to
compare the monomials of $D_n$.

\begin{lemma}\label{lem:effective-outer-monomial}
If $n\ge6401$ is an integer, $m=\lfloor\sqrt n\rfloor$, and $d$ is
an integer with $1\le d\le m-1$ such that
\[
 h:=\frac{2\pi d}{\sqrt3\sqrt{n-1/12-d^2}}\ge4(n-1/12)^{1/6},
\]
then
\[
 p(n,d)e^{dh}>\sum_{\substack{1\le\ell\le m\\\ell\ne d}}p(n,\ell)e^{\ell h}.
\]
\end{lemma}
\begin{proof}
Set $y_0=\sqrt{n-1/12}>80$, $j=m-d$ (so $1\le j\le m-1$), and
$f_j=\sqrt{y_0^2-d^2}$. Since
$\sqrt n-y_0=\frac{1/12}{\sqrt n+y_0}<\frac1{24y_0}<1$, we have
$d\le m-1\le\sqrt n-1<y_0$ and $m>\sqrt n-1>y_0-1$. Together with
$h\ge4y_0^{1/3}$ and $y_0>80$, these imply
\[
 \begin{gathered}
 f_j^2\ge y_0^2-(m-1)^2\ge2m-\frac{13}{12}>2y_0-\frac{37}{12}>1.96y_0,
 \qquad f_j=\frac{2\pi d}{\sqrt3 h}<\frac{\pi}{2\sqrt3}y_0^{2/3},\\
 d>0.978y_0,
 \qquad j<0.42y_0^{1/3}.
 \end{gathered}
\]
Indeed, $f_j/d=\frac{2\pi}{\sqrt3 h}\le\frac{\pi}{2\sqrt3 y_0^{1/3}}$
gives $d/y_0=(1+f_j^2/d^2)^{-1/2}\ge(1+\pi^2/(12y_0^{2/3}))^{-1/2}>0.978$,
and
$j=(m-y_0)+(y_0-d)<\frac1{12y_0}+\frac{f_j^2}{y_0+d}<0.42y_0^{1/3}$.

We first bound $p(n,d)$ from below. Set $t=n-d^2-d-1$. Since
$d\le m-1$, we have $t\ge m-1>y_0-2>78$; moreover,
$d+1\le m\le2(m-13/12)\le2(t-1/12)$ implies
$f_j^2=(t-1/12)+(d+1)\le3(t-1/12)$. By \eqref{eq:outer-p2-diff},
\[
 0\le p_2(n-d^2)-p(n,d)\le2\sum_{\nu=0}^t p_2(\nu).
\]
For $25\le x\le t$, logarithmic differentiation gives
\[
 \frac{d}{dx}\log B(\sqrt{x-1/12})
 =\frac{\pi}{\sqrt{3(x-1/12)}}-\frac5{4(x-1/12)}
 \ge\frac{2\pi}{3\sqrt3\sqrt{t-1/12}}=:\alpha_0,
\]
since $\frac{5}{4\sqrt{25-1/12}}<\frac{\pi}{3\sqrt3}$. Integrating
from $\nu$ to $t$ yields
$B(\sqrt{\nu-1/12})\le B(\sqrt{t-1/12})e^{-\alpha_0(t-\nu)}$ for
$25\le\nu\le t$. We also have
\[
 \frac{2\pi}{\sqrt3}\bigl(f_j-\sqrt{t-1/12}\bigr)
 =\frac{2\pi}{\sqrt3}\frac{d+1}{f_j+\sqrt{t-1/12}}
 >\frac{\pi d}{\sqrt3 f_j}=\frac h2
\]
and
$\alpha_0^{-1}=\frac{3\sqrt3\sqrt{t-1/12}}{2\pi}<\frac{3\sqrt3 f_j}{2\pi}
=\frac{3d}{h}<\frac{3y_0}{h}$. Applying
Lemma~\ref{lem:two-colour-uniform}, summing the geometric series
with $(1-e^{-\alpha_0})^{-1}\le1+\alpha_0^{-1}$, using
$f_j/\sqrt{t-1/12}\le\sqrt3$, and dividing by $0.9B(f_j)$ therefore
gives
\[
 \frac2{0.9B(f_j)}\sum_{\nu=25}^t p_2(\nu)
 \le\frac{2.2}{0.9}3^{5/4}\left(1+\frac{3y_0}{h}\right)e^{-h/2}
 \le\frac{2.2}{0.9}3^{5/4}\left(1+\frac34y_0^{2/3}\right)e^{-2y_0^{1/3}}
 <0.026053.
\]
For $0\le\nu<25$, \eqref{eq:outer-p2-bound} gives
$p_2(\nu)\le e^{10\pi/\sqrt3}$, whence
\[
 \frac2{0.9B(f_j)}\sum_{\nu=0}^{24}p_2(\nu)
 \le\frac{200\cdot3^{3/4}}{0.9}f_j^{5/2}e^{2\pi(5-f_j)/\sqrt3}
 <4\cdot10^{-7}.
\]
Both upper bounds decrease for $y_0\ge80$: the first has
logarithmic derivative $\frac{(3/2)w}{1+(3/4)w^2}-2<0$ with respect
to $w=y_0^{1/3}$, while the second decreases with
$f_j>\sqrt{1.96y_0}>12.5$ since $\frac5{2f_j}-\frac{2\pi}{\sqrt3}<0$;
evaluating the first at $y_0=80$ and the second at
$f_j=\sqrt{1.96\cdot80}$ gives the displayed constants. Since
$n-d^2=f_j^2+1/12>100$, the lower bound in
Lemma~\ref{lem:two-colour-uniform} yields
\[
 p(n,d)\ge0.9(1-0.026053-4\cdot10^{-7})B(f_j)>0.875B(f_j).
\]

Now index the remaining terms $\ell\in\{1,\ldots,m\}\setminus\{d\}$
by $k=m-\ell$, so $0\le k\le m-1$ and $k\ne j$. Put
$t_k=n-(m-k)^2$; whenever $t_k\ge1$, put
$f_k=\sqrt{t_k-1/12}=\sqrt{y_0^2-(m-k)^2}$. Using
$h=\frac{2\pi d}{\sqrt3 f_j}$ and $f_k^2-f_j^2=2d(k-j)-(k-j)^2$,
we obtain the exact identity
\[
 H(k):=hd+\frac{2\pi f_j}{\sqrt3}-h(m-k)-\frac{2\pi f_k}{\sqrt3}
 =\frac{\pi}{\sqrt3 f_j}\bigl((k-j)^2+(f_k-f_j)^2\bigr).
\]
Since $p(n,m-k)=p_{2,m-k}(t_k)\le p_2(t_k)$,
Lemma~\ref{lem:two-colour-uniform} shows that for $t_k\ge25$,
\begin{equation}\label{eq:outer-monomial-ratio}
 \frac{p(n,m-k)e^{(m-k)h}}{p(n,d)e^{dh}}
 \le\frac{1.1}{0.875}\left(\frac{f_j}{f_k}\right)^{5/2}e^{-H(k)}
 <1.26\left(\frac{f_j}{f_k}\right)^{5/2}e^{-H(k)}.
\end{equation}
We sum the ratios $\frac{p(n,m-k)e^{(m-k)h}}{p(n,d)e^{dh}}$ in three
ranges of $k$.

For $k>j$, we have $t_k>n-d^2>100>25$ and $f_k>f_j$. Extend
$f(u)=\sqrt{y_0^2-(m-u)^2}$ and
$H(u)=h(u-j)+\frac{2\pi}{\sqrt3}(f_j-f(u))$ to $u\in(m-y_0,m+y_0)$.
Then $H(j)=H'(j)=0$,
\[
 H''(u)=\frac{2\pi y_0^2}{\sqrt3 f(u)^3}>0,
\]
and for $0\le w\le1$ (when $j+1\le m-1$),
$f(j+w)^2-f_j^2=w(2d-w)\le2y_0w\le\frac{50}{49}wf_j^2$.  Using
$f_j<\frac{\pi}{2\sqrt3}y_0^{2/3}$ (which gives
$\frac{2\pi y_0^2}{\sqrt3 f_j^3}>\frac{48}{\pi^2}$) and Taylor's
formula with integral remainder, we find that whenever $j<m-1$,
\[
 H(j+1)=\int_0^1(1-w)H''(j+w)\,dw
 \ge\frac{2\pi y_0^2}{\sqrt3 f_j^3}
    \int_0^1\frac{1-w}{(1+(50/49)w)^{3/2}}\,dw
 >\frac{16}{\pi^2},
\]
since the integral equals $2/(1+\sqrt{99/49})^2>1/3$. By convexity
of $H$, we have $H(k)\ge(k-j)H(j+1)$ for all $k>j$; together with
$f_k>f_j$ and \eqref{eq:outer-monomial-ratio}, this bounds the sum
over $k>j$ by
\[
 \frac{1.26}{e^{16/\pi^2}-1}<0.311.
\]

For $k<j$ with $t_k\ge1$ and $f_k\ge f_j/2$, we have
$t_k=f_k^2+1/12>0.49y_0+1/12>25$, so \eqref{eq:outer-monomial-ratio}
applies. Put $\xi_k=f_k/f_j\in[1/2,1)$. Concavity of $f$ on
$(m-y_0,m+y_0)$ gives $f_k\le f_j-\frac{d}{f_j}(j-k)$, whence
$\frac{2\pi f_j}{\sqrt3}(1-\xi_k)\ge h(j-k)$, while
$\xi_k\ge1/2$ gives $\frac52\log(1/\xi_k)\le5(1-\xi_k)$. From
$H(k)\ge\frac{\pi f_j}{\sqrt3}(1-\xi_k)^2$ and $j-k\ge1$ we obtain
$H(k)\ge\frac{h}{2}(1-\xi_k)$, so
$\frac52\log(1/\xi_k)\le\frac{10}{h}H(k)\le\frac5{2y_0^{1/3}}H(k)$.
Moreover, $f(u)\le f_j$ on $[k,j]$ gives
$H''(u)\ge\frac{48}{\pi^2}$, and integrating twice from $j$ to $k$
yields $H(k)\ge\frac{24}{\pi^2}(j-k)^2$. Thus, setting
\[
 \tau_*=\frac{24}{\pi^2}\left(1-\frac5{2y_0^{1/3}}\right)>0,
\]
the bound \eqref{eq:outer-monomial-ratio} is at most
$1.26e^{-\tau_*(j-k)^2}$. Using $(j-k)^2\ge1+3(j-k-1)$ for
integers $k<j$, the sum over these indices $k$ is at most
\[
 \frac{1.26e^{-\tau_*}}{1-e^{-3\tau_*}}<0.477,
\]
since $\tau_*$ increases with $y_0\ge80$ and the bound is less than
$0.476245$ at $y_0=80$.

It remains to consider $k<j$ with $t_k\ge1$ and $f_k<f_j/2$, as well
as the possible case $t_k=0$. When $t_k\ge1$ and $f_k<f_j/2$, put
$\xi_k=f_k/f_j\in[\sqrt{11/12}/f_j,1/2)$.
Lemma~\ref{lem:two-colour-uniform} bounds
$\frac{p(n,m-k)e^{(m-k)h}}{p(n,d)e^{dh}}$ by
\[
 \frac{2.6}{0.875}\xi_k^{-5/2}
 \exp\!\left(-\frac{\pi f_j}{\sqrt3}(1-\xi_k)^2\right).
\]
Since $f_j>12.5$ and $\xi(1-\xi)$ increases on $[0,1/2]$, for
$\sqrt{11/12}/f_j\le\xi\le1/2$ we have
\[
 \frac{2\pi f_j}{\sqrt3}\xi(1-\xi)
 \ge\frac{2\pi}{\sqrt3}\sqrt{11/12}
       \left(1-\frac{\sqrt{11/12}}{12.5}\right)>\frac52.
\]
Thus
\[
 \frac{d}{d\xi}\log\left(
   \xi^{-5/2}e^{-\pi f_j(1-\xi)^2/\sqrt3}\right)
 =-\frac5{2\xi}+\frac{2\pi f_j}{\sqrt3}(1-\xi)>0.
\]
Evaluating at $\xi=1/2$ therefore bounds the displayed ratio by
\[
 \frac{2.6}{0.875}\,2^{5/2}e^{-\pi f_j/(4\sqrt3)}
 <16.82e^{-\pi f_j/(4\sqrt3)}.
\]
If $t_k=0$, then $k=0$, $n=m^2$, and $p(n,m)=1$. Since $d<y_0$,
direct algebra gives
\[
 h(y_0-d)-\frac{2\pi f_j}{\sqrt3}
 =-\frac{2\pi y_0}{\sqrt3(y_0+d)}f_j
 <-\frac{\pi f_j}{\sqrt3},
\]
while
$h(m-y_0)<\frac{h}{12y_0}=\frac{2\pi d}{12\sqrt3 y_0 f_j}
<\frac{\pi}{6\sqrt3 f_j}$. Therefore, in this case,
\[
 \frac{p(n,m)e^{mh}}{p(n,d)e^{dh}}
 <\frac{4\cdot3^{3/4}}{0.875}f_j^{5/2}
  e^{-\pi f_j/\sqrt3+\pi/(6\sqrt3 f_j)}
 <16.82e^{-\pi f_j/(4\sqrt3)};
\]
indeed, the ratio of the middle expression to the right-hand side
has logarithmic derivative
$\frac5{2f_j}-\frac{3\pi}{4\sqrt3}-\frac{\pi}{6\sqrt3 f_j^2}<0$
for $f_j>12.5$ and is less than $1$ at $f_j=12.5$.
Since there are at most $j<0.42y_0^{1/3}$ indices $k<j$, the sum
over all remaining cases is at most
\[
 16.82\cdot0.42y_0^{1/3}e^{-\pi\sqrt{1.96y_0}/(4\sqrt3)}<0.105,
\]
which decreases for $y_0\ge80$ (having logarithmic derivative
$\frac1{3y_0}-\frac{\pi\sqrt{1.96}}{8\sqrt3\sqrt{y_0}}<0$) and is
less than $0.104106$ at $y_0=80$.

Summing the three parts gives
$0.311+0.477+0.105<1$, proving the lemma.
\end{proof}

Proposition~\ref{thm:effective-outer-poisson} and
Lemma~\ref{lem:effective-outer-monomial} now supply an ordered
sequence of alternating signs for $D_n(-x)$ at large $x$, and we
also compare the phase at its first point with $\Psi_n$.

\begin{corollary}\label{cor:effective-outer-signs}
If $n\ge6401$ is an integer and $m=\lfloor\sqrt n\rfloor$, then
there exist an integer $d_*$ with $0.92\sqrt{n-1/12}<d_*\le m-1$
and strictly increasing numbers $x_d>0$ ($d_*\le d\le m-1$) such
that $\operatorname{sgn}D_n(-x_d)=(-1)^d$. Moreover,
\[
 \log(n-1/12)\le\log x_{d_*}\le\log(n-1/12)+2,
 \qquad
 \left|\Psi_n(\sqrt{x_{d_*}})-\pi d_*\right|<0.48,
\]
and $D_n$ has at least $m-d_*$ distinct negative zeros in
$(-\infty,-x_{d_*})$.
\end{corollary}
\begin{proof}
Put $y_0=\sqrt{n-1/12}>80$ and $h_*=\log(n-1/12)=2\log y_0>\log6400$.
For $h>0$, with $\rho(\pi)=\sqrt{(h+i\pi)^2+4\pi^2/3}$ (the
principal square root), define
\[
 \phi(h)=y_0\Ima\rho(\pi)+\arg\rho(\pi),
 \qquad
 \delta(h)=\frac{y_0 h}{\sqrt{h^2+4\pi^2/3}},
\]
and set $d_*=\lceil\phi(h_*)/\pi\rceil$.
The function $\delta$ increases strictly on $(0,\infty)$, and
$h_*+2<4y_0^{1/3}$ for $y_0>80$ (since $4y_0^{1/3}-2\log y_0-2$
increases for $y_0>4$ and is positive at $y_0=80$).
For $h>2\pi/\sqrt3$, the estimate from
\eqref{eq:outer-quadratic-derivative} proved in
Proposition~\ref{thm:effective-outer-poisson} gives
$0<\rho(0)-\Rea\rho(\pi)\le\frac{2\pi^4}{3h(h^2+4\pi^2/3)}$.
Since $\Rea\rho(\pi)>h$ and $\Ima\rho(\pi)=h\pi/\Rea\rho(\pi)<\pi$
(by \eqref{eq:outer-quadratic}), we have $0<\arg\rho(\pi)<\pi/h$
and
\[
 0<\frac{\phi(h)}{\pi}-\delta(h)
 =\frac{y_0 h\bigl(\rho(0)-\Rea\rho(\pi)\bigr)}{\rho(0)\Rea\rho(\pi)}
  +\frac{\arg\rho(\pi)}{\pi}
 \le\frac{2\pi^4 y_0}{3h(h^2+4\pi^2/3)^{3/2}}+\frac1h.
\]
Differentiating $\delta$ yields
$\delta'(h)=\frac{4\pi^2}{3}y_0(h^2+4\pi^2/3)^{-3/2}$ and
$\delta'''(h)=16\pi^2 y_0(h^2-\pi^2/3)(h^2+4\pi^2/3)^{-7/2}>0$ for
$h>\pi/\sqrt3$, so $\delta'$ is convex on $(\pi/\sqrt3,\infty)$.
Integrating $\delta'$ over $[h_*,h_*+2]$ and applying the midpoint
inequality therefore gives
\[
 \delta(h_*+2)-\frac{\phi(h_*)}{\pi}
 \ge\frac{8\pi^2 y_0}{3((h_*+1)^2+4\pi^2/3)^{3/2}}
   -\frac{2\pi^4 y_0}{3h_*(h_*^2+4\pi^2/3)^{3/2}}-\frac1{h_*}.
\]
With $t\ge\log6400$ and $e^{t/2}$ in place of $h_*$ and $y_0$, the
right-hand side equals
\[
 \frac{4\pi^2 e^{t/2}}{3t^3}
 \left\{\frac{2t^3}{((t+1)^2+4\pi^2/3)^{3/2}}
       -\frac{\pi^2 t^2}{2(t^2+4\pi^2/3)^{3/2}}\right\}
 -\frac1t.
\]
Here $e^{t/2}/t^3$ increases for $t>6$, the first term in braces
has positive derivative
$\frac{6t^2(t+1+4\pi^2/3)}{((t+1)^2+4\pi^2/3)^{5/2}}$, the
subtracted term has negative derivative
$\frac{\pi^2 t(8\pi^2/3-t^2)}{2(t^2+4\pi^2/3)^{5/2}}$ for
$t\ge\log6400$, and $-1/t$ increases. At $t=\log6400$, the
expression in braces is positive and the full right-hand side
exceeds $1.0543$; hence $\delta(h_*+2)-\phi(h_*)/\pi>1.05$, which
implies $d_*<\delta(h_*+2)$.
Furthermore,
\[
 y_0-\delta(h_*+2)
 >\frac{2\pi^2 y_0}{3((h_*+2)^2+4\pi^2/3)}>4.07,
\]
since the function
$\frac{2\pi^2e^{t/2}}{3((t+2)^2+4\pi^2/3)}$
has derivative of sign $t^2+4\pi^2/3-4>0$ and exceeds $4.07$
at $t=\log6400$.
As $m>\sqrt n-1>y_0-1$, this gives $d_*<y_0-4.07<m-1$, so
$d_*\le m-1$. Also, $d_*>\delta(h_*)>0.92y_0$ because
$t/\sqrt{t^2+4\pi^2/3}>0.92$ for $t\ge\log6400$.

For each integer $d$ with $d_*\le d\le m-1$, we have $d<y_0$, so we
may define
\[
 t_d=\frac{2\pi d}{\sqrt3\sqrt{y_0^2-d^2}},
\]
which satisfies $\delta(t_d)=d$ and increases strictly with $d$.
Since $\delta(t_d)=d\ge d_*>\delta(h_*)$, we have $t_d>h_*$.
If $t_d\le4y_0^{1/3}$, then $\phi(h_*)\le\pi d=\pi\delta(t_d)<\phi(t_d)$,
so by continuity there is a smallest $h^{(d)}\in[h_*,t_d]$ with
$\phi(h^{(d)})=\pi d$; set $x_d=e^{h^{(d)}}$.
If $t_d>4y_0^{1/3}$, set $x_d=e^{t_d}$.
To see that $x_d$ increases strictly with $d$: if
$d_*\le d_1<d_2$ and $t_{d_2}\le4y_0^{1/3}$, then
$\phi(h_*)\le\pi d_1<\pi d_2=\phi(h^{(d_2)})$, so the intermediate
value theorem yields a point in $[h_*,h^{(d_2)})$ where
$\phi=\pi d_1$, whence $h^{(d_1)}<h^{(d_2)}$; if
$t_{d_1}\le4y_0^{1/3}<t_{d_2}$, then
$h^{(d_1)}\le t_{d_1}\le4y_0^{1/3}<t_{d_2}$; and if
$4y_0^{1/3}<t_{d_1}<t_{d_2}$, then $x_{d_1}<x_{d_2}$ immediately.
For $d=d_*$, the inequality $d_*<\delta(h_*+2)$ implies
$t_{d_*}<h_*+2<4y_0^{1/3}$, so
$h_*\le\log x_{d_*}=h^{(d_*)}\le t_{d_*}<h_*+2$ and
$\phi(\log x_{d_*})=\pi d_*$.

When $t_d\le4y_0^{1/3}$, writing $h=h^{(d)}\in[h_*,4y_0^{1/3}]$ and
$r=\sqrt{x_d}=e^{h/2}$, we have $\phi(h)=\pi d$, so the
approximating term in Proposition~\ref{thm:effective-outer-poisson}
equals $(-1)^d\frac{|\rho(\pi)|}{h}\mathcal A_{y_0}(h)$, which has
sign $(-1)^d$ and modulus greater than $\mathcal A_{y_0}(h)$ (since
$|\rho(\pi)|\ge\Rea\rho(\pi)>h$). Thus
Proposition~\ref{thm:effective-outer-poisson} gives
$\operatorname{sgn}D_n(-x_d)=(-1)^d$.
When $t_d>4y_0^{1/3}$, Lemma~\ref{lem:effective-outer-monomial}
(with $h=t_d$) gives
$p(n,d)x_d^d>\sum_{\ell\ne d}p(n,\ell)x_d^\ell$, so again
$\operatorname{sgn}D_n(-x_d)=(-1)^d$.
The alternating signs at $-x_{m-1}<\cdots<-x_{d_*}$, together with
$\operatorname{sgn}D_n(x)=(-1)^m$ as $x\to-\infty$ (since
$p(n,m)>0$), yield at least $m-d_*$ distinct negative zeros of $D_n$
in $(-\infty,-x_{d_*})$.

It remains to compare $\Psi_n(\sqrt{x_{d_*}})$ with $\pi d_*$. For
any $h\ge h_*$, put $\varepsilon=e^{-h/2}\le y_0^{-1}$. Setting
$r=e^{h/2}$ in \eqref{eq:dilog-inversion} and multiplying by $4$
gives
\[
 4s(e^{h/2})^2=(h+i\pi)^2+\frac{4\pi^2}{3}+8\Li_2(i\varepsilon)
              =\rho(\pi)^2+8\Li_2(i\varepsilon).
\]
Since $\Rea\bigl(4s(e^{h/2})^2\bigr)=h^2+\frac{\pi^2}{3}+2\Li_2(-\varepsilon^2)>h^2$
and $\Rea(\rho(\pi)^2)=h^2+\frac{\pi^2}{3}>h^2$, the identity
$(\Rea w)^2=\Rea(w^2)+(\Ima w)^2$ implies that both $2s(e^{h/2})$
and $\rho(\pi)$ have real part strictly greater than $h$.
Because $|8\Li_2(i\varepsilon)|\le8\varepsilon/(1-\varepsilon)$, we
obtain
\[
 |2s(e^{h/2})-\rho(\pi)|
 =\frac{|8\Li_2(i\varepsilon)|}{|2s(e^{h/2})+\rho(\pi)|}
 \le\frac{4\varepsilon}{h(1-\varepsilon)}.
\]
Since the line segment joining $\rho(\pi)$ and $2s(e^{h/2})$ lies
in $\{\Rea w>h\}$, integrating $1/w$ along this segment yields
\[
 |\arg s(e^{h/2})-\arg\rho(\pi)|
 =\bigl|\Ima\bigl(\log(2s(e^{h/2}))-\log\rho(\pi)\bigr)\bigr|
 \le\frac{|2s(e^{h/2})-\rho(\pi)|}{h}
 \le\frac{4\varepsilon}{h^2(1-\varepsilon)}.
\]
From \eqref{eq:amplitude-phase} and
$\arctan(e^{h/2})=\frac\pi2-\arctan\varepsilon$, we have
\[
 \Psi_n(e^{h/2})
 =\sqrt n\,\Ima\bigl(2s(e^{h/2})\bigr)
  +\arg s(e^{h/2})-\frac12\arctan\varepsilon.
\]
Using $\sqrt n-y_0<\frac1{24y_0}$, $|\Ima\rho(\pi)|<\pi$,
$\arctan\varepsilon<\varepsilon\le y_0^{-1}$, and $h\ge h_*=2\log y_0$,
we find
\[
 |\Psi_n(e^{h/2})-\phi(h)|
 \le\frac{4\sqrt{y_0^2+1/12}}{y_0(2\log y_0)(1-y_0^{-1})}
   +\frac{\pi}{24y_0}
   +\frac4{y_0(2\log y_0)^2(1-y_0^{-1})}
   +\frac1{2y_0}<0.48.
\]
Each of the four terms decreases for $y_0\ge80$, and their sum at
$y_0=80$ is less than $0.470736$. Evaluating at $h=\log x_{d_*}$,
where $\phi(h)=\pi d_*$, proves
$|\Psi_n(\sqrt{x_{d_*}})-\pi d_*|<0.48$.
\end{proof}

We can now combine the sign changes from all three ranges.

\begin{proof}[Proof of Theorem~\ref{thm:main}]
Set $m=\lfloor\sqrt n\rfloor$. First suppose that $n\ge6401$, and
put $r_-=16\pi^2/n$. By Corollary~\ref{cor:effective-bessel-sign}
(which applies for $1\le k\le8$ since $n\ge6401>6400$), we have
$\operatorname{sgn}D_n\bigl(-(k\pi)^4/(16n^2)\bigr)=(-1)^k$ for
$1\le k\le8$. Since
$-r_-^2=-(8\pi)^4/(16n^2)<\cdots<-\pi^4/(16n^2)<0$, the
intermediate value theorem yields at least seven distinct negative
zeros of $D_n$ in $(-r_-^2,0)$, and $D_n(-r_-^2)>0$.

Note that $r_-=16\pi^2/n\le16\pi^2/6401<0.025$. For $0<r\le0.025$,
\eqref{eq:unit-basic-bounds} gives $r\le a(r)^2\le r(1+0.282r)$,
while \eqref{eq:parameters} gives
$0<a(r)^2-b(r)^2=-\frac12\Li_2(-r^2)\le r^2/2$. Hence
$1-r/2\le b(r)^2/r\le1+0.282r$, which implies
$|b(r)/\sqrt r-1|\le r/2$ (using $\sqrt{1+x}-1\le x/2$ and
$1-\sqrt{1-r/2}=\frac{r/2}{1+\sqrt{1-r/2}}\le r/2$).
Moreover, by \eqref{eq:parameters} and $\arctan y\ge y-y^3/3$ for
$y\ge0$, we have
$\Ima(s(r)^2)=2\int_0^r\frac{\arctan y}{y}\,dy\ge2r(1-r^2/9)>r$.
Since $s(r)^2$ lies in the open first quadrant,
\[
 0<\frac\pi4-\arg s(r)
  =\frac12\arctan\!\left(\frac{\Rea(s(r)^2)}{\Ima(s(r)^2)}\right)
  \le\frac{\Rea(s(r)^2)}{2\Ima(s(r)^2)}
  <\frac r4,
\]
while $0<\frac12\arctan r\le r/2$; consequently,
\[
 \left|\arg s(r)-\frac\pi4+\frac12\arctan r\right|<r.
\]
Evaluating \eqref{eq:amplitude-phase} at $r=r_-$ (where
$2\sqrt{nr_-}=8\pi$) therefore gives
\[
 |\Psi_n(r_-)-8\pi|
 \le8\pi\left|\frac{b(r_-)}{\sqrt{r_-}}-1\right|
   +\left|\arg s(r_-)-\frac\pi4+\frac12\arctan r_-\right|
 <(4\pi+1)r_-<0.35.
\]

For the other endpoint, let $d_*$ and $x_{d_*}$ be given by
Corollary~\ref{cor:effective-outer-signs}, and set
$r_+=\sqrt{x_{d_*}}$. Then $D_n$ has at least $m-d_*$ distinct
negative zeros in $(-\infty,-r_+^2)$, with
\[
 \operatorname{sgn}D_n(-r_+^2)=(-1)^{d_*},\qquad
 |\Psi_n(r_+)-\pi d_*|<0.48.
\]
Since $\sqrt{n-1/12}>80$, Corollary~\ref{cor:effective-outer-signs}
also gives $d_*>0.92\sqrt{n-1/12}>9$ (so $d_*\ge10$) and
$\frac12\log(n-1/12)\le\log r_+\le\frac12\log(n-1/12)+1<\frac12\log n+1$,
whence $0<r_-<1<r_+<e\sqrt n$. On $[r_-,1]$ we have
$\sqrt{nr}\ge4\pi$, so Proposition~\ref{lem:effective-unit-left}
applies, while on $[1,r_+]$ Proposition~\ref{thm:effective-right}
applies. Together, these results show that
$\operatorname{sgn}D_n(-r^2)=(-1)^k$ whenever $r\in[r_-,r_+]$ and
$\Psi_n(r)=k\pi$ ($k\in\mathbb Z$).

By the endpoint phase bounds, $\Psi_n(r_-)<8\pi+0.35<9\pi$ and
$\Psi_n(r_+)>\pi d_*-0.48>(d_*-1)\pi$. For each integer
$k\in\{9,\ldots,d_*-1\}$, continuity of $\Psi_n$ on $[r_-,r_+]$
implies that $\{r\in[r_-,r_+]:\Psi_n(r)=k\pi\}$ is nonempty and
compact; let $r_k$ be its minimum. Clearly $r_-<r_k<r_+$. If
$9\le k<k'\le d_*-1$, then $\Psi_n(r_-)<k\pi<k'\pi=\Psi_n(r_{k'})$,
so the intermediate value theorem yields a point in $(r_-,r_{k'})$
where $\Psi_n$ equals $k\pi$, whence $r_k<r_{k'}$. Thus
$-r_+^2<-r_{d_*-1}^2<\cdots<-r_9^2<-r_-^2$, and $D_n$ takes the
alternating signs $(-1)^{d_*},(-1)^{d_*-1},\ldots,(-1)^9,(-1)^8$ at
these $d_*-7$ points, yielding at least $d_*-8$ distinct negative
zeros in $(-r_+^2,-r_-^2)$. Since $0<r_-<r_+$, the three open
intervals $(-\infty,-r_+^2)$, $(-r_+^2,-r_-^2)$, and $(-r_-^2,0)$
are pairwise disjoint and together contain at least
\[
 (m-d_*)+(d_*-8)+7=m-1
\]
distinct negative zeros of $D_n$. Since $\deg D_n=m$, $D_n(0)=0$,
and $D_n'(0)=p(n,1)=n>0$, the origin and these $m-1$ negative zeros
account for all $m$ zeros of $D_n$, all of which are simple.

It remains to treat $1\le n\le6400$. Write
\[
 Q_n(x)=\frac{D_n(x)}{x}=\sum_{i=0}^{m-1}p(n,i+1)x^i.
\]
If $m=1$, then $Q_n(x)=n>0$ has no zeros.
If $m=2$, then $Q_n(x)=n+p(n,2)x$ (with $p(n,2)>0$) has one simple
negative zero. For $m\ge3$ (that is, $9\le n\le6400$), it suffices
to exhibit positive rational numbers $y_2,\ldots,y_{m-1}$ such that
\begin{equation}\label{eq:finite-signs}
 y_2<\cdots<y_{m-1},\qquad
 \operatorname{sgn}Q_n(-y_k)=(-1)^{k-1}
 \quad(2\le k\le m-1).
\end{equation}
Indeed, these signs, together with $Q_n(0)=n>0$ and
$\operatorname{sgn}Q_n(x)=(-1)^{m-1}$ as $x\to-\infty$ (since
$p(n,m)>0$), yield $m-1$ distinct negative zeros of $Q_n$, all of
which are simple because $\deg Q_n=m-1$.

We verify \eqref{eq:finite-signs} for every $9\le n\le6400$ by exact
integer arithmetic in Python. Starting from $(q;q)_0^{-2}=1$, the
series coefficients of $(q;q)_d^{-2}$ through degree $6400-d^2$ are
computed recursively from those of $(q;q)_{d-1}^{-2}$ by multiplying
twice by $(1-q^d)^{-1}=\sum_{j\ge0}q^{jd}$ using only exact integer
additions, and \eqref{eq:F} gives $p(n,d)$ for all $1\le n\le6400$
and $1\le d\le\lfloor\sqrt n\rfloor$. As an independent check, for
every $1\le n\le6400$, the identity
$\sum_{d=1}^{\lfloor\sqrt n\rfloor}p(n,d)=p(n)$ is verified against
$p(n)$ computed from Euler's pentagonal number theorem
\cite[Corollary~1.7]{Andrews}.

For each $9\le n\le6400$ (with $m=\lfloor\sqrt n\rfloor$), numerical
approximations are used solely to propose candidate rational points
$y_k=u_k/w_k$ (with positive integers $u_k,w_k$, where $w_k$ is a
power of $2$); neither the accuracy of these approximations nor any
floating-point calculation enters the verification. For a candidate
point $y=u/w$, let $\ell=m-1$ and initialize $H_\ell=p(n,m)$. The
descending recurrence
\[
 H_i=-uH_{i+1}+p(n,i+1)w^{\ell-i}
 \qquad(i=\ell-1,\ldots,0)
\]
uses only exact integer arithmetic and gives
$H_0=w^\ell Q_n(-u/w)$, so $\operatorname{sgn}H_0=\operatorname{sgn}Q_n(-u/w)$.
Strict ordering $y_k<y_{k+1}$ between consecutive points is likewise
verified by exact integer cross-multiplication $u_kw_{k+1}<u_{k+1}w_k$.

This exact integer computation confirms all
\[
 \sum_{n=1}^{6400}\max\{0,\lfloor\sqrt n\rfloor-2\}=325{,}403
\]
strict sign conditions in \eqref{eq:finite-signs}, all required
strict orderings, and all $6400$ identities
$\sum_{d=1}^{\lfloor\sqrt n\rfloor}p(n,d)=p(n)$. Thus $Q_n$ has
$\lfloor\sqrt n\rfloor-1$ simple negative zeros for every $1\le n\le6400$. Together with
the simple zero of $D_n$ at the origin, this completes the proof.
\end{proof}

\medskip
\noindent\textbf{Code and formal verification.}
The Lean~4 formalization of Theorem~\ref{thm:main} and the Python script
(\texttt{review/generate\_step05\_certificates.py}) are available at
\url{https://github.com/tfyf103/durfee-ccs-lean}.
The Python script generates and independently verifies the exact integer
certificates for $1\le n\le6400$; the Lean~4 proof checks these
certificates directly in Lean's kernel without trusting external
computation or introducing additional axioms.

\section*{Acknowledgments}
We thank the developers and contributors of Lean~4 and Mathlib for
providing the infrastructure used in the formal verification of
Theorem~\ref{thm:main}.
We used ChatGPT to assist with language
editing, exposition, literature searches, and mathematical discussion during the preparation
of this paper.  We take full responsibility for the content of this work.

\end{document}